\documentclass{gtart_a}
\pdfoutput=1
\usepackage{pinlabel}

\title{Une nouvelle preuve du th\'eor\`eme de point fixe de Handel}
\author{Patrice Le Calvez}
\givenname{Patrice}
\surname{Le Calvez}
\address{Laboratoire Analyse\\G\'eom\'etrie et
Applications\\C.N.R.S.-U.M.R 7539\\Institut 
Galil\'ee\\\newline
Universit\'e Paris 13\\Avenue
J.-B.Cl\'ement\\93430 Villetaneuse\\France}
\email{lecalvez@math.univ-paris13.fr}
\urladdr{}

\volumenumber{10}
\issuenumber{}
\publicationyear{2006}
\papernumber{52}
\lognumber{0753}
\startpage{2299}
\endpage{2349}

\doi{}
\MR{}
\Zbl{}


\received{1 June 2006}
\revised{}
\accepted{28 October 2006}
\proposed{Leonid Polterovich}
\seconded{Yasha Eliashberg, Benson Farb}
\published{8 December 2006}
\publishedonline{8 December 2006}
\corresponding{}
\editor{CPR}
\version{}

\arxivreference{}  %%% not in arXiv

\AtBeginDocument{\let\bar\wbar\let\tilde\wtilde\let\hat\what\let\widetilde\wwtilde\let\overline\wbar}
\def\co{\mskip0.7mu\colon\thinspace}

\def\bee{\begin{enumerate}}
\def\bei{\begin{itemize}}
\def\ene{\end{enumerate}}
\def\eni{\end{itemize}}

\renewcommand\theenumi{\roman{enumi}}
\def\demo{\proof[D\'emonstration]}

\makeatletter
\def\cnewtheorem#1[#2]#3{\newtheorem{#1}{#3}[section]
\expandafter\let\csname c@#1\endcsname\c@thm}

\newtheorem{thm}{Th\'eor\`eme}[section]
\cnewtheorem{prop}[thm]{Proposition}
\cnewtheorem{lem}[thm]{Lemme}

\makeatother  %  move after \newtheorem block

\def\N{{\mathbb N}}

\def\D{{\bf D}}

\begin{document}

\begin{asciiabstract}
M Handel has proved in [Topology 38 (1999) 235--264] a fixed point
theorem for an orientation preserving homeomorphism of the open unit
disk, that may be extended to the closed disk and that satisfies a
linking property of orbits. We give here a new proof of Handel's fixed
point theorem, based on Brouwer theory and some plane topology
arguments. We will slightly improve the theorem by proving the
existence of a simple closed curve of index $1$. This index result was
known to be true under an additional hypothesis and has been used by
different authors (J Franks [NYJM 2 (1996) 1--19, Trans.AMS 348 (1996)
2637--2662] S Matsumoto [Topol. Appl. 104 (2000) 191--214]) to study
homeomorphisms of surfaces.
\end{asciiabstract}

\begin{webabstract}
M Handel has proved in [Topology 38 (1999) 235--264] a fixed point
theorem for an orientation preserving homeomorphism of the open unit
disk, that may be extended to the closed disk and that satisfies a
linking property of orbits. We give here a new proof of Handel's fixed
point theorem, based on Brouwer theory and some plane topology
arguments. We will slightly improve the theorem by proving the
existence of a simple closed curve of index 1. This index result was
known to be true under an additional hypothesis and has been used by
different authors (J Franks [NYJM 2 (1996) 1--19, Trans.AMS 348 (1996)
2637--2662] S Matsumoto [Topol. Appl. 104 (2000) 191--214]) to study
homeomorphisms of surfaces.
\end{webabstract}

\begin{abstract}
M Handel has proved in \cite{Ha} a fixed point theorem for an
orientation preserving homeomorphism of the open unit disk, that may
be extended to the closed disk and that satisfies a linking property
of orbits. We give here a new proof of Handel's fixed point theorem,
based on Brouwer theory and some plane topology arguments. We will
slightly improve the theorem by proving the existence of a simple
closed curve of index $1$. This index result was known to be true
under an additional hypothesis and has been used by different authors
(J Franks \cite{Fr2,Fr3}, S Matsumoto \cite{M}) to study
homeomorphisms of surfaces.
\end{abstract}

\maketitle

\setcounter{section}{-1}
\section{Introduction}\label{sec0}

M Handel a \'enonc\'e et prouv\'e dans \cite{Ha} un
th\'eor\`eme d'existence de point fixe pour un 
hom\'eo\-mor\-phis\-me $f$ pr\'eservant
l'orientation du disque $\D=\{z\in\C\,\vert\, \vert z\vert<1\}$, qui se
prolonge au bord et qui v\'erifie une propri\'et\'e d'enlacement
d'orbites qui s'\'enonce sous diff\'erentes formes, par exemple
sous la forme suivante: il existe $n$ points $z_{i}$ dans $\D$ tels que
$$\lim_{k\to-\infty} f^k(z_{i})=\alpha_{i},\enskip \lim_{k\to+\infty}
f^k(z_{i})=\omega_{i},$$ o\`u les $2n$ points $\alpha_i$ et $\omega_i$ sont
des points distincts de $S^1=\{z\in\C\,\vert\, \vert z\vert=1\}$ et
o\`u les droites $\Delta_{i}$ passant par $\alpha_{i}$ et
$\omega_{i}$ d\'elimitent un polygone convexe compact \`a $n$
c\^ot\'es inclus dans $\D$.
Il s'agit d'un r\'esultat d'une grande importance pour l'\'etude des
hom\'eomorphismes de surfaces, qui a maintes fois \'et\'e utilis\'e
dans le domaine (voir J Franks \cite{Fr2,Fr3} ou S Matsumoto
\cite{M}). Dans la plupart des cas o\`u ce th\'eor\`eme est
utilis\'e, l'hom\'eomorphisme $f$ v\'erifie une
propri\'et\'e suppl\'ementaire qui permet d'obtenir une conclusion plus
forte: il existe une courbe ferm\'ee simple $C\subset \D$ ne
contenant pas de point fixe de $f$ et dont l'indice est
$1$. Rappelons que si
$\Gamma\co  s\mapsto\Gamma(s)$ est un param\' etrage de $C$ d\' efini
sur $S^{1}$, {\it l'indice} $i(f,C)$ est
\'egal au degr\'e de l'application
$$s\mapsto {f(\Gamma(s))-\Gamma(s)\over \vert f(\Gamma(s))-\Gamma(s)\vert}.$$
Si $i(f,C)\not=0$,
on sait alors que $f$ a un point fixe dans la composante connexe de
$\D\setminus C$ qui est simplement connexe.

\begin{figure}[ht!]\small
\labellist
\pinlabel $\alpha_1$ [tl] at 94 7
\pinlabel $\alpha_2$ [l] at 126 63
\pinlabel $\alpha_3$ [b] at 60 126
\pinlabel $\alpha_4$ [r] at 8 33
\pinlabel $\omega_1$ [b] at 91 120
\pinlabel $\omega_2$ [br] at 14 102
\pinlabel $\omega_3$ [tr] at 27 11
\pinlabel $\omega_4$ [tl] at 118 33
\endlabellist
\cl{
\psfig{file=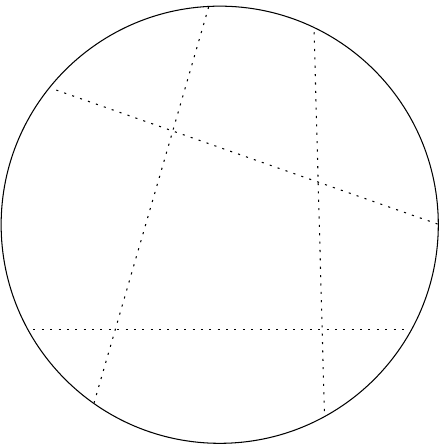}}
\caption{}\label{fig1}
\end{figure}

La preuve de Handel de ce th\' eor\`eme de point fixe est
consid\'er\'ee
comme difficile.  Elle utilise en
particulier une g\' en\' eralisation de la classification de
Nielsen-Thurston des ho\-m\'eo\-mor\-phis\-mes de surfaces, au cas
du disque $\D$ priv\' e d'une r\'eunion finie d'orbites
sans point limite. Elle utilise \' egalement la th\' eorie de Brouwer des
hom\' eomorphismes du disque  et une version homotopique
de cette th\' eorie. Nous donnons dans cet article une nouvelle preuve
du th\'eor\`eme du point fixe de Handel qui 
n'utilise que la th\' eorie de Brouwer
``classique'' et o\`u n'apparaissent donc aucun 
argument g\' eom\' etrique (existence d'une
m\'etrique de courbure
$-1$,
groupes fuchsiens,\ldots). M\^eme si la preuve n'est pas tr\`es courte,
elle n'utilise, outre une partie de la th\'eorie de Brouwer, que des
arguments simples de topologie plane. Nous prouverons que les hypoth\`eses du
th\'eor\`eme impliquent l'existence d'une {\it cha\^ine ferm\'ee  de
disques libres}, c'est-\`a-dire d'une famille $(D_i)_{i\in\Z/r\Z}$
de parties ouvertes disjointes, connexes et simplement
connexes, v\'erifiant $f(D_{i})\cap 
D_{i}=\emptyset$, telle que $\bigcup_{k\geq 
1}f^{k}(D_i)\cap D_{i+1}\not=\emptyset$, pour tout
$i\in\Z/r\Z$.
Le lemme de Franks \cite{Fr1} \'enonce qu'il 
existe alors une courbe ferm\'ee simple $C\subset 
\D$ ne
contenant pas de point fixe de $f$ et dont l'indice est
$1$. Ainsi le r\'esultat d'indice est vrai sous les hypoth\`eses du
th\'eor\`eme du point fixe de Handel, il n'y a pas besoin d'ajouter
des hypoth\`eses. Nous verrons \'egalement que l'hypoth\`ese de
prolongement au bord peut \^etre l\'eg\`erement affaiblie.

\'Enon\c cons donc le th\'eor\`eme, dont les 
hypoth\`eses sont illustr\'ees sur la \fullref{fig1}:

\begin{thm}\label{thm0.1}Soit $f$ un
hom\'eomorphisme de $\D$ qui pr\'eserve l'orientation et $n$ un
entier sup\'erieur ou \'egal \`a $3$. On suppose que les hypoth\`eses
suivantes sont satisfaites:

\bee
\item
il existe une famille $(z_i)_{i\in\Z/n\Z}$
dans $\D$ et deux familles 
$(\alpha_i)_{i\in\Z/n\Z}$ et 
$(\omega_i)_{i\in\Z/n\Z}$
dans $S^{1}$, telles que
$\lim_{k\to-\infty} f^k(z_{i})=\alpha_{i}$ et $\lim_{k\to+\infty}
f^k(z_{i})=\omega_{i}$, pour tout $i\in\ Z/n\Z$;

\item
les $2n$ points 
$\alpha_i$, $\omega_i$, $i\in\Z/n\Z$, sont
distincts;

\item
le seul, parmi ces points, qui se trouve
sur l'intervalle ouvert du cercle orient\' e 
$S^1$ joignant $\omega_{i-1}$ \`a $\omega_i$
est
$\alpha_{i+1}$;

\item$f$ se prolonge en un hom\' eomorphisme de
$\D\cup\left(\bigcup_{i\in \Z/n\Z} \{\alpha_i,\omega_i\}\right)$.

\ene

Il existe alors une courbe ferm\'ee simple $C\subset \D$ ne
contenant pas de point fixe de $f$ et dont l'indice est
$1$. 
\end{thm}

\eject
{\bf Remarques}

(1)\qua Le prolongement d\' efini par {\rm(iv)} fixe
n\'ecessairement les points
$\alpha_i$ et
$\omega_i$. Cette hypoth\`ese de prolongement est indispensable. En
effet, nous rappelerons dans la prochaine section que les orbites
d'un hom\'eomorphisme $f$ de $\D$, qui pr\'eserve l'orientation et qui
n'a pas de point fixe, n'ont pas de point limite (en d'autres termes
tendent vers le bout de $\D$ en $-\infty$ et 
$+\infty$). Ainsi, si $n\geq 3$ est fix\' e,
si $(z_i)_{i\in\Z/n\Z}$ est une famille arbitraire de points de $\D$ d'orbites
distinctes, si $(\alpha_i)_{i\in\Z/n\Z}$ et $(\omega_i)_{i\in\Z/n\Z}$
sont deux familles arbitraires dans $S^{1}$ qui v\'erifient les
hypoth\`eses {\rm(ii)} et {\rm(iii)}, il n'est pas difficile  de construire
un hom\'eomorphisme $h$ de $\D$ pr\'eservant 
l'orientation, pour que l'hypoth\`ese {\rm(i)} 
soit
satisfaite par la famille  $(h(z_i))_{i\in\Z/n\Z}$ et
l'hom\'eomorphisme
$h\circ f\circ h^{-1}$.  L'hypoth\`ese de 
prolongement qui est suppos\' ee dans l'article
original de Handel est un peu plus forte: $f$ se 
prolonge en un hom\' eomorphisme de $\overline 
{\bf
D}$.

\medskip
(2)\qua  La propri\'et\'e d'enlacement (ici
l'hypoth\`ese {\rm(iii)}) peut \^etre
formul\'ee dif\-f\'e\-re\-ment. Ainsi, dans l'article de Handel, les
$z_i$ sont index\' es par $\{1,\dots,n\}$, et
l'hypoth\`ese {\rm(iii)} est remplac\' ee par les deux assertions suivantes:

{\sl {\rm(iii$'$)}\qua tout 
intervalle de $S^1$ d\' elimit\' e par deux points
$\alpha$-limites (resp.\ $\omega$-limites) 
contient un point $\omega$-limite (resp.\ 
$\alpha$-limite);

{\rm(iii$''$)}\qua pour tout 
$i\in\{1,\dots,n\}$, il existe 
$i'\in\{1,\dots,n\}$ tel que
$\alpha_i$ et $\omega_i$ sont s\' epar\' es sur 
$S^{1}$ par $\alpha_{i'}$ et $\omega_{i'}$.}

\medskip
Expliquons rapidement pourquoi ces deux conditions impliquent que
la propri\'et\'e {\rm(iii)} est vraie (\`a orientation pr\`es) si on ne 
garde qu'un certain nombre des $z_i$.
Notons
$\Delta_i$ la droite affine r\'eelle orient\' ee, 
engendr\' ee par le vecteur joignant
$\alpha_i$ \`a $\omega_i$, et
$\Pi_i^+$ (resp.\
$\Pi_i^-$) le demi-plan affine situ\' e 
strictement \`a gauche (resp.\ \`a droite) de
$\Delta_i$. Quitte \`a perturber l\' eg\`erement les points
$\alpha_i$ et
$\omega_i$, on peut supposer que l'intersection de trois droites $\Delta_i$ est
toujours vide. Il suffit de
trouver une partie
$I$ de $\{1,\dots,n\}$ telle que l'intersection
$\bigcap_{i\in I}
\Pi_i^+$ ou $\bigcap_{i\in I}
\Pi_i^-$ soit relativement compacte dans $\D$, 
puis de choisir $I'\subset I$, minimale parmi les
parties de $I$ satisfaisant cette propri\' et\' e.  L'assertion {\rm
(iii$'$)} implique que la fonction 
$\nu\co \C\setminus\bigcup_{1\leq i\leq 
n}\Delta_i\to\N$,  qui
\`a $z$ associe le nombre de demi-plans $\Pi_i^+$ 
contenant $z$, prend exactement deux
valeurs au voisinage de
$S^1$ et que ces deux valeurs sont des entiers 
cons\' ecutifs. L'assertion {\rm(iii$''$)} implique 
qu'il
existe au moins deux droites $\Delta_i$ et
$\Delta_{i'}$ qui s'intersectent dans $\D$. Au 
voisinage du point d'intersection, la fonction 
$\nu$
prend trois valeurs cons\' ecutives. On en d\' 
eduit que le maximum ou le minimum de la fonction 
$\nu$ sur
$\D\setminus\bigcup_{1\leq i\leq n}\Delta_i$ 
n'est pas l'une des deux valeurs pr\`es du bord. 
Dans
le premier cas (resp.\ second cas), on consid\' 
ere une composante connexe $U$ de
$\D\setminus\bigcup_{1\leq i\leq n}\Delta_i$ o\`u 
ce maximum (resp.\ minimum) est atteint. Elle
est relativement compacte dans $\D$ et on peut \' ecrire
$U=\bigcap_{i\in I} \Pi_i^+$ (resp.\ $U=\bigcap_{i\in I} \Pi_i^-$), o\`u
$I\subset\{1,\dots,n\}$.

Nous allons continuer cette introduction en rappelant certaines applications
importantes du th\'eor\`eme de Handel \`a l'\'etude des
hom\'eomorphismes de surfaces isotopes \`a l'identit\'e. Si
$(F_t)_{t\in[0,1]}$ est une isotopie issue de 
l'identit\' e sur une surface $M$, on peut
d\' efinir la {\it trajectoire} $\gamma_z 
\co t\mapsto F_t(z)$ de tout point $z\in M$. Si $z$ 
est un
point p\' eriodique, de p\' eriode $q\geq 1$, on 
obtient alors un lacet $\prod_{0\leq i\leq
q-1}\gamma_{F^i(z)}$ par assemblage des 
trajectoires, dont on note $\kappa(z_i)\in 
H_1(M,\Z)$ la classe
d'homologie. Les deux r\' esultats qui suivent
sont des cons\' equences directes du th\' eor\`eme de point fixe.
Le premier est attribu\' e par Franks dans
\cite{Fr3} \`a M Betsvina et Handel, le second 
est d\^u \`a Franks \cite{Fr2}. Nous rappelerons
\' egalement les grandes lignes des preuves que l'on trouve dans
\cite{Fr2} et \cite{Fr3}.

\begin{thm}\label{thm0.2}Soit $M$ une partie 
ouverte de la
sph\`ere $S^2$, dont le com\-pl\' e\-men\-tai\-re a au moins trois composantes
connexes,
et $\widetilde M$ son rev\^etement universel.
Soit
$(F_t)_{t\in[0,1]}$ une isotopie sur $M$ issue de l'identit\' e, et
$(f_t)_{t\in[0,1]}$ l'isotopie relev\'ee \`a 
$\widetilde M$ qui est issue de l'identit\' e.
S'il existe un point p\' eriodique $z_0$, de p\' eriode $q\geq 1$, tel que
$\kappa(z_0)=0$, alors on peut trouver une courbe ferm\'ee simple $C\subset
\widetilde M$ ne
contenant pas de point fixe de $f_{1}$ et dont l'indice $i(f_{1},C)$ est
$1$. En particulier, $F_1$ a au moins un point fixe.
\end{thm}

\proof[Id\'ee de la preuve] On 
peut munir $M$ d'une structure riemannienne 
compl\`ete de
courbure $-1$ et identifier le rev\^etement 
universel $\widetilde M$ au disque de Poincar\' e 
$\D$. On
sait alors que $f_1$ se prolonge en un hom\' 
eomorphisme de $\overline{\bf D}$ qui fixe tout 
point de
$S^1$ (voir Handel et Thurston \cite{HaT}).  Notons 
$\Gamma$ l'unique g\' eod\' esique
ferm\' ee qui est librement homotope \`a $\prod_{0\leq i\leq
q-1}\gamma_{F_{1}^i(z_0)}$. Tout rel\`evement de $\Gamma$ dans $\bf D$ est
une g\' eod\' esique de
$\D$ qui joint un point $\alpha\in S^1$ \`a un point $\omega\in S^1$.
De plus, il existe un
ant\' ec\' edent $\wtilde z_0$ de $z_0$ tel que $\lim_{k\to-\infty}
f_{1}^k(z_{0})=\alpha$ et $\lim_{k\to+\infty}
f_{1}^k(z_{0})=\omega$. Puisque
$\Gamma$ est homologue \`a z\' ero dans $M$, il existe une fonction
$\Lambda\co M\setminus\Gamma\to\Z$, \`a support relativement compact, tel que
$\Lambda(z')-\Lambda(z)$ repr\' esente le nombre 
d'intersection alg\' ebrique $\Gamma\wedge\Gamma'$
entre
$\Gamma$ et un arc quelconque $\Gamma'$ qui joint 
$z$ \`a $z'$. L'un au moins des deux nombres
$\max\Lambda$ ou
$\min\Lambda$ doit \^etre non nul. Supposons par 
exemple que $\max\Lambda\not=0$ et fixons une
composante connexe $U$ de
$M\setminus\Gamma$ o\`u $\Lambda$ atteint son 
maximum. Il n'est pas difficile de voir que $U$ 
est
un disque ouvert, que sa fronti\`ere est une 
courbe ferm\' ee simple, r\' eunion de $p$ 
sous-segments de
$\Gamma$, et que $U$ est localement \`a gauche de 
chacun de ces segments. Toute composante connexe
de la pr\' eimage de $U$ dans $\D$ est 
l'intersection (relativement compacte) de $p$
demi-plans
g\' eod\' esiques, \`a gauche de g\' eod\' 
esiques relevant $\Gamma$. Parmi ces demi-plans, 
on
choisit une  famille dont l'intersection est 
relativement compacte dans $\D$, et qui est 
minimale
pour cette propri\' et\' e. Les conditions du 
\fullref{thm0.1} sont alors v\' erifi\' 
ees.\endproof

\begin{thm}\label{thm0.3}
Soit $M$ une surface compacte orientable de
genre $g\geq 1$ et $\widetilde M$ son rev\^etement universel.
Soit
$(F_t)_{t\in[0,1]}$ une isotopie sur $M$ issue de l'identit\' e et
$(f_t)_{t\in[0,1]}$ l'isotopie relev\'ee \`a 
$\widetilde M$ qui est issue de l'identit\' e.
On suppose qu'il existe $2g+1$ points p\' 
eriodiques $z_i$, $0\leq i\leq 2g$, tels que
l'int\' erieur de l'enveloppe convexe des $\kappa (z_i)\in H_1(M,\R)$
contient $0$. Il
existe alors une courbe ferm\'ee simple $C\subset
\widetilde M$ ne
contenant pas de point fixe de $f_{1}$ et dont l'indice $i(f_{1},C)$ est
$1$. En particulier, $F_1$ a au moins deux points fixes.
\end{thm}

\proof[Id\' ee de la preuve] Dans 
le cas o\`u $g\geq 2$, la preuve est tr\`es 
proche de
celle du \fullref{thm0.2}. On munit
$M$ d'une structure riemannienne de courbure $-1$ et on identifie
$\widetilde M$ \`a $\D$. On sait, l\`a encore, 
que $f_1$ se prolonge en un hom\' eomorphisme
de $\overline{\bf D}$ qui fixe tout point de
$S^1$. On note $\Gamma_i$ la g\' eod\' esique
ferm\' ee homotope \`a $\prod_{0\leq i\leq
q_i-1}\gamma_{F_{1}^i(z_0)}$, o\`u $q_i$ est la 
p\' eriode de $z_i$. On peut alors trouver des 
entiers
$m_i>0$ tels que $\sum_{i=0}^{2g} m_i\kappa(z_i)=0$. La multi-courbe
$\Gamma=\sum_{i=0}^{2g}m_i\Gamma_i$ est homologue 
\`a z\' ero et on lui associe une fonction duale
$\Lambda\co M\setminus\Gamma\to\Z$, comme dans la 
preuve du \fullref{thm0.2}. L\`a encore, on peut
supposer que $\max\Lambda\not=0$. Pour montrer le r\'esultat d'indice,
il suffit de prouver que toute composante connexe $U$ de
$M\setminus\Gamma$ o\`u $\Lambda$ atteint son maximum est un disque
ouvert, puis de reprendre les
arguments donn\' es dans la preuve pr\' ec\' edente. Puisque les classes
$[\Gamma_i]$ engendrent $H_1(M,\R)$, 
l'application $\iota_* \co  H_1(U,\R)\to H_1(M,\R)$
induite par l'inclusion $\iota \co U\to M$ doit 
\^etre nulle. On en d\' eduit que l'intersection 
alg\' ebrique
entre deux lacets quelconques de $U$ est nul et donc que $U$ est de genre
z\'ero. Il reste \`a montrer
que son compl\' ementaire est connexe. Dans le cas contraire, il
existerait une partition non triviale $\{0,\dots ,2g\}=I_1\sqcup I_2$, telle
que $\sum_{i\in I_1} m_i[\Gamma_i]=\sum_{i\in I_2} m_i[\Gamma_i]=0$, ce qui
contredirait le fait que
$0$ est \`a l'int\' erieur de l'enveloppe convexe des $[\Gamma_i]$.

Si $g=1$, la preuve est plus simple. On peut \' 
ecrire $M=\C/(\Z+i\Z)$ et on sait alors que $f_1$
commute avec les translations enti\`eres. On en d\' eduit que
  $f_1-{\rm Id}_{\C}$ est une application born\' 
ee et donc que $f_1$ se prolonge en un hom\' 
eomorphisme
sur le compactifi\' e de
$\C$ obtenu en rajoutant les directions \`a l'infini. En d'autres
termes, l'hom\' eomorphisme $\wtilde f_{1}=h\circ
f_{1}\circ h^{-1}$ de $\D$, conjugu\' e par
$\displaystyle h\co \C\to\D, z\mapsto {z\over 1+\vert z\vert}$, se prolonge en
l'identit\' e sur $S^1$. Remarquons maintenant 
que $\wtilde f_{1}$ v\' erifie les 
hypoth\`eses du \fullref{thm0.1}. avec $n=3$.

 Dans chacun des cas, on sait que $f_{1}$ 
(et donc que $F_{1}$) a au moins un point fixe.
Si $F_{1}$ n'avait qu'un seul point fixe, la formule de Lefschetz nous
dirait que son
indice de Lefschetz est $2-2g$. Tout point fixe de $f_{1}$ serait
donc d'indice $2-2g$, ce qui impliquerait que l'indice de toute
courbe ferm\'ee simple $C$ ne
contenant pas de point fixe de $f_{1}$ serait \'egal \`a $m(2-2g)\leq
0$, l'entier $m$ \'etant le nombre de points fixes dans la composante connexe
simplement connexe de
$\widetilde M\setminus C$. Le m\^eme raisonnement, en utilisant cette
fois-ci la
formule de Lefschetz-Nielsen, nous dit en fait que $F_{1}$ a au moins deux
points fixes qui se rel\`event en des points fixes de $f_{1}$. \endproof

En utilisant alors un simple argument de perturbation, Franks donne
dans \cite{Fr2} la version plus g\' en\' erale
suivante, o\`u l'ensemble de rotation d\' esigne 
la partie convexe et compacte form\' ee des 
vecteurs de
rotation (ou enlacements asymptotiques) de toutes 
les mesures bor\' eliennes de probabilit\' e
invariantes par
$F_1$ (voir S Schwartzman \cite{Sc}).

\begin{thm}\label{thm0.4}
Soit $M$ une surface compacte orientable de
genre $g\geq 1$ et $\widetilde M$ son rev\^etement universel.
Soit
$(F_t)_{t\in[0,1]}$ une isotopie sur $M$ issue de l'identit\' e et
$(f_t)_{t\in[0,1]}$ l'isotopie relev\'ee \`a 
$\widetilde M$ qui est issue de l'identit\' e.
Si l'ensemble de rotation $\rho(f_1)\subset H_1(M,\R)$ d\' efini par l'isotopie
contient
$0$ dans son int\' erieur, il
existe alors une courbe ferm\'ee simple $C\subset
\widetilde M$ ne
contenant pas de point fixe de $f_{1}$ et dont l'indice $i(f_{1},C)$ est
$1$.\end{thm}

Les r\'esultats pr\'ec\'edents sont d'une grande 
utilit\'e pour montrer l'existence de points
fixes ou de points p\' eriodiques. Int\' eressons 
nous par exemple au cas o\`u $M$ est une surface
compacte munie d'une forme symplectique et 
l'isotopie $(F_t)_{t\in[0,1]}$ est hamiltonienne
(i.e.\ induite par un champ de vecteurs 
hamiltonien de classe $C^1$ d\' ependant du 
temps). Le vecteur
de rotation de la mesure naturellement d\' efinie 
par la forme est alors \' egal \`a $0$. Un r\'
esultat, d\^u
ind\' ependamment \`a A Floer
\cite{Flo} et \`a J\,C Sikorav \cite{Si}, nous dit que $F_1$ a au moins trois points
fixes $z$ contractiles (i.e.\ tels que le lacet 
$\gamma_z$ soit homotope \`a z\' ero). Dans 
\cite{Fr2},
Franks donne une preuve de ce r\' esultat qui se d\' eduit du
\fullref{thm0.4}, toujours par des arguments de 
perturbation. La m\' ethode de Franks a \' et\' e
g\' en\' eralis\' ee au cas des hom\' eomorphismes. Matsumoto \cite{M}
utilise \' egalement le
th\' eor\`eme de point fixe de Handel pour donner 
la version topologique suivante de la conjecture
d'Arnold pour les hom\' eomorphismes de surfaces:

\begin{thm}\label{thm0.5}Soit
$M$ une surface compacte orientable de genre 
$g\geq 1$ et $(F_t)_{t\in[0,1]}$ une isotopie 
issue
de l'identit\' e. On suppose que $F_1$ pr\' 
eserve la mesure de Lebesgue $\mu$ et que le
vecteur de rotation de cette mesure est nul. Il 
existe alors au moins trois points fixes
contractiles.
\end{thm}

Nous rappelerons dans la section qui suit certains points de la th\' eorie des
hom\' eo\-mor\-phis\-mes du plan initi\' ee par Brouwer. En
particulier, nous donnerons des crit\`eres
d'existence de courbes ferm\' ees simple
d'indice
$1$ pour un hom\'eo\-mor\-phis\-me de $\D$ qui pr\'eserve
l'orientation. Nous rappelerons dans la \fullref{sec2} la notion de
d\'ecomposition en briques libre, qui sera l'outil principal pour
montrer que ces crit\`eres sont v\'erifi\'es. Nous
donnerons \`a la fin de la \fullref{sec2} le plan de la preuve du
\fullref{thm0.1}. Nous
d\'etaillerons cette preuve \`a partir de la \fullref{sec3}. Nous allons
conclure
cette introduction en introduisant les principales d\' efinitions
de cet article.

On \' ecrira respectivement $\overline X$, ${\rm 
Int} (X)$, $\partial X$ pour l'adh\' erence, 
l'int\' erieur
et la fronti\`ere d'une partie $X$ d'un espace 
topologique $E$. On dira qu'une partie $Y\subset 
E$
{\it s\' epare} deux parties
connexes $X_1$ et $X_2$ si $Y$ est disjoint de 
$X_1$ et $X_2$ et si $X_1$ et $X_2$ sont dans des 
composantes connexes
diff\' erentes de $E\setminus Y$.

On appelera {\it disque ouvert} (resp.\ {\it disque ferm\'e}) d'une
surface
$M$ toute partie hom\'eo\-morphe \`a
${\bf D}$
(resp.\ \`a $\overline \D$).

On appelera {\it arc} une application continue $\gamma$ d'un
intervalle $I\subset\R$ vers
une surface
$M$, ou plus pr\'ecis\'ement, une classe d'\'equivalence d'applications
par reparam\'etrage pr\'eservant l'orientation. On utilisera souvent,
par
abus de langage, le mot arc pour d\'esigner l'image de l'arc,
en faisant bien attention qu'il n'y ait pas d'ambiguit\'e.

Si $\gamma_1\co I_1\to M$ et
$\gamma_2\co I_2\to M$ sont deux arcs
tels que $$\sup I_1\in I_1,\enskip\enskip\inf I_2\in
I_2,\enskip\enskip \gamma_1(\sup
I_1)=\gamma_2(\inf I_2),$$ on
notera $\gamma_1\gamma_2$ l'arc naturellement d\'efini par assemblage
de
$\gamma_1$ et $\gamma_2$. Si $J$ est un intervalle de $\Z$ et si
$(\gamma_j)_{j\in J}$ est une
famille d'arcs $\gamma_j\co I_j\to M$ telle que
$$\sup I_j\in I_j,\enskip\enskip\inf
I_{j+1}\in I_{j+1},\enskip\enskip\gamma_j(\sup
I_j)=\gamma_{j+1}(\inf I_{j+1}),$$ si $j$ et $j+1$
appartiennent \`a $J$, on \'ecrira
$\displaystyle\prod_{j\in J}\gamma_j$ l'arc d\'efini par assemblage
des
$\gamma_j$.

On appelera {\it droite} (resp.\ {\it demi-droite n\'egative}, {\it
demi-droite positive}, {\it segment}, {\it
cercle}) d'une surface $M$
tout plongement topologique propre de
$\R$ (resp.\ $]-\infty,0]$,
$[0,+\infty[$, $[0,1]$, $S^1$), plus pr\'ecis\'ement, toute classe
d'\'equivalence d'un tel plongement
par reparam\'etrage pr\'eservant l'orientation. Un tel objet est alors
d\'etermin\'e par son image et une
orientation. On utilisera, l\`a encore, le mot droite,
demi-droite,
segment ou cercle, pour
d\'esigner l'image du plongement. Si $\gamma$ est 
un segment ou une demi-droite, on notera ${\rm
int}(\gamma)$ l'arc priv\' e des ses extr\' emit\' es.

On notera ${\rm
Homeo}^+(\D)$ l'ensemble des hom\' eomorphismes du disque
$\D$
  qui pr\' eservent l'orien\-ta\-tion. Pour tout $f\in{\rm
Homeo}^+(\D)$ on notera ${\rm Fix}(f)$ l'ensemble des points fixes de $f$.
On dira qu'une partie $X\subset\D$ est {\it libre} si $f(X)\cap
X=\emptyset$.

On identifiera $\R^{2}$ \`a $\C$, via l'isomorphisme naturel
$(x,y)\mapsto x+iy$.

\section{Rappels sur la th\' eorie de Brouwer}\label{sec1}

Commen\c cons par rappeler le {\it lemme de translation de
Brouwer} qui est le point fondamental de la th\' 
eorie (voir \cite{Br}, \cite{Bn}, \cite{Fa} ou
\cite{G}):

\begin{prop}\label{prop1.1}Soit
$f\in{\rm Homeo}^+(\D)$. On peut trouver un cercle d'indice 1
dans chacun des deux cas suivants:

\bei
\item
il existe un segment
$\gamma$ qui joint un point $z\in{\rm 
Fix}(f^2)\setminus {\rm Fix}(f) $ \`a $f(z)$ et 
tel que
$f(\gamma)\cap \gamma=\{z,f(z)\}$;

\item
il existe un segment
$\gamma$ qui joint un point $z\not\in{\rm Fix}(f^2)$ \`a $f(z)$, tel que
$f(\gamma)\cap \gamma=\{f(z)\}$, et un entier 
$k\geq 2$ tel que $f^k(\gamma)\cap
\gamma\not=\emptyset$.
\eni
\end{prop}

\begin{figure}[ht!]\small
\vspace{2mm}\labellist
\pinlabel $\gamma$  [b] at 59 57
\pinlabel $\gamma$  [b] at 280 57
\pinlabel $f(z)$  [l] at 119 42
\pinlabel $f(z)$  [l] at 342 41
\pinlabel $z$  [r] at 218 41
\pinlabel $z$  [r] at 1 42
\pinlabel $f(\gamma)$  [t] at 60 31
\pinlabel $f(\gamma)$  [tl] at 300 35
\pinlabel $f^2(z)$  [r] at 296 0
\endlabellist
\cl{\psfig{file=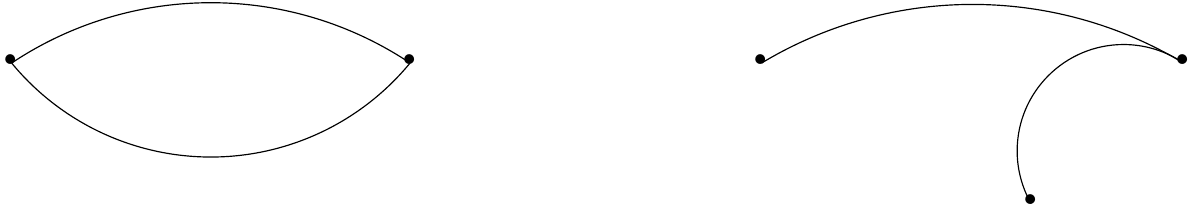}}
\vspace{-2mm}\caption{}\label{fig2}
\end{figure}

  Un segment
$\gamma$ qui joint un point $z\not\in{\rm Fix}(f)$ \`a $f(z)$ et tel que
$$\begin{cases}
  f(\gamma)\cap \gamma=\{z,f(z)\} & \text{si $f^2(z)=z$}\\
f(\gamma)\cap \gamma=\{z\} & \text{si
$f^2(z)\not=z$,}
\end{cases}
$$
s'appelle un {\it
arc de translation} de
$f$ (voir \fullref{fig2}).

Il n'est pas difficile de montrer que tout point 
qui n'est pas fixe appartient \`a un
arc de translation (nous en rappelerons 
d'ailleurs la preuve en \fullref{sec3}). On en d\' 
eduit
l'\' enonc\' e suivant:

\begin{prop}\label{prop1.2}Soit
$f\in{\rm Homeo}^+(\D)$. S'il existe
un point p\' eriodique qui n'est pas fixe, alors 
il existe un cercle d'indice $1$.
\end{prop}

\demo
Supposons que $z'$ soit un point p\' eriodique de
p\' eriode
$q\geq 2$ et consid\' erons un arc de translation 
$\gamma$ joignant un point $z$ \`a $f(z)$ et
contenant $z'$. Remarquons maintenant que l'une 
des deux hypoth\`eses de la \fullref{prop1.1}
est n\' ecessairement v\' erifi\' ee.
\endproof

\'Enon\c cons maintenant
le {\it lemme de Franks}, introduit dans 
\cite{Fr1}, qui \' etend le crit\`ere pr\' ec\' 
edent. Pour cela,
rappelons qu'une {\it cha\^\i ne ferm\' ee} de 
$f\in{\rm Homeo}^+(\D)$ est une famille 
$(X_i)_{i\in\Z/r\Z}$
de parties de $\D$ telle que, pour tout
$i\in\Z/r\Z$,  $\bigcup_{k\geq 1}f^{k}(X_i)\cap 
X_{i+1}\not=\emptyset$. L'entier
$r$ est alors la {\it longueur} de la cha\^\i ne.

\begin{prop}\label{prop1.3}Soit  $f\in{\rm
Homeo}^+(\D)$. S'il existe une
cha\^\i ne ferm\' ee de disques ouverts libres 
disjoints deux \`a deux, alors il existe
un cercle d'indice 1.\end{prop}

\demo
Supposons qu'il existe une cha\^\i ne ferm\' ee 
de disques ouverts libres disjoints
deux \`a deux. Parmi ces cha\^\i nes, choisissons en une de longueur
minimale, not\'ee $(D_i)_{i\in\Z/r\Z}$. Pour tout
$i\in\Z/r\Z$, choisissons un point $z_i\in D_i$ 
dont l'orbite positive rencontre
$D_{i+1}$,
en notant $f^{k_i}(z_i)$ le premier point o\`u 
cette orbite rencontre $D_{i+1}$. Puisque $r$ est
minimal, on peut remarquer que tout point 
$f^k(z_i)$, $0<k<k_i$, est en dehors de chaque 
disque de
la cha\^\i ne. On peut construire une isotopie $(h_t)_{t\in[0,1]}$
dans ${\rm Homeo}^+(\D)$, issue de l'identit\' e, telle que le
support de chaque
$h_t$ est inclus dans
$\bigcup_{i\in\Z/r\Z}D_i$, et telle que $h_1$ 
envoie $f^{k_{i}}(z_{i})$ sur $z_{i+1}$. 
Remarquons
alors que $h_1\circ f$ a un point p\' eriodique 
de p\' eriode $\sum_{i\in\Z/r\Z}k_i\geq 2$.
D'apr\`es la \fullref{prop1.2}, il
existe un cercle
$C\subset\D\setminus{\rm Fix}( h_1\circ f)$ tel que $i(h_1\circ f,C)=1$.
Remarquons, puisque les disques $D_i$
sont libres, que ${\rm Fix}( h_t\circ f)={\rm 
Fix}( f)$, pour tout $t\in[0,1]$. On en d\' eduit 
que la
fonction $t\mapsto i(h_t\circ f,C)$ est bien d\' efinie. Elle est
continue et \`a valeurs enti\`eres, elle est donc 
constante \' egale \`a $1$. \endproof

 Nous dirons qu'un hom\' eomorphisme 
$f\in{\rm Homeo}^+(\D)$ est {\it r\' ecurrent} 
s'il
v\' erifie les hypoth\`eses du lemme de Franks, 
c'est-\`a-dire s'il existe une cha\^\i ne ferm\' 
ee de disques
ouverts libres disjoints deux \`a deux. Nous 
pouvons alors \' enoncer deux classiques
crit\`eres de
r\' ecurrence:

\begin{prop}\label{prop1.4}Un hom\' eomorphisme $f\in{\rm
Homeo}^+(\D)$ est r\' ecurrent s'il v\' erifie l'une des conditions suivantes:

\bee
\item
il existe une cha\^\i ne
ferm\' ee de parties libres, connexes par arcs et disjointes deux \`a deux;

\item
il existe un arc de translation
$\gamma$ qui joint un point $z\not\in{\rm 
Fix}(f^2)$ \`a $f(z)$ et un entier $k\geq 2$ tel 
que
$f^k(\gamma)\cap
\gamma\not=\emptyset$.\ene\end{prop}

\demo
Supposons que la condition {\rm(i)}
est v\' erifi\' ee.
Parmi
les cha\^\i nes ferm\' ees de parties libres 
connexes par arcs et disjointes deux \`a deux, 
choisissons en une
de longueur minimale, not\'ee $(X_i)_{i\in\Z/r\Z}$. Pour tout
$i\in\Z/r\Z$, choisissons un point $z_i\in X_i$ 
dont l'orbite positive rencontre $X_{i+1}$
et notons $f^{k_i}(z_i)$ le premier point o\`u 
cette orbite rencontre $X_{i+1}$.
D'apr\`es la minimalit\'e de $r$, deux cas sont possibles:

\bei
\item $r=1$ et
$f^{k_1}(z_1)=z_1$;
\item
$f^{k_i}(z_i)\not=z_{i+1}$, pour tout $i\in\Z/r\Z$.
\eni

Dans le premier cas, $z_1$ est un point p\' 
eriodique de p\' eriode $k_1>1$. Il suffit alors 
de choisir
un disque ouvert contenant
$z_1$ et libre pour obtenir une cha\^\i ne ferm\' ee de longueur $1$:
l'existence d'un point p\' eriodique qui
n'est pas fixe implique la r\' ecurrence. Dans le second cas, puisque chaque
$X_i$ est connexe par arcs, on peut trouver un
segment
$\gamma_i\subset X_i$ joignant $f^{k_i}(z_i)$ \`a 
$z_{i+1}$. Si on remplace chaque $\gamma_i$ par
un disque ouvert $D_i$, voisinage assez proche de 
$\gamma_{i}$, on obtient notre cha\^\i ne.

Supposons maintenant que la condition {\rm(ii)} 
est v\' erifi\' ee. On peut toujours supposer que
$f^{k'}(\gamma)\cap
\gamma=\emptyset$ si $2\leq k'<k$ et donc que 
$f^k(\gamma)$ ne contient pas $f(z)$.
Deux cas sont alors possibles:

\bei
\item
$f^k(\gamma)\cap\gamma$ 
se r\' eduit \`a un seul point $z'$ et ce point 
est
$f^{k+1}(z)$;

\item
il existe un point 
$z''\not=f(z)$ dont l'image par $f^k$ appartient 
\`a
$\gamma$.
\eni

Pla\c cons nous dans le premier cas. Si $z'=z$, alors
$f$ a un point p\' eriodique de p\' eriode
$k+1\geq3$ et $f$ est donc r\' ecurrent. Si 
$z'\not=z$, le sous-segment de $\gamma$ qui joint 
$z$ \`a
$z'$ est libre et d\' efinit une cha\^\i ne 
ferm\' ee de longueur $1$, la condition {\rm(i)} 
est v\' erifi\' ee.
Pla\c cons nous maintenant dans le second cas. 
Choisissons $z'\not=f(z)$ sur $\gamma$ assez 
proche de $f(z)$ pour que
le sous-segment de $\gamma$ qui joint $z'$ \`a
$f(z)$ ne contienne pas $z''$ et soit disjoint de 
$f^{k}(\gamma)$. Le sous-segment de $\gamma$ qui
joint $z$ \`a
$z'$ est libre et d\' efinit une cha\^\i ne 
ferm\' ee de longueur $1$, la condition {\rm(i)} 
est
encore v\' erifi\' ee.\endproof

 Notre propri\' et\' e de r\' ecurrence, 
qui est clairement invariante par conjugaison, 
est en fait
tr\`es faible. Consid\'erons, dans le plan complexe $\C$,
l'hom\'eomorphisme $g\co  z\mapsto-{1\over 2}z$. Le 
demi-plan de Poincar\' e $\Im (z)>0$ est un 
domaine p\' eriodique
de $g$, de
p\' eriode $2$. Ainsi, tout
hom\' eomorphisme
$f\in{\rm Homeo}^+(\D)$ qui est conjugu\' e \`a $g$ est r\'ecurrent.
Un tel hom\' eomorphisme a une dynamique Nord-Sud 
entre le bout de $\D$ et un point
fixe.

\'Enon\c cons pour finir le crit\`ere suivant, 
d\^u \`a L Guillou et F Le Roux (voir
\cite[page 39]{LeR}), qui est
un peu plus d\' elicat \`a montrer:

\begin{prop}\label{prop1.5}Soit
$f\in{\rm Homeo}^+ (\D)$. Supposons qu'il existe 
une cha\^\i ne ferm\' ee $(X_i)_{i\in\Z/r\Z}$
de parties libres dont les
int\' erieurs sont disjoints deux \`a deux et qui 
v\' erifient la propri\' et\' e
sui\-van\-te:  si $z$ et $z'$ sont deux points
distincts de
$X_i$, il existe un segment joignant $z$ \`a $z'$ 
tel que ${\rm int}(\gamma)\subset {\rm Int}(X_i)$.
Alors $f$ est r\' ecurrent.\end{prop}

Pour d\'emontrer le \fullref{thm0.1}, nous prouverons le r\'esultat
plus fort suivant:

\begin{thm}\label{thm1.6}Si $f\in{\rm Homeo}^+(\D)$ v\' erifie les
hypoth\`eses {\rm(i)}, {\rm(ii)}, {\rm(iii)}, 
{\rm(iv)} du \fullref{thm0.1}, alors $f$ est r\' 
ecurrent.\end{thm}

Le principe de la
preuve est le suivant. Nous \'enoncerons un certain nombre
d'hypoth\`eses qui, si elles ne sont pas
satisfaites par $f$,  obligent $f$ \`a v\' erifier l'une des deux
assertions {\rm(i)} ou {\rm(ii)} de la \fullref{prop1.4} (et
donc \`a \^etre r\' ecurrent). Dans le cas o\`u 
toutes ces conditions sont satisfaites, nous 
montrerons que
$f$ v\' erifie les hypoth\`eses de la \fullref{prop1.5}, et est donc \'
egalement r\' ecurrent. Nous construirons
pour cela une {\it d\' ecomposition en briques 
libre} particuli\`ere de ${\bf D}\setminus{\rm 
Fix}(f)$.
Ces objets ont \' et\' e introduits par M Flucher
\cite{Flu} puis d\'evelopp\'es par A Sauzet
\cite{Sa} (voir \' egalement \cite
{LeR} et \cite{LeC}). Nous allons rappeler dans la \fullref{sec2} la
d\'efinition ainsi que les principales propri\'et\'es des
d\'ecompositions en briques.

\section{D\'ecomposition en briques}\label{sec2}

\begin{figure}[ht!]
\cl{\psfig{file=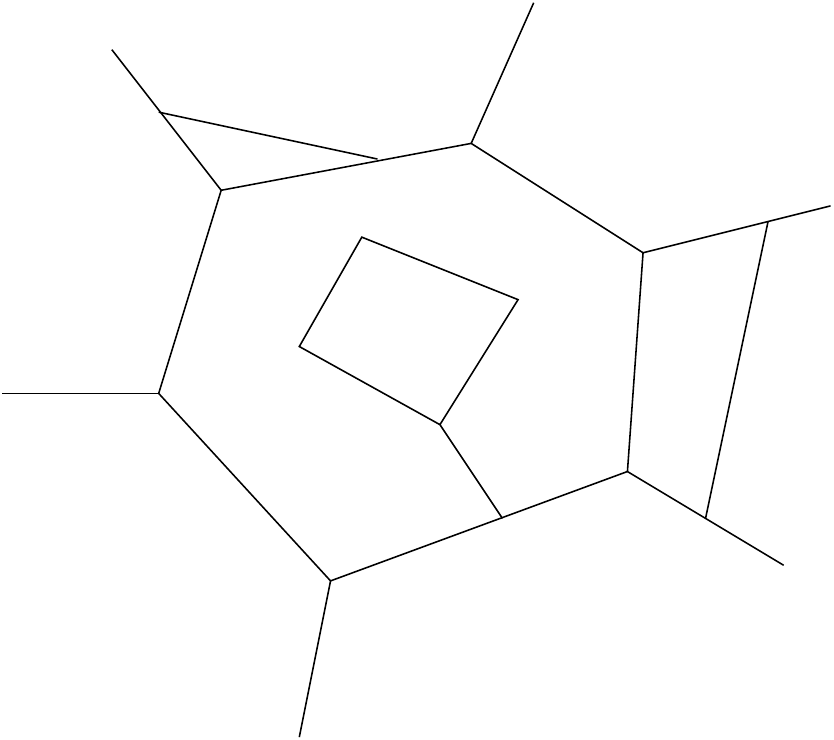}}
\caption{}\label{fig3}
\end{figure}

Une
{\it d\'ecomposition en briques} d'une surface orientable $M$ (non
n\' ecessairement connexe) est la donn\'ee d'un ensemble
stratifi\'e
$\Sigma({\cal D})$ de dimension un, appel\'e {\it 
squelette} de la d\'ecomposition
$\cal D$, avec une sous-vari\'et\'e de dimension z\'ero $S$ tel que de
tout {\it sommet} $s\in S$ sont issues (localement)
exactement trois ar\^etes. Les {\it briques} sont 
les adh\' erences des composantes
connexes de
$M\setminus\Sigma({\cal D})$ et les {\it 
ar\^etes} les adh\' erences des composantes
connexes de
$\Sigma({\cal D})\setminus S$. On
notera
$A$ l'ensemble des ar\^etes et $B$ l'ensemble des 
briques et on \'ecrira ${\cal D}=(S,A,B)$ pour 
d\'esigner
une d\'ecomposition en briques. La \fullref{fig3} donne un exemple de
d\'ecomposition en briques.

Soit ${\cal D}=(S,A,B)$ une d\'ecomposition en 
briques de $M$. On dit qu'une partie $X$ de $B$ 
est
{\it connexe} si, pour toutes briques $b, b'\in 
X$, il existe une suite $(b_i)_{0\leq
i\leq n}$, o\`u
$b_0=b$ et
$b_n=b'$, telle que $b_i$ et $b_{i+1}$ sont {\it adjacentes} si
$i\in\{0,\dots,n-1\}$, c'est-\`a-dire contiennent 
une ar\^ete commune. On identifiera par la suite 
une
partie
$X$ de
$B$ et la partie ferm\'ee de
$M$, r\'eunion des briques de $X$. Remarquons que
$\partial X$ est une vari\'et\'e topologique de 
dimension 1, remarquons \' egalement que la
connexit\' e de $X\subset B$ \' equivaut \`a la 
connexit\' e de $X\subset M$ ainsi qu'\`a celle 
de ${\rm
Int}(X)\subset M$.

On dit que la d\' ecomposition ${\cal D}'$ est 
une {\it sous-d\'ecomposition} de ${\cal D}$ si
$\Sigma({\cal D}')$ est inclus dans $\Sigma({\cal 
D})$. Tout sommet de ${\cal D}'$ est un sommet de
$\cal D$, toute ar\^ete de
${\cal D}'$ est r\'eunion d'ar\^etes de ${\cal 
D}$ et toute brique de ${\cal D}'$ r\'eunion de
briques de ${\cal D}$. Si $(X_i)_{i\in I}$ est 
une partition de $B$ en parties connexes, alors
$\bigcup_{i\in I} \partial X_i$ est le squelette 
d'une sous-d\' ecomposition ${\cal D}'$ de $\cal 
D$
dont les briques sont les $X_i$. En particulier 
si $(X_i)_{i\in I}$ est la partition en points, 
les
briques de ${\cal D}'$ et les briques de $\cal D$ 
sont les m\^emes. Si ${\cal D}'$ et $\cal D$
co\"\i ncident, on dira que $\cal D$ est une d\' 
ecomposition {\it remplie}: ceci signifie que 
toute
ar\^ete de la d\' ecomposition est incluse dans exactement deux
briques (ainsi la d\'ecomposition repr\'esent\'ee sur la \fullref{fig3}
n'est pas remplie).

Si $f$ est un hom\' eomorphisme de $M$, on 
d\'efinit naturellement une application\break $\varphi\co 
{\cal
P}(B)\to {\cal P}(B)$ en posant:
$$\eqalign{ \varphi(X)&=\{b\in 
B\enskip\vert\enskip {\rm il \enskip 
existe}\enskip b'\in
X\enskip{\rm tel\enskip que}\enskip b\cap f(b')\not=\emptyset\}\cr
&=\{b\in B\enskip\vert\enskip b\cap f(X)\not=\emptyset\enskip\}.\cr}$$
Remarquons que $\varphi(X)$ est connexe si c'est 
le cas de $X$. Remarquons \'egalement
que pour toute famille $(X_i)_{i\in I}$ de parties de $B$, on a:
$$\varphi(\bigcup_{i\in I} X_i)=\bigcup_{i\in I}
\varphi(X_i),\enskip\varphi(\bigcap_{i\in I}
X_i)\subset\bigcap_{i\in I}
\varphi(X_i).$$
On d\'efinit de fa\c con analogue une application 
$\varphi_-\co  {\cal P}(B)\to {\cal P}(B)$ en
posant:
$$\eqalign{ \varphi_-(X)&=\{b\in 
B\enskip\vert\enskip {\rm il \enskip 
existe}\enskip b'\in
X\enskip{\rm tel\enskip que}\enskip b\cap f^{-1}(b')\not=\emptyset\}\cr
&=\{b\in B\enskip\vert\enskip b\cap f^{-1}(X)\not=\emptyset\}.\cr}$$
On d\' efinit le {\it futur} $b_{\geq}$ et le {\it pass\' e} $b_{\leq}$ d'une
brique $b$ en posant:
$$b_{\geq}=\bigcup_{k\geq 0}\varphi^k(\{b\}),\enskip b_{\leq}=\bigcup_{k\geq
0}\varphi_-^k(\{b\}).$$
On d\' efinit de m\^eme le {\it futur strict} 
$b_{>}$ et le {\it pass\' e strict} $b_{<}$ en
posant:
$$b_{>}=\bigcup_{k\geq 1}\varphi^k(\{b\}),\enskip
b_{<}=\bigcup_{k\geq 1}\varphi_-^k(\{b\}).$$
Remarquons que les conditions suivantes sont \' equivalentes:
$$b\in b_{>},\enskip b_{>}=b_{\geq},\enskip b\in
b_{<},\enskip b_{<}=b_{\leq},\enskip b_{<}\cap
b_{\geq}\not=\emptyset,\enskip b_{\leq}\cap b_{>}\not=\emptyset.$$
L'existence d'une brique $b\in B$ pour laquelle 
ces conditions sont v\' erifi\' ees \' equivaut 
\`a
l'existence d'une cha\^\i ne ferm\' ee constitu\' ee de briques.

On dira qu'une d\'ecomposition en briques ${\cal
D}=(S,A,B)$ est {\it libre} si toute brique de 
${\cal D}$ est libre. Si $f$ n'a pas de point 
fixe, il est toujours possible, quitte \`a 
choisir les briques assez
petites, de construire une d\'ecomposition en briques libre. On d\' eduit
imm\' ediatement de la \fullref{prop1.5}:

\begin{prop}\label{prop2.1}Soit  $f\in{\rm
Homeo}^+(\D)$ un hom\' eomorphisme diff\' erent 
de l'identit\' e et ${\cal D}=(S,A,B)$ une d\' 
ecomposition en
briques libre de $\D\setminus{\rm Fix}(f)$. S'il 
existe $b\in B$ tel que $b\in b_{>}$, alors
$f$ est r\' ecurrent.\end{prop}

  Rappelons la d\'efinition de {\it 
d\'ecomposition libre maximale}, notion 
introduite par
Sauzet dans sa th\`ese \cite{Sa}. Soit $f$ un 
hom\'eomorphisme sans point fixe d'une
surface $M$ et $\cal D$ une d\'ecomposition en 
briques libre. On peut trouver une
sous-d\'ecomposition
${\cal D}'$ de $\cal D$ qui est libre et telle 
que toute sous-d\'ecomposition stricte de ${\cal 
D}'$
n'est plus libre. On dit que ${\cal D}'$ est une 
d\'ecomposition en briques libre maximale.
Remarquons que toute d\'ecomposition en briques 
libre maximale est remplie. Remarquons \' 
egalement
que la r\' eunion de deux briques adjacentes 
n'est pas libre. Dans sa th\`ese, Sauzet montre, 
dans le
cas o\`u
$M=\D$, que toute brique $b$ est adjacente \`a 
une brique qui rencontre $f(b)$ et \`a une brique
qui rencontre
$f^{-1}(b)$. Nous allons prouver un r\' esultat similaire si $M$ est le
compl\' ementaire de l'ensemble des points fixes 
d'un hom\' eomorphisme non trivial $f\in{\rm
Homeo}^+(\D)$. Plus pr\' ecis\' ement:

\begin{prop}\label{prop2.2}Soit $f\in{\rm
Homeo}^+(\D)$ un hom\' eomorphisme diff\' erent 
de l'identit\' e et ${\cal D}=(S,A,B)$ une d\' 
ecomposition
en briques libre maximale de $\D\setminus{\rm Fix}(f)$. Alors, l'une
au moins des assertions suivantes est vraie:

\bee
\item il existe une 
cha\^\i ne ferm\' ee de briques de $B$ (et $f$ 
est donc
r\' ecurrent);

\item toute brique 
$b$ est adjacente \`a une brique qui rencontre 
$f(b)$ et
\`a une brique qui rencontre $f^{-1}(b)$.
\ene

\noindent Dans le second cas, on en d\' eduit que les parties $b_{\geq}$,
$b_{>}$, $b_{\leq}$ et $b_{<}$ sont connexes.
\end{prop}

\demo
Expliquons d'abord le dernier
point. Si {\rm(ii)} est v\' erifi\' e, les ensembles
$\{b\}\cup\varphi(\{b\}$ et 
$\{b\}\cup\varphi_-(\{b\})$ sont connexes et on 
peut \' ecrire
$$b_{\geq}=\bigcup_{k\geq 0} \varphi^k(\{b\})=\bigcup_{k\geq 0}
\varphi^k\left(\{b\}\cup\varphi(\{b\})\right).$$ Ainsi, $b_{\geq}$ est
connexe, \'etant
l'union de parties
connexes $X_k=\varphi^k\left(\{b\}\cup\varphi(\{b\})\right)$
telles que
$X_k\cap X_{k+1}\not=\emptyset$. De m\^eme on peut \'ecrire
$$b_{>}=\bigcup_{k\geq 1}
\varphi^k\left(\{b\}\cup\varphi(\{b\})\right)$$ ainsi que $$
b_{\leq}=\bigcup_{k\geq 0} 
\varphi_-^k\left(\{b\}\cup\varphi_-(\{b\})\right),\enskip
b_{<}=\bigcup_{k\geq 1} \varphi_-^k\left(\{b\}\cup\varphi_-(\{b\})\right).$$
Prouvons maintenant que l'une des assertions {\rm 
(i)} ou {\rm(ii)} est vraie. Commen\c cons par
com\-pac\-ti\-fi\-er $\D$
en ajoutant le bout
$\infty$ et
prolongeons
$f$ en un ho\-m\' eo\-mor\-phis\-me de la 
sph\`ere $\D\sqcup\{\infty\}$ qui fixe $\infty$. 
On sait,
d'apr\`es M. Brown et J. Kister
\cite{BK}, que toute composante connexe de
$\D\setminus {\rm Fix}(f)$ est fixe (et n'est 
donc pas libre). Ainsi la fronti\`ere de toute 
brique
$b$ est non vide. C'est une sous-vari\'et\'e topologique de dimension
un, et donc une r\'eunion de cercles et de droites de $\D\setminus{\rm
Fix}(f)$.

 Commen\c cons par le cas o\`u
$\partial b$ admet au moins deux composantes connexes
$\Gamma_1$ et
$\Gamma_2$. Fixons une ar\^ete
$a_1\subset \Gamma_1$ (resp.\ $a_2\subset 
\Gamma_2$) et notons $b_1$ (resp.\ $b_2$) la 
brique adjacente
\`a
$b$ qui contient $a_1$ (resp.\ $a_{2}$). On sait 
que ni $b\cup b_1$, ni $b\cup b_2$ n'est libre,
puisque $\cal D$ est une d\' ecomposition libre maximale. Ainsi,
$b_{1}$ et $b_{2}$ rencontrent $f(b)\cup 
f^{-1}(b)$. Si $b$ ne v\' erifie pas la condition
{\rm(ii)}, l'un des ensembles $f(b)$ ou 
$f^{-1}(b)$ rencontre \`a la fois $b_1$ et $b_2$.
On peut
supposer, sans perte de g\' en\' eralit\' e, que 
c'est $f(b)$. L'ensemble connexe
$\varphi(\{b\})$ contient donc $b_1$ et $b_2$. 
Nous pouvons donc construire deux segments 
$\gamma$ et
$\gamma'$ joignant un point $z_1\in{\rm 
int}(a_1)$ \`a un point $z_2\in{\rm int}(a_2)$, 
le premier
prenant ses valeurs, sauf aux extr\' emit\' es, 
dans l'int\' erieur de $\varphi(\{b\})$, le second
prenant ses valeurs, sauf aux extr\' emit\' es, 
dans l'int\' erieur de $b$. Le r\' eunion des deux
segments est un cercle $C\subset b_{\geq}$ qui 
rencontre transversalement chaque arc $\Gamma_1$ 
et $\Gamma_2$ en un
unique point. Les arcs $\Gamma_1$ et $\Gamma_2$ 
ne peuvent donc pas \^etre des cercles, ce sont
des droites de $\D\setminus{\rm Fix}(f)$. On peut alors trouver deux points
d'accumulation (dans
$\D\sqcup\{\infty\}$) de $\Gamma_1$ qui sont s\' 
epar\' es par $C$. Il s'agit bien s\^ur de points 
fixes de
$f$ et ils appartiennent donc \' egalement \`a 
l'adh\' erence (dans $\D\sqcup\{\infty\}$) de 
$f^{-1}(b)$. On
en d\' eduit que $f^{-1}(b)$ rencontre $C$ et 
donc que $b_{<}\cap b_{\geq}\not=\emptyset$. 
Ainsi $f$
v\' erifie la \fullref{prop1.5} et est donc r\' 
ecurrent. Remarquons que nous venons \' egalement 
de montrer
que $f$ est r\' ecurrent si
$\partial b$ a au moins trois composantes 
connexes, puisque le raisonnement pr\' ec\' edent 
s'applique \`a
au moins un couple de composantes connexes de $\partial b$.

Supposons maintenant que $\partial b$ a une seule 
composante connexe $\Gamma$ et que c'est un
cercle. Supposons, l\`a encore, que $b$ ne v\' erifie pas
la condition {\rm(ii)} et montrons qu'il existe 
une cha\^\i ne ferm\' ee de briques de $B$. Nous 
sommes
dans une situation d\' ej\`a \' etudi\' ee par 
Sauzet. Il y a un nombre fini de briques 
adjacentes \`a $b$
et on peut supposer que l'image par $f$ de 
chacune de ces briques rencontre $b$. Si $b'$ est
adjacente \`a
$b$, l'image par $f$ de l'ar\^ete commune \`a
$b$ et
$b'$ est disjointe de $b$ puisque $b$ est libre. 
Ainsi $f(b')$ rencontre n\' ecessairement une
autre brique adjacente \`a $b$. On en d\' eduit 
imm\' ediatement l'existence d'une cha\^\i ne 
ferm\' ee de
briques adjacentes \`a $b$.

Il reste \`a \' etudier le cas le plus difficile, 
o\`u $\partial b$ a une seule composante connexe 
$\Gamma$
qui est une droite de $\D\setminus{\rm Fix}(f)$. 
Nous sommes dans une situation \' egalement \' 
etudi\' ee par
Sauzet (voir \cite{Sa} ou \cite{LeC}). Nous 
allons supposer que toute brique adjacente \`a 
$b$ rencontre $f^{-1}(b)$
et en d\' eduire que $f$ est r\' ecurrent. Nous 
notons $W$ la composante connexe de
${\bf D}\setminus{\rm Fix}(f)$ qui contient $b$. Nous savons d'apr\`es
\cite{BK} qu'elle
est invariante par $f$ et qu'elle contient \'egalement $f(b)$ et $f^{-1}(b)$.
On peut bien s\^ur supposer que
$f(b)$ et
$f^{-1}(b)$ sont disjoints, puisque dans le cas 
contraire, $f$ est r\' ecurrent. Soit
$b_0$ une brique adjacente \`a
$b$. On sait que $f(b_{0})$ rencontre \`a la fois 
$b$ et $f(b)$. On peut donc trouver un segment
$\gamma_0$ joignant un point de $b$ \`a un point de $f(b)$ et contenu dans
$f(b_{0})$. Quitte \`a restreindre ce segment,
on peut supposer qu'il joint un point 
$z_0\in\partial b$ \`a un point $z'_0\in\partial 
f(b)$ et que
${\rm int}(\gamma_0)\subset W\setminus (b\cup f(b))$. On peut \'ecrire
$\Gamma=\Gamma_0^-\Gamma_0^+$ comme assemblage de 
deux demi-droites ayant $z_0$ comme extr\'emit\'e 
commune.
Fixons alors $z_{0}^-\in{\rm int}(\Gamma_0^-)$ et $z_{0}^+\in{\rm
int}(\Gamma_0^+)$ et remarquons que tout point de
$W\setminus (b\cup f(b))$, proche de $z_{0}^-$, 
est s\'epar\'e dans $W\setminus
(b\cup f(b))$ par ${\rm int}(\gamma_0)$ de tout point de
$W\setminus (b\cup f(b))$ proche de $z_{0}^+$. En effet, dans le
cas contraire, on pourrait construire un cercle 
$C\subset W\setminus (b\cup f(b))$ qui intersecte
$\gamma_{0}$ en un unique point et ceci de fa\c con transverse. Le
cercle $C$ s\'eparerait donc $b$ et $f(b)$ dans 
$\D$. Notons $\overline b$ et $\overline{f(b)}$
les adh\'erences de $b$ et $f(b)$ dans $\D$. La courbe $C$
s\'eparerait
$\overline b$ et $\overline{f(b)}$ dans $\D$ puisque $b$ et $f(b)$
sont ferm\'es dans $W$. Notons maintenant que 
$\overline b$ contient un point fixe (il suffit 
de choisir un point
d'accumulation de $\Gamma$). Ce point fixe 
appartient \'egalement \`a $\overline{f(b)}$.
Ceci contredit la propri\'et\'e de s\'eparation. Dans le cas o\`u
$f^{-1}(b)\cap \gamma_{0}\not=\emptyset$, l'hom\'eomorphisme $f$ est
r\'ecurrent. En effet, $\gamma_{0}$ est libre, puisqu'il est contenu
dans $f(b_{0})$. On peut donc \'epaissir $\gamma_{0}$ pour construire
un disque libre $D_{0}\subset W$ dont l'int\'erieur est disjoint de $b$ et
$f(b)$ et ce disque v\'erifie $f(D_{0})\cap
b=\emptyset$ si $f^{-1}(b)\cap 
\gamma_{0}\not=\emptyset$. Puisque $f(b)\cap 
D_{0}\not=\emptyset$, on peut appliquer la 
\fullref{prop1.5} pour conclure que
$f$ est r\'ecurrent. Nous supposerons donc que
$f^{-1}(b)\cap \gamma_{0}=\emptyset$. Puisque 
$f^{-1}(b)$ est connexe, on peut toujours 
supposer, quitte \`a
intervertir $\Gamma_0^-$ et $\Gamma_0^+$, qu'un point dans
$W\setminus (b\cup f(b))$ qui est  proche de 
$z_{0}^-$ est s\'epar\'e de $f^{-1}(b)$ dans
$W\setminus (b\cup f(b))$ par ${\rm int}(\gamma_0)$. Par hypoth\`ese,
toute brique $b'$ adjacente \`a $b$
intersecte
$f^{-1}(b)$, elle rencontre donc $\gamma_0$ si 
elle rencontre $\Gamma_0^-$. Puisque $\gamma_0$ 
est
une partie compacte de
$W$, il n'y a donc qu'un nombre fini de briques 
adjacentes \`a $b$ qui rencontrent $\Gamma_0^-$.
Il existe donc une brique $b_1$, adjacente \`a 
$b$, qui rencontre $\Gamma_0^-$ en des points 
arbitrairement \'eloign\'es de
$z_0$. On construit alors de fa\c con similaire un segment $\gamma_1\subset
f(b_{1})$ joignant un point
$z_1\in\partial b$ \`a un point $z'_1\in\partial f(b)$ tel que
${\rm int}(\gamma_1)\subset W\setminus (b\cup 
f(b))$. On \'ecrit alors 
$\Gamma=\Gamma_1^-\Gamma_1^+$ comme
assemblage de deux demi-droites ayant $z_1$ comme 
extr\'emit\'e commune en faisant le choix 
suivant:
$\Gamma_1^-\cap\Gamma_0^-$ est une demi-droite. On fixe alors $z_{1}^-\in{\rm
int}(\Gamma_1^-)$ et $z_{1}^+\in{\rm
int}(\Gamma_1^+)$. Si $f^{-1}(b)$ rencontre
$\gamma_1$, alors $f$ est r\'ecurrent. Sinon, la partie connexe
$b_{1}\cup f^{-1}(b)$ est disjointe de $\gamma_1$ et rencontre
$\Gamma_{1}^-$. On en d\'eduit que tout point de 
$W\setminus b\cup f(b)$ qui est
proche de $z_{1}^{+}$ est s\'epar\'e de 
$f^{-1}(b)$ dans $W\setminus (b\cup f(b))$ par 
${\rm
int}(\gamma_1)$. On sait donc qu'il n'y a qu'un nombre fini de
briques adjacentes \`a $b$ qui rencontrent 
$\Gamma_1^+$. On vient de prouver qu'il n'y a 
qu'un nombre fini de
briques adjacentes \`a
$b$. On remarque alors que le raisonnement fait 
dans le cas o\`u $\Gamma$ est un cercle est 
encore valable.
\endproof

Nous allons conclure cette section en donnant le plan de la
d\'emonstration du
\fullref{thm0.1}. Pour prouver ce th\'eor\`eme, la question qui vient
naturellement \`a l'esprit est la suivante:

\medskip
{\sl Si $f\in{\rm Homeo}^+(\D)$ v\'erifie les
hypoth\`eses {\rm(i)}, {\rm(ii)}, {\rm(iii)}, 
{\rm(iv)} du \fullref{thm0.1}, et si $\cal D$ est 
une
d\'ecomposition en briques libre maximale de 
$\D\setminus{\rm Fix}(f)$, existe-t-il une 
cha\^\i ne
ferm\'ee constitu\'ee de briques de la d\'ecomposition?}

\medskip
Nous ne savons pas r\'epondre \`a cette question. 
Nous verrons, par contre, qu'il est
possible de construire une d\'ecomposition en 
briques libre maximale $\cal D$ de 
$\D\setminus{\rm
Fix}(f)$, pour laquelle la r\'eponse est oui. Pour obtenir cette
cha\^\i ne, nous nous servirons de la connexit\'e 
du futur, du pass\'e, du futur strict et du 
pass\'e strict de chaque brique. Il
faudra cependant rendre utilisables, dans 
l'\'etude de la d\'ecomposition, les hypoth\`eses 
dynamiques et
combinatoires du \fullref{thm0.1}. Ces quatre 
hypoth\`eses portent sur des ensembles discrets, 
les
orbites des
$z_i$. La premi\`ere \'etape est de construire 
des ensembles unidimensionnels sur lesquels, on 
peut
traduire ces hypoth\`eses.

Nous commencerons dans la \fullref{sec3} par
construire pour tout $i\in\Z/n\Z$ et tout $k\in \Z$ un arc de
translation $\gamma_{i}^k$ de $f^k(z_{i})$, param\'etr\'e par
$t\in[k,k+1]$, de telle fa\c con que l'arc
$\Gamma_{i}=\prod_{k\in\Z} \gamma^k_i$ param\'etr\'e par $t\in\R$
v\'erifie:

\bei\sl
\item
$\lim_{t\to-\infty}\Gamma_{i}(t)=\alpha_{i}$;

\item
$\lim_{t\to+\infty}\Gamma_{i}(t)=\omega_{i}$;

\item pour tous r\'eels $t$ et $t'$, on a
$f(\Gamma_i(t))\not=\Gamma_i(t')$ si $t'\leq t$. 
\eni

L'arc $\Gamma_{i}$ contient l'orbite $O_{i}$ du point
$z_{i}$. Si chaque $\gamma_{i}^k$ est un arc simple, il n'en est
pas de m\^eme de chaque $\Gamma_{i}$. Dans la \fullref{sec4}, nous verrons que
l'on peut modifier les segments $\gamma^k_{i}$ pour avoir une
situation ``g\'en\'erique'' au sens suivant:

\bei\sl
\item
tout point $f^k(z_i)$ 
est un point simple de $\Gamma_i$;
\item
les arcs $\Gamma_i$ n'ont pas de points triples;
\item
$\Gamma_i\cap O_{i'}=\emptyset$, si $i\not=i'$;
\item
$\Gamma_i\cap
\Gamma_{i'}$ est fini, si $i\not=i'$;
\item
$\Gamma_i\cap
\Gamma_{i'}$ ne contient aucun point double ni de $\Gamma_i$, ni de
$\Gamma_{i'}$, si $i\not=i'$;
\item
$\Gamma_i\cap 
\Gamma_{i'}\cap\Gamma_{i''}=\emptyset$ si $i$,
$i'$ et $i''$ sont distincts.
\eni

En utilisant les hypoth\`eses combinatoires {\rm
(ii)} et {\rm(iii)} du \fullref{thm0.1}, nous d\'emontre\-rons
dans la \fullref{sec5} l'existence d'une {\it configuration de
connexit\'e critique} $(t_{i})_{i\in\Z/n\Z}$ au sens suivant:

\bei\sl
\item[{\bf (A)}]les ensembles
$\Gamma_i(]-\infty,t_{i}[)$, $i\in\Z/n\Z$, sont disjoints deux \`a
deux;

\item[{\bf (B)}]pour tout 
$i\in\Z/n\Z$, il existe $i'\not\in\{i-1,i\}$ tel 
que
$$\Gamma_i(]-\infty,t_{i}])\cap\Gamma_{i'}(]-\infty,t_{i'}])\not=\emptyset.$$
\eni

Il est facile de voir qu'il existe un entier $K>0$ tel que si $i$ et
$i'$ sont distincts, alors
$$\Gamma_{i}(]-\infty,-K])\cap \Gamma_{i'}=\Gamma_{i}([K, +\infty[)\cap
\Gamma_{i'}=\emptyset.$$
On en d\'eduit que $\vert t_{i}\vert<K$ pour tout
$i\in\Z/n\Z$. La
propri\'et\'e {\bf (B)} implique alors que l'ensemble
$\Gamma_{i-1}(]-\infty,-K])$ est s\'epar\'e dans $\D$ de
$\Gamma_{i-1}([K,+\infty[)$ par
$\Gamma_i(]-\infty,t_{i}])\cup\Gamma_{i'}(]-\infty,t_{i'}])$. Nous
verrons que l'on peut de plus supposer que pour tout $i\in\Z/n\Z$, il
existe $i'\not=i$ tel que
$\Gamma_{i}(t_{i})\in\Gamma_{i'}(]-\infty,t_{i'}[)$.

Dans la \fullref{sec6}, nous construirons \`a partir des ensembles
$\Gamma_{i}(]-\infty,t_{i}])$ une famille d'arbres
$(A_{i}^-)_{i\in\Z/n\Z}$, o\`u $A_{i}^-\subset
\Gamma_{i}(]-\infty,t_{i}])$ sera obtenu en enlevant des
boucles de $
\Gamma_{i}(]-\infty,t_{i}])$. L'arbre $A_{i}^-$ contiendra tous
les points $f^k(z_{i})$, $k\leq t_{i}$, l'extr\'emit\'e $\Gamma_{i}(t_{i})$
ainsi que toutes les extr\'emit\'es $\Gamma_{i'}(t_{i'})$ qui sont
sur $\Gamma_{i}(]-\infty,t_{i}])$. Ainsi:

\bei\sl
\item[{\bf (A$'$)}]les ensembles
$A_{i}^-\setminus\{\Gamma_{i}(t_{i})\}$ sont disjoints deux \`a
deux;

\item[{\bf (B$'$)}]pour tout 
$i\in\Z/n\Z$, il existe $i'\not\in\{i-1,i\}$ tel 
que
$A_{i}^-\cap A_{i'}^-\not=\emptyset$.
\eni

Nous construirons de fa\c con analogue une famille d'arbres
$(A_{i}^+)_{i\in\Z/n\Z}$, o\`u $A_{i}^+\subset
\Gamma_{i}([K,+\infty[)$ contiendra tous
les points $f^k(z_{i})$, $k>K$. Chaque arbre $A_{i}^+$ sera
disjoint de tous les arbres $A_{i'}^-$, ainsi que de
tous les $A_{i'}^+$, $i'\not=i$. La propri\'et\'e  {\bf
(B$'$)} implique qu'il existe une composante connexe de $A_{i}^-\cup
A_{i'}^-$ qui contient tous les points $f^k(z_{i-1})$, $k\leq -K$, et
une autre composante connexe qui contient l'arbre $A_{i-1}^{+}$ et
donc tous les points $f^k(z_{i-1})$, $k> K$.

 Le caract\`ere ``presque disjoint'' de tous ces arbres,
illustr\'e par la propri\'et\'e {\bf (A$'$)}, va nous permettre de
construire dans la \fullref{sec7} une d\'ecomposition en
briques libre maximale model\'ee sur ces arbres. Cette construction
se fera en plusieurs \'etapes, la premi\`ere \'etape consistant \`a
\'epaissir les arbres, ce qui sera possible justement gr\^ace \`a la
propri\'et\'e  {\bf (A$'$)}. La propri\'et\'e essentielle v\'erifi\'ee
par la d\'ecomposition finale ${\cal D}''=(S'',A'',B'')$ sera la suivante:

\medskip{\sl
Pour tout
$i\in\Z/n\Z$, il existe une famille de briques $(b''{}_{i}^l)_{l\not=0}$
telle que:

\bei
\item
$A_{i}^-\subset\bigcup_{l<0} b_{i}''{}^l$;
\item
$A_{i}^+\subset\bigcup_{l>0} b_{i}''{}^l$;
\item
$\Gamma_{i}(t_{i})\in b_{i}''{}^{-1}$;
\item
$b''{}_i^{l'}\in 
(b''{}_i^{l})_{\geq}$  si $0<l<l'$, si
$l<l'=-1$ ou si $l<0<l'$.
\eni}

Gr\^ace \`a la propri\'et\'e {\bf (B$'$)}, ou plut\^ot \`a la propri\'et\'e de
s\'eparation qu'elle implique, on montrera dans la \fullref{sec8} que le futur
$(b''{}_{i-1}^{-1})_{\geq}$ de $b''{}_{i-1}^{-1}$, qui peut
\^etre suppos\'e
connexe d'apr\`es la \fullref{prop2.2}, et qui contient
$b''{}_{i-1}^{-1}$, doit contenir une brique
$b''{}_{i}^{l}$, $l<0$, ou une brique $b''{}_{i'}^{l'}$, $l'<0$.
Ainsi $(b''{}_{i-1}^{-1})_{\geq}$ contient 
$b''{}_{i}^{-1}$ ou $b''{}_{i'}^{-1}$.  Dans le 
cas o\`u les briques
$b''{}_{i}^{-1}$, $i\in\Z/n\Z$, sont toutes distinctes, on en d\'eduit
imm\'ediatement qu'il existe une cha\^\i ne ferm\'ee de briques. Il
restera \`a \'etudier le cas o\`u deux au moins de ces briques sont
\'egales. Un argument de connexit\'e du m\^eme type permettra de
prouver que cette brique $b''$ v\'erifie $b''\in b''_{>}$.

Il faut noter que les propositions que nous verrons dans la suite de
l'article sont du m\^eme type que la \fullref{prop2.2}. En d'autres
termes, si dans la construction pr\'ec\'edente, on ne peut pas passer
d'une \'etape \`a l'\'etape suivante, c'est que $f$ est r\'ecurrent.

\section{Construction des arcs $\Gamma_i$}\label{sec3}

Nous fixons jusqu'\`a la fin de l'article un hom\'eomorphisme
$f\in{\rm Homeo}^+(\D)$ qui v\'erifie les 
hypoth\`eses du \fullref{thm0.1}.  Nous notons 
alors $O_i$
l'orbite de $z_i$ et posons $z_i^k=f^k(z_i)$. 
Nous prouverons dans cette section le r\'esultat
suivant:

\begin{prop}\label{prop3.1}L'une au moins
des assertions suivantes est vraie.

\bee
\item
L'hom\'eomorphisme $f$ est r\'ecurrent.

\item
Pour tout 
$i\in\Z/n\Z$, on peut construire une suite
$(\gamma_i^k)_{k\in\Z}$ de segments tels que:
\item[]$\bullet$\qua
chaque $\gamma_i^k$ est 
un arc de translation qui joint $z_{i}^k$ \`a
$z_{i}^{k+1}$;

\item[]$\bullet$\qua
$f(\gamma_i^k)\cap\gamma_i^{k'}=\emptyset$ si $k'<k$;

\item[]$\bullet$\qua
la suite $(\gamma_i^k)_{k\leq 0}$ converge
vers $\{\alpha_i\}$ pour la
distance de Hausdorff;

\item[]$\bullet$\qua
la suite $(\gamma_i^{k})_{k\geq 0}$ converge
vers $\{\omega_i\}$.
\ene\end{prop}

Pour simplifier les notations, nous omettrons dans cette section
l'indice $i$. Commen\c cons par un lemme.

\begin{lem}\label{lem3.2}On peut construire une suite
$(\gamma'{}^k)_{k\in\Z}$ de segments disjoints deux \`a deux tels que:

\bei\item$z^k\in\gamma'{}^k$;
\item$\gamma'{}^k\cap {\rm Fix}(f)=\emptyset$;
\item$f(\gamma'{}^k)\cap \gamma'{}^k\not=\emptyset$;
\item la suite $(\gamma'{}^k)_{k\leq 0}$ converge
vers $\{\alpha\}$;
\item la suite 
$(\gamma'{}^k)_{k\geq 0}$ converge vers 
$\{\omega\}$.
\eni\end{lem}

\demo
On construit ais\'ement un hom\'eomorphisme
$h\co
\D\to]-1,1[^2$ tel que
$$\eqalign{\lim p_1(z')=-1&\Leftrightarrow \enskip\lim
h^{-1}(z')=\alpha,\cr
\lim
p_1(z')=1&\Leftrightarrow \enskip\lim h^{-1}(z')=\omega,\cr}$$
o\`u $p_1$, $p_2$ sont les deux projections 
d\'efinies sur $]0,1[^2$. Par un simple argument 
de
r\'ecurrence, utilisant en particulier le fait que
$\lim_{k\to -\infty} z^k=\alpha$ et  $\lim_{k\to +\infty}
z^k=\omega$, on peut toujours modifier $h$ pour 
que $p_1$ soit injective sur $h(O)$, ce qu'on
supposera dor\'enavant. Puisque
$p_1$ est injective sur $h(O)$, puisque $\lim_{k\to -\infty}
p_{1}(h(z^k))=-1$ et puisque $\lim_{k\to +\infty}
p_{1}(h(z^k))=1$, on peut, quitte \`a composer $h$ par un hom\'eomor\-phisme de
$]-1,1[^2$ fixant chaque verticale, choisir la suite $p_2(h(z^k))_{k\in\Z}$
arbitraire. On supposera donc que cette suite est 
strictement croissante. Pour simplifier 
l'\'ecriture, on omettra $h$ dans la suite;  on 
supposera en
fait
$p_1$ et
$p_2$ d\'efinies sur $\D$, via la carte globale
$h$.

Notons
$I_k$ l'intervalle ouvert de $]-1,1[$ d\'elimit\'e par $p_1(z^k)$ et
$p_1(z^{k+1})$ et posons $U^k=I_k\times]p_2(z^k),
p_2(z^{k+1})[$. Remarquons que les suites $(\overline
U^k)_{k\leq 0}$ et $(\overline
U^k)_{k\geq 0}$ sont form\'ees de disques ferm\'es de $\D$ et
convergent respectivement vers $\{\alpha\}$
et
$\{\omega\}$.

L'ensemble des points de $\D$ dont l'image par 
$f$ est sur la m\^eme verticale qu'un
\'el\'ement de $O$ est un $F_{\sigma}$ 
d'int\'erieur vide, on peut donc choisir un point 
$z'{}^0\in U^0$
en dehors de cet ensemble. Pour les m\^emes raisons, on
peut choisir $z'{}^1\in U^{1}$ en dehors de cet ensemble, qui
de plus n'appartient pas \`a la verticale passant par $z'{}^{0}$, ni
\`a l'image par $f$ de cette verticale, ni \`a 
son image inverse. En r\'eit\'erant cet argument, 
on peut construire par un
simple proc\'ed\'e de r\'ecurrence, une suite $(z'{}^k)_{k
\in\Z}$ telle que

\bei
\item  $z'{}^k\in U^k$;
\item $p_1(f(z'{}^k))\not 
=p_1(z^{k'})$ pour tous entiers $k$, $k'$;
\item $p_1(z'{}^k)\not =p_1(z'{}^{k'})$ si $k\not=k'$;
\item  $p_1(f(z'{}^k))\not=
p_1(z'{}^{k'})$ si $k\not=k'$.
\eni

  Puisque, par hypoth\`ese, $f$ se prolonge en un hom\'eomorphisme de
$\D\cup\{\alpha,\omega\}$, et puisque les suites  $(
z'{}^k)_{k\geq 0}$ et $(
z'{}^k)_{k\leq 0}$ convergent respectivement vers 
$\omega$ et $\alpha$, nous en d\'eduisons que
les suites  $(f(
z'{}^k))_{k\geq 0}$ et $(f(
z'{}^k))_{k\geq 0}$ convergent \'egalement 
respectivement vers $\omega$ et $\alpha$.
Quitte \`a composer  $h$ par un ho\-m\'eo\-mor\-phis\-me de
$]-1,1[^2$ fixant chaque verticale, on peut donc supposer que
$p_2(z^k)<p_2(f(z'{}^k))<p_2(z^{k+1})$. On peut alors
construire une suite $(I'{}^k)_{k\in\Z}$ 
d'intervalles ouverts relativement compacts de
$]-1,1[$, telle que:

\bei
\item  $I'{}^k$
contient $I^k$;

\item  $U'{}^k=I'{}^k\times]p_2(z^k),
p_2(z^{k+1})[$ contient $z'{}^k$ et $f(z'{}^k)$;

\item les suites $(\overline
U'{}^k)_{k\leq 0}$ et $(\overline
U'{}^k)_{k\geq 0}$ de disques ferm\'es de $\D$, convergent
respectivement vers $\{\alpha\}$
et
$\{\omega\}$.
\eni

Chaque segment $\gamma'{}^k$ que nous allons construire sera contenu
dans $U'_k\cup\{z^k\}$.
Ces segments seront donc disjoints et 
v\'erifieront les deux conditions de convergence 
demand\'ees. Il
reste \`a construire les $\gamma'{}^k$ pour qu'ils v\'erifient les trois autres
conditions.

Supposons d'abord qu'il n'y ait qu'un nombre fini de points fixes dans
$U'{}^k$. Nous pouvons alors supposer que 
$z'{}^k$ n'est pas fixe et trouver un segment 
issu de
$z^k$, \`a valeurs dans
$U'{}^k\cup\{z^k\}$, disjoint de ${\rm Fix}(f)$ 
et contenant \`a la fois $z'{}^k$ et $f(z'{}^k)$. 
Ce
segment convient.

\begin{figure}[ht!]\small
\vspace{2mm}\labellist
\pinlabel $\gamma$ [br] at 31 41
\pinlabel $\gamma'$ [t] at 38 30
\pinlabel $\gamma''$ [t] at 30 3
\pinlabel $z^k$  [r] at 2 12
\pinlabel $f(z^k)$  [l] at 116 36
\pinlabel $f(\gamma)$  [bl] at 80 41
\pinlabel $f(\gamma')$  [t] at 80 33
\pinlabel* $f(\gamma'')$  [tl] at 104 6
\endlabellist
\cl{\psfig{file=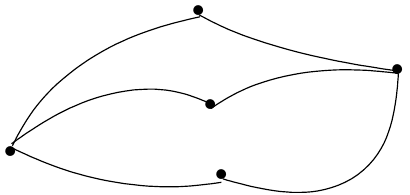,scale=120}}
\caption{}\label{fig4}
\end{figure}

Supposons maintenant qu'il y ait une infinit\'e de points fixes dans
$U'{}^k$. On peut construire dans
$U'{}^k\cup\{z^k\}$ trois segments issus de 
$z^k$, n'ayant que ce point en commun, et 
aboutissant
chacun en un point fixe. Quitte \`a restreindre 
chaque segment, on peut supposer que 
l'extr\'emit\'e
finale est l'unique point fixe rencontr\'e. Dans 
le cas o\`u l'un des segments rencontre son
image par $f$ en un point qui n'est pas 
l'extr\'emit\'e fixe, on peut trouver un
sous-segment d'extr\'emit\'e $z^k$, qui ne 
rencontre pas ${\rm Fix}(f)$ et qui n'est pas
libre. Ce sous-segment $\gamma'{}^k$ convient. Supposons maintenant que chacun
des trois segments ne rencontre son image qu'en 
son extr\'emit\'e fixe. Remarquons que si
l'image par $f$ de chacun des trois segments \'etait disjointe des
deux autres segments, alors $f$
devrait renverser l'orientation (voir \fullref{fig4}).

Puisque
$f$ pr\'eserve l'orientation, on peut donc 
trouver deux segments, parmi les trois, dont la 
r\'eunion
n'est pas libre. Enlevons \`a chacun de nos deux 
segments un petit voisinage de l'extr\'emit\'e 
fixe, la
r\'eunion des deux segments diminu\'es est un 
segment $\gamma'{}^{k}$ qui convient.
\endproof

On en d\'eduit alors:

\begin{lem}\label{lem3.3}On peut construire une suite
$(D'{}^k)_{k\in\Z}$ de disques ferm\'es disjoints deux \`a deux tels que:

\bei\item $z^k\in\partial
D'{}^k$;

\item $D'{}^k\cap {\rm Fix}(f)=\emptyset$;

\item $f(D'{}^k)\cap D'{}^k\not=\emptyset$;

\item $f({\rm 
Int}(D'{}^k))\cap {\rm Int }(D'{}^k)=\emptyset$;

\item la suite $(D'{}^k)_{k\leq 0}$ converge vers
$\{\alpha\}$:

\item la suite 
$(D'{}^k)_{k\geq 0}$ converge vers $\{\omega\}$.
\eni\end{lem}

\demo On 
peut \'epaissir chaque arc $\gamma'{}^k$ 
construit dans
le \fullref{lem3.2} pour construire une suite
$(D'{}^k)_{k\in\Z}$ de disques ferm\'es disjoints deux \`a deux
v\'erifiant toutes les conditions demand\'ees sauf la quatri\`eme.
Quitte \`a diminuer les disques que l'on vient de construire, on peut
cependant supposer que l'int\'erieur de
chaque
$D'{}^k$ est libre et donc que $f(D'{}^k)\cap D'{}^k=f(\partial D'{}^k)\cap
\partial D'{}^k$. En effet, on
peut \'ecrire
$D'{}^k\setminus\{z^k\}=\bigsqcup_{t\in]0,1]}\big(\partial 
D_t^{k}\setminus\{z^k\}\big)$,  o\`u
$t\mapsto D_t^{k}$ est une fonction continue 
croissante de disques v\'erifiant:

\bei
\item
$z^k\in \partial D_t^{k}$;

\item  $D_1^{k}=D'{}^k$;

\item  $\lim_{t\to 0} D_t^{k} =\{z^k\}$;

\item ${\rm 
Int}(D_t^{k})=\bigsqcup_{t'<t}\big(\partial 
D_{t'}^{k}\setminus\{z^k\}\big)$.
\eni

Il existe alors un plus
petit r\'eel
$t\in]0,1]$ tel que $D_t^{k}$ ne soit pas libre. 
On peut v\'erifier que son int\'erieur est
n\'ecessairement libre. \endproof

Continuons notre s\'erie de lemmes.

\begin{lem}\label{lem3.4}L'une au moins
des assertions suivantes est vraie.

\bee\item L'hom\'eomorphisme $f$ est r\'ecurrent.

\item On peut construire une suite
$(D^k)_{k\in\Z}$ de disques ferm\'es disjoints deux \`a deux tels que:

\item[]$\bullet$ $z^k\in{\rm Int}( D^k)$;

\item[]$\bullet$ $D^k\cap {\rm Fix}(f)=\emptyset$;

\item[]$\bullet$ $f(D^k)\cap D^k\not=\emptyset$;

\item[]$\bullet$ $f(D^k)\cap 
D^{k'}=f^2(D^k)\cap D^{k'}=\emptyset$ si $k'<k$;

\item[]$\bullet$ la suite $(D^{k})_{k\leq 0}$ converge vers
$\{\alpha\}$;

\item[]$\bullet$ la suite $(D^k)_{k\geq 
0}$ converge vers $\{\omega\}$.\ene\end{lem}

\demo
Consid\'erons une suite
$(D'{}^k)_{k\in\Z}$ d\'efinie par le \fullref{lem3.3}. 
Tout segment $\gamma$ joignant le point $z^k$ \`a 
un point $z'\in\partial D'{}^k$  et
contenu, sauf en ses extr\'emit\'es, dans l'int\'erieur du disque
$D'{}^k$, est libre. En effet $z_{k}$ n'\'etant 
pas fixe, le segment $f(\gamma)$ ne
peut rencontrer $\gamma$ que si l'extr\'emit\'e 
$z'$ est fixe, \'egale \`a $f(z^k)=z^{k+1}$ ou 
\'egale \`a
$f^{-1}(z^k)=z^{k-1}$. Le premier cas est impossible puisque
$D'{}^k\cap{\rm Fix}(f)=\emptyset$, les deux 
autres cas le sont \'egalement pusique
$D'{}^{k+1}$ et $D'{}^{k-1}$ sont disjoints de $D'{}^{k}$ et
contiennent respectivement $z^{k+1}$ et $z^{k-1}$.

Choisissons
$z''{}^k\in \partial D'{}^k\cap f^{-1}(\partial 
D'{}^k)$, puis deux segments $\gamma_-^{k}$ et
$\gamma_+^{k}$, le premier joignant $z''{}^k$ \`a 
$z^k$, le second joignant $z^k$ \`a $f(z''{}
^k)$, tous deux contenus dans l'int\'erieur de 
$D'{}^k$, sauf aux extr\'emit\'es, et ne se 
rencontrant
qu'en
$z^k$. Fixons $k'<k$. L'orbite positive de chacun des segments
$\gamma_-^{k'}$ et $\gamma_+^{k'}$
rencontre $\gamma_-^{k}$ et $\gamma_+^{k}$ (au 
moins au point $z^{k}$). Puisque tous ces
segments sont libres, si l'orbite positive de $\gamma_-^{k}$ ou de
$\gamma_+^{k}$ rencontre $\gamma_-^{k'}$ ou 
$\gamma_+^{k'}$, alors $f$ est r\'ecurrent, 
d'apr\`es la
\fullref{prop1.4}. Ainsi l'assertion {\rm(i)} est vraie. Supposons
maintenant  que l'orbite positive de $\gamma_-^{k}$ et de
$\gamma_+^{k}$ ne rencontre ni $\gamma_-^{k'}$ ni $\gamma_+^{k'}$,
d\`es que $k'<k$. Puisque
$$\lim_{k\to-\infty} \gamma_-^{k}\gamma_+^{k}=\alpha, \enskip
\lim_{k\to+\infty} \gamma_-^{k}\gamma_+^{k}=\omega,$$ on peut trouver
un disque ferm\'e $D^0$,
voisinage de
$\gamma_-^{0}\gamma_+^{0}$, tel que

\bei
\item
$D^0\cap{\rm Fix}(f)=\emptyset$;
\item
$D^0\cap\gamma_-^{k}\gamma_+^{k}=f(D^0)\cap\gamma_-^{k}\gamma_+^{k}=f^2(D^0)\cap\gamma_-^{k}\gamma_+^{k}=\emptyset$,
si $k<0$;
\item
$D^0\cap\gamma_-^{k}\gamma_+^{k}=f^{-1}(D^0)\cap\gamma_-^{k}\gamma_+^{k}=f^{-2}(D^0)\cap\gamma_-^{k}\gamma_+^{k}=\emptyset$, 
si
$k>0$.
\eni
On a alors n\'ecessairement:
\bei
\item
$z^{0}\in{\rm Int}(D^{0})$;
\item
$f(D^0)\cap D^{0}\not=\emptyset$.
\eni
On peut choisir ensuite un disque ferm\'e $D^1$,
voisinage de
$\gamma_-^{1}\gamma_+^{1}$, tel que
\bei
\item
$D^1\cap{\rm Fix}(f)=\emptyset$;
\item
$D^1\cap\gamma_-^{k}\gamma_+^{k}=f(D^1)\cap\gamma_-^{k}\gamma_+^{k}=f^2(D^1)\cap\gamma_-^{k}\gamma_+^{k}=\emptyset$,
si $k<0$;
\item
$D^1\cap D^{0}=f(D^1)\cap D^{0}=f^2(D^1)\cap D^{0}=\emptyset$;
\item
$D^1\cap\gamma_-^{k}\gamma_+^{k}=f^{-1}(D^1)\cap\gamma_-^{k}\gamma_+^{k}=f^{-2}(D^1)\cap\gamma_-^{k}\gamma_+^{k}=\emptyset$,
si $k>1$.
\eni
On a alors:
\bei
\item
$z^{1}\in{\rm Int}(D^{1})$;
\item
$f(D^1)\cap D^{1}\not=\emptyset$.
\eni
Par un simple proc\'ed\'e de r\'ecurrence on construit de
cette fa\c con notre
suite $(D^k)_{k\in\Z}$. \endproof

Nous pouvons maintenant prouver la \fullref{prop3.1}
en utilisant des arguments classiques dans le
domaine.

\proof[D\'emonstration de la \fullref{prop3.1}]
On peut bien s\^ur supposer que
l'assertion {\rm(ii)} du \fullref{lem3.4} est v\'erifi\'ee et consid\'erer une suite
$(D^k)_{k\in\Z}$ d\'efinie par ce lemme. L\`a 
encore, en diminuant $D^k$ si n\'ecessaire, on
peut supposer que l'int\'erieur de chaque $D^k$ est libre. Fixons
$z'''{}^k\in\partial D^k\cap f^{-1}(\partial 
D^k)$, puis deux segments $\gamma'{}_-^{k}$ et
$\gamma'{}_+^{k}$ le premier joignant $z'''{}^k$ 
\`a $z^k$, le second joignant $z^k$ \`a $f(z'''{}
^k)$ contenus dans l'int\'erieur de $D^k$, sauf 
en une extr\'emit\'e, et ne se rencontrant qu'en
$z^k$.  Le segment 
$\gamma'{}_-^{k}\gamma'{}_+^{k}$ est un arc de 
translation. Si l'une
des hypoth\`eses de la \fullref{prop1.1} (le lemme 
de Brouwer) est v\'erifi\'ee, alors $f$ est 
r\'ecurrent,
d'apr\`es la \fullref{prop1.4}. Supposons donc que 
pour tout $k\in\Z$, aucune de ces
hypoth\`eses n'est v\'erif\'ee. Remarquons  que
$\gamma^k=\gamma'{}_+^{k}f(\gamma'{}_-^{k})$ est 
alors un arc de translation joignant $z^k$ \`a
$z^{k+1}$. Puisque $\gamma^k\subset D^k\cup f(D^k)$ et
puisque $f$ se prolonge en un hom\'eomorphisme de
$\D\cup\{\alpha,\omega\}$, les deux propri\'et\'es de convergence de la
\fullref{prop3.1} sont
v\'erifi\'ees. Si $k'<k$, les ensembles $D^{k'}\cup f(D^{k'})$ et $f( D^k)\cup
f^2(D^k)$ sont disjoints. Ainsi toutes les 
conclusions de la proposition sont v\'erifi\'ees.
\endproof

\medskip
\noindent{\bf Remarque}\qua Les trois 
premi\`eres assertions du \fullref{lem3.2} impliquent
que chaque orbite
$O_i$ est contenue dans une composante connexe de ${\bf
D}\setminus{\rm Fix}(f)$ (on retrouve un cas particulier du
th\'eor\`eme prouv\'e dans \cite{BK}). Les deux
autres assertions ainsi que les hypoth\`eses 
combinatoires sur les $\alpha_i$ et $\omega_i$ 
v\'erif\'ees
par $f$ impliquent en fait que les orbites
$O_i$ sont contenues dans la m\^eme composante 
connexe de ${\bf D}\setminus{\rm Fix}(f)$.

Puisque nous voulons montrer que $f$ est r\'ecurrent, nous supposerons
dor\'enavant que l'assertion {\rm(ii)} de la \fullref{prop3.1} est
v\'erifi\'ee. Chaque segment $\gamma_{i}^k$ sera param\'etr\'e par
$t\in[k,k+1]$. On a donc un param\'etrage par $t\in\R$ de chaque arc
$\Gamma_{i}=\prod_{k\in\Z} \gamma^k_i$ et on a
$\lim_{t\to-\infty}\Gamma_{i}(t)=\alpha$,
$\lim_{t\to+\infty}\Gamma_{i}(t)=\omega$.  L'arc $\Gamma_i$
n'est pas n\'ecessairement simple, contrairement  aux arcs restreints
$\Gamma_i\vert_{[k,k+1]}$. Le fait que chaque
$\gamma_i^k$ est un arc de translation qui joint 
$z_i^k$ \`a $z_i^{k+1}$ et dont l'image par $f$ ne
rencontre aucun segment $\gamma_i^{k'}$, $k'<k$, nous donne:

\begin{prop}\label{prop3.5}
Pour tous r\'eels $t$ et $t'$, on a
$f(\Gamma_i(t))\not=\Gamma_i(t')$ si $t'\leq t$.\end{prop}

\section{Mise sous forme g\'en\'erique}\label{sec4}

Nous allons expliquer dans cette section comment perturber les
$\gamma_i^k$, pour que les arcs $\Gamma_i$ soient
  en position ``g\'en\'erique''.

\begin{prop}\label{prop4.1}
On peut toujours supposer que les arcs
$\gamma_i^k$ v\'erifient, en plus des 
propri\'et\'es \'enonc\'ees dans la \fullref{prop3.1}, les
conditions suivantes:

\bei\item
$\gamma_i^k\cap O_{i}=\{z_i^k,z_i^{k+1}\}$;
\item$\gamma_i^k\cap O_{i'}=\emptyset$ si $i'\not=i$;
\item$\gamma_i^k\cap\gamma_{i'}^{k'}$ est fini si 
$(i',k')\not\in\{(i,k-1), (i,k),
(i,k+1)\}$;
\item$\gamma_i^k\cap\gamma_{i'}^{k'}\cap\gamma_{i''}^{k''}=\emptyset$ 
si $(i,k)$,
$(i',k')$ et
$(i'',k'')$ sont distincts.\eni\end{prop}

Commen\c cons par \'enoncer deux r\'esultats pr\'eliminaires:

\begin{lem}\label{lem4.2}Soit $\gamma$ un arc de translation
d'un point
$z\not\in{\rm Fix}(f^2)$ et
$\gamma'\subset{\rm int}(\gamma)$ un segment. Il 
existe alors un voisinage $U$ de
$\gamma'$ tel que tout segment joignant $z$ \`a 
$f(z)$ \`a valeurs dans $\gamma\cup U$ est un arc 
de
translation de $z$.\end{lem}

\demo
Puisque
$$f(\gamma')\cap \gamma'=f(\gamma')\cap \gamma=\gamma'\cap
f(\gamma)=\emptyset,$$ on peut choisir $U$ de telle fa\c con que
$$f(U)\cap U=f(U)\cap \gamma=U\cap
f(\gamma)=\emptyset.$$
On en d\'eduit alors
$$f(U\cup\gamma)\cap (U\cup\gamma)=f(\gamma)\cap \gamma=\{f(z)\}.\eqno{\qed}$$

\begin{lem}\label{lem4.3}On peut construire sur
$\R^2$ une famille $(X_p)_{p\in \N}$ de parties 
d\'enombrables denses, qui v\'erifient les
propri\'et\'es suivantes, o\`u $\Xi_p$ d\'esigne 
l'en\-sem\-ble des droites affines contenant au 
moins deux points de
$X_p$:

\bei
\item si $p$ et $p'$ sont 
distincts, alors $\Xi_p\cap\Xi_{p'}=\emptyset$;

\item si $p$, $p'$ et $p''$ sont distincts, alors
$\Delta\cap\Delta'\cap\Delta''=\emptyset$ pour tout
$(\Delta,\Delta',\Delta'')\in\Xi_p\times\Xi_{p'}\times\Xi_{p''}$.
\eni\end{lem}

\demo
On pose 
$X_0=\Q^2$. On suppose ensuite que la famille
$(X_p)_{0\leq p\leq P}$ a \'et\'e cons\-tru\-ite jusqu'\`a l'ordre $P$ et
v\'erifie les deux propri\'et\'es du \fullref{lem4.3}. On va
montrer qu'on peut adjoindre un ensemble
$X_{P+1}$ pour que ces deux propri\'et\'es soient encore v\'erifi\'ees.
Remarquons que l'ensemble
$$Y_P=\bigcup_{0\leq p<p'\leq P} 
\bigcup_{\Delta\in\Xi_p}\bigcup_{\Delta'\in\Xi_{p'}}
\Delta\cap\Delta'$$ est d\'enombrable puisque la 
premi\`ere condition est vraie jusqu'au rang $P$.
Notons ${\rm Vect}_{\Q}(A)$ le $\Q$-espace vectoriel engendr\'e par une partie
$A\subset\R^2$ et $\R A$ l'ensemble des multiples r\'eels d'\'elements de $A$.
Si $A_{1}$ et $A_{2}$ sont deux parties de $\R^{2}$, notons
$A_{1}+A_{2}$ l'ensemble des points $z\in\R^{2}$ qui s'\'ecrivent
$z=z_{1}+z_{2}$, o\`u $z_{1}\in A_{1}$ et $z_{2}\in A_{2}$.
L'ensemble
$$Z_P={\rm Vect}_{\Q}\left(Y_P\cup\left(\bigcup_{0\leq p\leq P}
X_p\right)\right)+\R\Q^2$$ est une r\'eunion
d\'enombrable de droites affines (\`a pentes 
rationnelles) et donc un $F_{\sigma}$ 
d'int\'erieur vide.
On choisit un point $x_{P+1}\not\in Z_P$ et on pose
$X_{P+1}=x_{P+1}+\Q^2$. Il est clair que $X_{P+1}$ est disjoint de
$Z_P$. On v\'erifie ais\'ement que la famille
$(X_p)_{0\leq p\leq P+1}$ v\'erifie les conditions demand\'ees. En
effet, toute droite $\Delta\in\Xi_{P+1}$ est \`a pente rationnelle et ne
contient donc aucun \'el\'ement de $X_{p}$, $p\leq P$. La premi\`ere
condition est donc v\'erifi\'ee. Pour les m\^emes 
raisons, $\Delta$ ne contient aucun
point de $Y_{P}$, ce qui implique que la seconde condition est
\'egalement v\'erifi\'ee.\endproof

Avant de prouver la \fullref{prop4.1}, remarquons d'abord que la famille
$(\gamma_i^k)_{i,k}$ est localement finie: seul 
un nombre fini de ces arcs rencontrent une partie
compacte donn\'ee de
$\D$. On utilisera implicitement ce fait dans la suite.
Remarquons ensuite que l'on peut construire, par r\'ecurrence, une suite
$(W_i^k)_{k\in\Z}$ de {\it disques de s\'ecurit\'e}:

\bei
\item chaque $W_i^k$ est un 
disque ouvert contenant $\gamma_i^k$;
\item $f(W_i^k)\cap W_i^{k'}=\emptyset$ si $k'<k$;
\item la suite $(\overline W_i^k)_{k\leq 0}$
converge vers $\{\alpha_i\}$;
\item la suite $(\overline 
W_i^k)_{k\geq 0}$ converge vers $\{\omega_i\}$.
\eni
\proof[D\'emonstration de la \fullref{prop4.1}]
Commen\c cons par perturber les 
$\gamma_i^k$ pour que les deux premi\`eres
assertions de la \fullref{prop4.1} soient v\'erifi\'ees. Fixons
$\gamma_i^{k}$. Les orbites $O_i$ sont disjointes 
deux \`a deux, ferm\'ees et discr\`etes. 
L'ensemble
$$\gamma_i^{k}\cap\left(\left(\bigcup_{i'}
O_{i'}\right)\setminus\{z_i^k,z_{i}^{k+1}\}\right)$$ est
donc fini et inclus dans un segment
$\gamma'\subset{\rm int}(\gamma_i^k)$. On choisit alors un
voisinage
$U$ de
$\gamma'$ satisfaisant le \fullref{lem4.2}, inclus dans 
le disque de s\'ecurit\'e $W_i^k$  et tel que
$$U\cap\left(\left(\bigcup_{i'}
O_{i'}\right)\setminus\{z_i^k,z_{i}^{k+1}\}\right)=\gamma_i^{k}\cap\left(\left(\bigcup_{i'}
O_{i'}\right)\setminus\{z_i^k,z_{i}^{k+1}\}\right).$$ On construit
facilement un segment $\wtilde\gamma_i^{k}$ dans $\gamma_i^k\cup U$
joignant $z_i^k$ \`a $z_i^{k+1}$ et \'evitant l'ensemble
$$\gamma_i^{k}\cap\left(\left(\bigcup_{i'}
O_{i'}\right)\setminus\{z_i^k,z_{i}^{k+1}\}\right).$$ Puisque dans
cette construction, on ne touche pas aux autres arcs, on peut faire
une simple r\'ecurrence et perturber tous les $\gamma_i^k$ pour qu'ils
v\'erifient les deux premi\`eres conditions de la
\fullref{prop4.1}. Nous supposerons maintenant que ces deux conditions
sont satisfaites et allons perturber de nouveaux les $\gamma_{i}^k$
pour que les deux derni\`eres conditions soient \'egalement
v\'erifi\'ees.

  On peut construire une
famille de disques ouverts $(\Delta_i^k)_{i,k}$, 
disjoints deux \`a deux, tels que
\bei
\item
$\Delta_i^k\subset W_i^k$;
\item
$z_i^k\in \Delta_i^k$;
\item
$\overline{\Delta_i^k}\cap\gamma_{i'}^{k'}=\emptyset$ 
si
$(i',k')\not\in\{(i,k-1),(i,k)\}$.
\eni

On choisit alors, pour tout $(i,k)$,  un
segment dans ${\rm int}(\gamma_i^k)$ avec une 
extr\'emit\'e dans $\Delta_i^k$ et
l'autre dans
$\Delta_{i+1}^k$ puis un voisinage $U_i^k\subset 
W_i^k$ de ce segment qui v\'erifie le \fullref{lem4.2}. 
On peut
supposer de plus que
$\overline \Delta_{i}^{k}\cap 
U_{i'}^{k'}=\emptyset$, si 
$(i',k')\not\in\{(i,k-1),(i,k)\}$.

On peut trouver une famille $(X_i^k)_{i,k}$ de parties d\'enombrables denses de
$\R^2$ qui v\'erifient les propri\'et\'es 
indiqu\'ees dans le \fullref{lem4.3}. On construit 
ensuite une famille
$(\wtilde\gamma_i^k)_{i,k}$ de segments, tels que
$\wtilde\gamma_i^k$ joint $z_i^k$ \`a
$z_i^{k+1}$, est \`a valeurs dans $\gamma^k_i\cup U^k_i$ et s'\'ecrit
$\wtilde\gamma_i^k=\wtilde\gamma_i^{k, 
1}\wtilde\gamma_i^{k, 2}\wtilde\gamma_i^{k,
3}$, o\`u
$\wtilde\gamma_i^{k, 1}$ (resp.\ 
$\wtilde\gamma_i^{k, 2}$, 
$\wtilde\gamma_i^{k, 3}$)
est un segment \`a valeurs dans $\Delta_i^k$, 
(resp.\ $U_i^k$, $\Delta_i^{k+1}$), et o\`u 
$\wtilde\gamma_i^{k,
2}$ est affine par morceaux avec des sommets dans $X_i^k$.

Si $(i',k')\not\in\{(i,k-1), (i,k),
(i,k+1)\}$, alors $\wtilde\gamma_i^k\cap\wtilde\gamma_{i'}^{k'}=
\wtilde\gamma_i^{k, 
2}\cap\wtilde\gamma_{i'}^{k', 2}$ est fini. De 
m\^eme, si
$(i,k)$,
$(i',k')$ et
$(i'',k'')$ sont distincts, alors
$$\wtilde\gamma_i^k\cap\wtilde\gamma_{i'}^{k'}\cap\wtilde\gamma_{i''}^{k''}=
\wtilde\gamma_i^{k, 
2}\cap\wtilde\gamma_{i'}^{k', 
2}\cap\wtilde\gamma_{i''}^{k'',
2}=\emptyset.\eqno{\qed}$$

On supposera dor\'enavant que les $\gamma_i^k$ 
v\'erifient les propri\'et\'es \'enonc\'ees dans les
propositions \ref{prop3.1} et \ref{prop4.1}.

\begin{prop}\label{prop4.4}Les arcs $\Gamma_i$ v\'erifient les 
propri\'et\'es
suivantes:
\bei
\item tout point $z_i^k$ est un point simple de $\Gamma_i$;
\item les arcs $\Gamma_i$ n'ont pas de points triples;
\item  $\Gamma_i\cap O_{i'}=\emptyset$, si $i\not=i'$;
\item $\Gamma_i\cap
\Gamma_{i'}$ est fini, si $i\not=i'$;
\item $\Gamma_i\cap
\Gamma_{i'}$ ne contient aucun point double ni de $\Gamma_i$, ni de
$\Gamma_{i'}$, si $i\not=i'$;
\item  $\Gamma_i\cap 
\Gamma_{i'}\cap\Gamma_{i''}=\emptyset$ si $i$,
$i'$ et $i''$ sont distincts.\eni\end{prop}

\demo
Il s'agit de 
cons\'equences imm\'ediates des points suivants:
\bei
\item les arcs 
$\gamma_i^k=\Gamma_i\vert_{[k,k+1]}$ sont des 
arcs simples;
\item les arcs $\gamma_i^k$ 
v\'erifient les conditions de la \fullref{prop4.1};
\item $\lim_{k\to-\infty} \gamma_i^k=\{\alpha_i\}$;
\item $\lim_{k\to+\infty} \gamma_i^k=\{\omega_i\}$;
\item les $2n$ points 
$\omega_i$, $\alpha_i$, $i\in\Z/n\Z$, sont 
distincts.\endproof\eni

\section{D\'emonstration du lemme topologique fondamental}\label{sec5}

Nous prouverons dans ce paragraphe le 
r\'esultat-cl\'e de l'article (\fullref{prop5.1}). 
Il s'agit
d'un r\'esultat purement topologique qui utilise que le fait que
$\Gamma_i$ joint $\alpha_i$ \`a $\omega_i$ ainsi 
que les hypoth\`eses combinatoires {\rm(ii)} et 
{\rm(iii)} du \fullref{thm0.1}.  Le r\'esultat est 
encore vrai sans les hypoth\`eses de 
g\'en\'ericit\'e
d\'ecrites par la \fullref{prop4.4}. Nous 
n'\'ecrirons cependant la preuve que dans le cas 
o\`u ces
hypoth\`eses sont v\'erifi\'ees, beaucoup d'arguments se simplifiant.

 Nous
noterons
$\Gamma_i^{\leq t}$ (resp.\
$\Gamma_i^{<t}$) la restriction de $\Gamma_i$ \`a
$]-\infty,t]$ (resp.\ $]-\infty,t[$) et nous 
d\'efinirons de fa\c con analogue $\Gamma_i^{\geq 
t}$ et
$\Gamma_i^{>t}$. Nous ferons, en prenant
les pr\'ecautions usuelles, l'abus de langage qui 
consiste \`a identifier un arc \`a son image. 
Pour des
questions d'homog\'en\'eit\'e d'\'ecriture, nous 
\'ecrirons $\Gamma_i^t$ au lieu de $\Gamma_i(t)$. 
Nous
d\'efinirons l'ensemble des {\it configurations} ${\cal
T}=\R^{\Z/n\Z}$ et pour tout $i_0\in\Z/n\Z$, nous \'ecrirons
$p_{i_0}\co T=(t_i)_{i\in\Z/n\Z}\mapsto t_{i_0}$ pour la
$i_0$-\`eme projection.

Si $K$ est un entier assez grand et si $i$ et
$i'$ sont distincts, alors
$$\Gamma_{i}(]-\infty,-K])\cap \Gamma_{i'}=\Gamma_{i}([K,+\infty[)\cap
\Gamma_{i'}=\emptyset.$$
Nous 
fixons un tel entier $K$ jusqu'\`a la fin de
l'article.

\begin{prop}\label{prop5.1}Il existe $T=(t_i)_{i\in\Z/n\Z}\in{\cal T}$ 
tel
que:
\bei
\item[{\bf (A)}]les ensembles
$\Gamma_i^{<t_i}$, $i\in\Z/n\Z$, sont disjoints deux \`a
deux;
\item[{\bf (B)}]pour tout 
$i\in\Z/n\Z$, il existe $i'\not\in\{i-1,i\}$ tel 
que
$\Gamma_i^{\leq t_i}\cap\Gamma_{i'}^{\leq t_{i'}}\not=\emptyset$.
\eni\end{prop}

Remarquons que si les $t_i$ sont tous tr\`es 
petits, la condition {\bf (A)} est v\'erifi\'ee 
mais
pas la condition {\bf (B)}. Si, par contre, les 
$t_i$ sont tous tr\`es grands, c'est {\bf (B)} qui
est v\'erifi\'ee, mais pas {\bf (A)}. Nous devons 
donc chercher une configuration critique,
correspondant \`a une situation de connexit\'e 
limite. Nous utiliserons un proc\'ed\'e de type
minimax.

Pour tout
$T=(t_i)_{i\in\Z/n\Z}\in{\cal T}$, notons $\Sigma(T)$ l'ensemble des
$i\in\Z/n\Z$ tels qu'il existe
$i'\not\in\{i-1,i\}$ v\'erifiant
$\Gamma_i^{\leq t_i}\cap\Gamma_{i'}^{\leq 
t_{i'}}\not=\emptyset$ et appelons {\it ordre} de 
$T$ le
cardinal de $\Sigma(T)$. D'apr\`es les hypoth\`eses {\rm
(ii)} et {\rm(iii)} du \fullref{thm0.1}, on peut remarquer 
que si $i\in\Sigma(T)$ et si $\Gamma_i^{\leq
t_i}$ rencontre
$\Gamma_{i'}^{\leq {t_{i'}}}$, 
$i'\not\in\{i-1,i\}$, alors les arcs 
$\Gamma_{i-1}^{\leq -K}$ et
$\Gamma_{i-1}^{\geq K}$ sont s\'epar\'es dans 
${\bf D}$ par $\Gamma_i^{\leq 
t_i}\cup\Gamma_{i'}^{\leq
{t_{i'}}}$. Notons \'egalement
${\cal T}_*$ l'ensemble des configurations v\'erifiant
${\bf (A)}$. Remarquons que c'est une partie ferm\'ee de $\cal T$.

\begin{lem}\label{lem5.2}Soit $T=(t_i)_{i\in\Z/n\Z}\in{\cal T}_*$.

\bee\item
Pour tout 
$i\in\Sigma(T)$, on a $t_i>-K$ et $t_{i-1}<K$.
\item
Les configurations
\begin{gather*}
T'=(t'_i)_{i\in\Z/n\Z}, \enskip T''=(t''_i)_{i\in\Z/n\Z},\enskip
T'''=(t'''_i)_{i\in\Z/n\Z},\quad \text{o\`u}\\
t'_i=\min(t_i,K),  \enskip t''_i=\max(t_i,-K),  \enskip
t'''_{i}=\min(t''_i,K)=\max(t'_i,-K),\end{gather*}
appartiennent \`a ${\cal T}_*$. De plus on a
$$\Sigma(T')=\Sigma(T'')=\Sigma(T''')=\Sigma(T).$$ \ene\end{lem}

\demo Soit 
$i\in\Sigma(T)$ et $i'\not\in\{i-1,i\}$ tel que
$\Gamma_i^{\leq t_i}\cap\Gamma_{i'}^{\leq t_{i'}}\not=\emptyset$.
L'\'egalit\'e $t_{i}>-K$ est une cons\'equence de la d\'efinition de
$K$. Prouvons par l'absurde que $t_{i-1}<K$. Dans le cas contraire,
l'arc
$\Gamma_{i-1}=\Gamma_{i-1}^{<t_{i-1}}\Gamma_{i-1}^{\geq 
t_{i-1}}$ serait disjoint de $\Gamma_i^{< t_i}$
et de $\Gamma_{i'}^{< t_{i'}}$, il s\'eparerait
ces deux arcs. Les deux arcs $\Gamma_i^{\leq 
t_i}$ et $\Gamma_{i'}^{\leq t_{i'}}$ auraient donc
une extr\'emit\'e commune 
$\Gamma_i^{t_i}=\Gamma_{i'}^{t_{i'}}$ situ\'ee 
sur $\Gamma_{i-1}$. Ceci
contredit la derni\`ere propri\'et\'e de la \fullref{prop4.4}.

Il est \'evident que $T'\in{\cal T}_{*}$ et que $\Sigma(T')\subset
\Sigma(T)$ puisque $t'_{i}\leq t_{i}$, pour tout 
$i\in\Z/n\Z$. Soit $i\in\Sigma(T)$ et 
$i'\not\in\{i-1,i\}$ tel que
$\Gamma_i^{\leq t_i}\cap\Gamma_{i'}^{\leq t_{i'}}\not=\emptyset$.
Puisque $\Gamma_{i}^{\geq K}\cap\Gamma_{i'}=\Gamma_{i'}^{\geq
K}\cap\Gamma_{i}=\emptyset$, on en d\'eduit que
$\Gamma_i^{\leq t'_i}\cap\Gamma_{i'}^{\leq t'_{i'}}\not=\emptyset$ et
donc que $i\in\Sigma(T')$.

Puisque $\Gamma_{i}^{\leq-K}\cap \Gamma_{i'}=\emptyset$, si
$i'\not=i$, on sait que
$T''\in{\cal T}_{*}$. On sait \'egalement que $\Sigma(T)\subset
\Sigma(T'')$ puisque $t_{i}\leq t''_{i}$, pour tout $i\in\Z/n\Z$. Soit
$i\in\Sigma(T'')$ et $i'\not\in\{i-1,i\}$ tel que
$\Gamma_i^{\leq t''_i}\cap\Gamma_{i'}^{\leq t''_{i'}}\not=\emptyset$.
Puisque $\Gamma_{i}^{\leq -K}\cap\Gamma_{i'}=\Gamma_{i'}^{\leq
-K}\cap\Gamma_{i}=\emptyset$, on en d\'eduit que $t''_{i}=t_{i}$ et
$t''_{i'}=t_{i'}$, puis que $i\in\Sigma(T)$.

On en d\'eduit finalement que $T'''\in{\cal T}_{*}$ et que
$\Sigma(T''')=\Sigma(T)$. Remarquons que $\vert t'''_{i}\vert\leq K$,
pour tout $i\in\Z/n\Z$. \endproof

D\'emontrer la \fullref{prop5.1}, c'est trouver une configuration
$T\in{\cal T}_*$ d'ordre
$n$. Commen\c cons d'abord par le r\'esultat suivant:

\begin{lem}\label{lem5.3}Il existe  une
configuration $T=(t_i)_{i\in\Z/n\Z}\in{\cal T}_*$
d'ordre au moins $n-1$.\end{lem}

\begin{figure}[ht!]\small\vspace{2mm}
\labellist
\pinlabel $\alpha_3$ [b] at 89 191
\pinlabel $\alpha_3$ [b] at 360 192
\pinlabel $\alpha_1$ [tl] at 145 12
\pinlabel $\alpha_1$ [tl] at 415 13
\pinlabel $\alpha_2$ [l] at 189 94
\pinlabel $\alpha_2$ [l] at 460 93
\pinlabel $\alpha_4$ [tr] at 8 54
\pinlabel $\alpha_4$ [tr] at 277 54
\pinlabel $\omega_4$ [tl] at 185 58
\pinlabel $\omega_4$ [tl] at 451 50
\pinlabel $\omega_3$ [tr] at 44 13
\pinlabel $\omega_3$ [tr] at 314 13
\pinlabel $\omega_2$ [br] at 25 158
\pinlabel $\omega_2$ [br] at 294 158
\pinlabel $\omega_1$ [bl] at 140 177
\pinlabel $\omega_1$ [bl] at 414 175
\endlabellist
\cl{\psfig{file=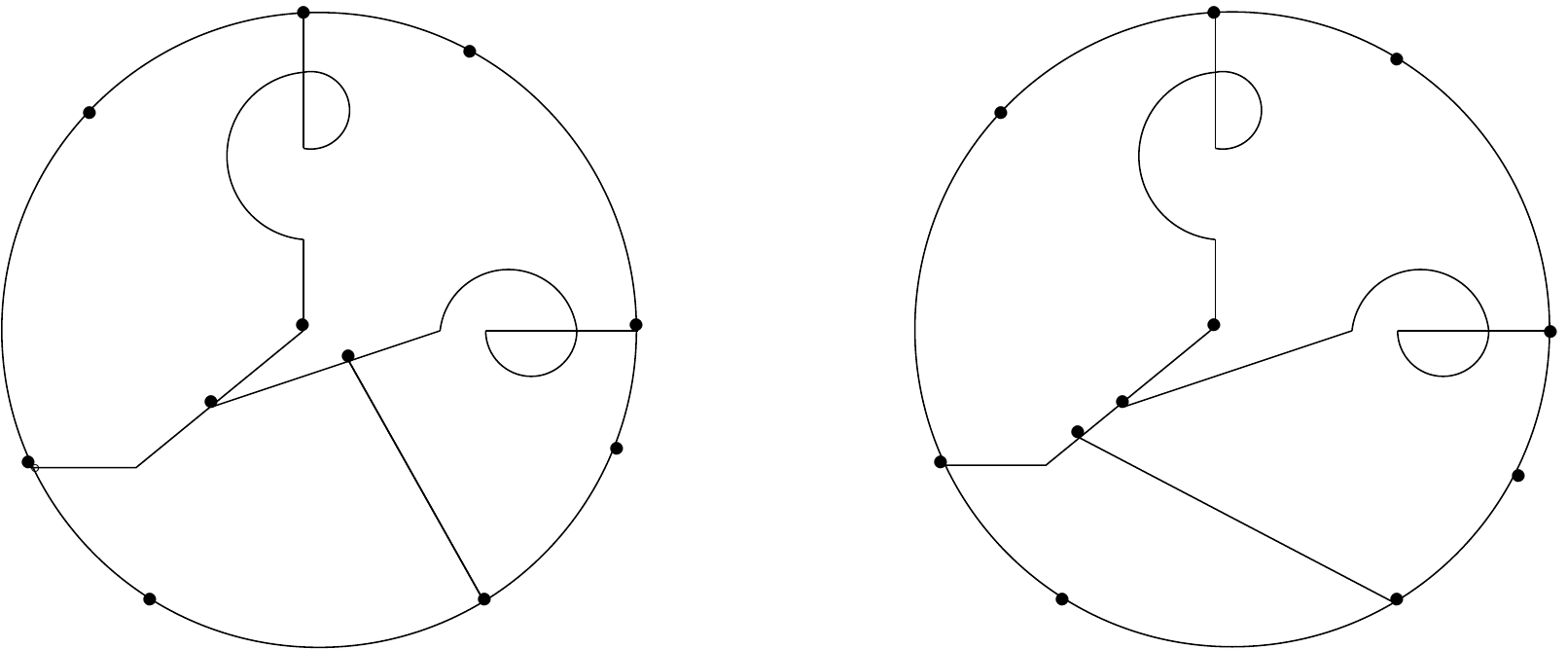,scale=70}}
\caption{}\label{fig5}
\end{figure}

\demo Fixons 
$i\in\Z/n\Z$. Utilisant les hypoth\`eses {\rm
(ii)} et {\rm(iii)} du \fullref{thm0.1}, on sait 
que $\Gamma_i^{\leq -K}$ et $\Gamma_i^{\geq K}$ 
sont
disjoints de $\Gamma_{i-1}$ et s\'epar\'es par 
cet arc. On peut donc d\'efinir le premier instant
$t_i$ o\`u
$\Gamma_{i}$ rencontre $\Gamma_{i-1}$. D'apr\`es 
la \fullref{prop4.4}, le point $\Gamma_i^{t_i}$ est
un point simple de $\Gamma_{i-1}$. On \'ecrit
$\Gamma_i^{t_i}=\Gamma_{i-1}^{t_{i-1}}$. Remarquons que
$\Gamma_{i}^{<t_i}\cap\Gamma_{i-1}^{<t_{i-1}}=\emptyset$ et que
$\Gamma_{i}^{\leq t_i}\cup\Gamma_{i-1}^{\leq 
t_{i-1}}$ s\'epare $\Gamma_{i-2}^{\leq -K}$ de
$\Gamma_{i-2}^{\geq K}$. On peut donc d\'efinir le premier instant
$t_{i-2}$ o\`u
$\Gamma_{i-2}$ rencontre $\Gamma_i^{\leq 
t_i}\cup\Gamma_{i-1}^{\leq t_{i-1}}$. On peut 
continuer le
raisonnement et d\'efinir par r\'ecurrence sur
$r\in\{2,\dots ,n-1\}$ le premier instant 
$t_{i-r}$ o\`u  $\Gamma_{i-r}$ rencontre 
$\bigcup_{0\leq
s\leq r-1}
\Gamma_{i-s}^{\leq t_{i-s}}$. La configuration 
$T=(t_i)_{i\in\Z/n\Z}$ v\'erifie {\bf (A)} et
$\Sigma(T)$ contient $\{i-n+2,\dots,i-1\}$. Remarquons maintenant que
$\Sigma(T)$ contient $i$ si $\Gamma_{i+1}^{ t_{i+1}}=\Gamma_{i-n+1}^{
t_{i-n+1}}\in \Gamma_{i}^{\leq t_{i}}$.
Si, par contre,  $\Gamma_{i+1}^{ t_{i+1}}\not\in \Gamma_{i}^{\leq t_{i}}$,
il existe $j\in\{i-n+2,\dots,i-1\}$ tel que $\Gamma_{i+1}^{ t_{i+1}}\in
\Gamma_{j}^{\leq t_{j}}$. On en d\'eduit que 
$i+1\in\Sigma(T)$. La configuration
$T$ est donc au moins d'ordre $n-1$ (les deux cas sont illustr\'es sur
la \fullref{fig5}).\endproof

Pour des raisons d'\'ecriture, il sera plus 
commode dans la suite de la \fullref{sec5}
d'indexer
$\cal T$ par
$j\in\{1,\dots,n\}$. Quitte \`a changer les indices, on peut
supposer, d'apr\`es le lemme pr\'ec\'edent, qu'il
existe
$T\in{\cal T}_*$ telle que
$j\in\Sigma(T)$, si $j\not=n$. D'apr\`es le \fullref{lem5.2}, on peut toujours
supposer que $\vert t_j\vert\leq K$, pour tout $j\in\{1,\dots,n\}$.  On notera
${\cal T}_{-1}\subset{\cal T}_*$ l'ensemble des 
configurations v\'erifiant ces deux 
propri\'et\'es.
Remarquons que c'est une partie compacte de $\cal 
T$.  On en d\'eduit que l'ensemble ${\cal T}_{0}$
form\'e des configurations $T\in{\cal T}_{-1}$ 
o\`u $p_{n}$ atteint son maximum est une
partie compacte non vide. On notera $t^*_{n}$ ce 
maximum. Nous voulons montrer qu'il existe une 
configuration d'ordre $n$ dans
${\cal T}_{-1}$. Ce sera une cons\'equence du 
\fullref{lem5.4} qui va suivre (et de la finitude de
$\{1,\dots,n\}$). Le principe est un raisonnement 
par l'absurde, on suppose qu'il n'y a pas de
configuration d'ordre $n$ dans ${\cal T}_{-1}$ et 
on construira alors par r\'ecurrence une suite
(infinie !)  d\'ecroissante de parties compactes 
$({\cal T}_{r})_{r\geq 0}$ de $\cal T$ d\'efinies
par les valeurs de
$r+1$ projections. Nous allons commencer par
initier la r\'ecurrence, pour faire ressortir les 
id\'ees, puis nous \'enoncerons et d\'emontrerons 
la
construction g\'en\'erale (\fullref{lem5.4}).

Fixons $T=(t_j)_{1\leq j\leq n}$ dans ${\cal T}_0$ et montrons le
r\'esultat suivant:

\medskip
{\sl Le point
$\Gamma_{n}^{t^*_{n}}$ appartient \`a un arc
$\Gamma_{j_1}^{\leq t_{j_1}}$,
$j_1\not=n$, cet entier $j_1$ est unique, il est \'egal \`a $n-1$ et on a
$\Gamma_{n}^{t^*_{n}}\in\Gamma_{n-1}^{<t_{n-1}}$.}

\medskip
D'apr\`es la
premi\`ere assertion du \fullref{lem5.2} appliqu\'e \`a $j=1\in\Sigma(T)$, on sait que
$t_n=t_n^*<K$. D'apr\`es la \fullref{prop4.4}, on sait qu'il existe
$t'_{n}\in]t^*_{n},K[$ tel que
$\Gamma_{n}\vert_{]t^*_{n},t'_{n}]}$
ne rencontre aucun arc $\Gamma_j$, $j\not=n$. La configuration $T'$ obtenue \`a
partir de $T$ en rempla\c cant $t^*_{n}$ par $t'_{n}$ n'est pas dans
${\cal T}_{-1}$, puisque
$t_n^*$ est le maximum de $p_n$ sur ${\cal 
T}_{-1}$. Or les coordonn\'ees de $T'$ sont toutes
entre $-K$ et $K$ et on a bien \'evidemment
$\Sigma(T)\subset\Sigma(T')$. On en d\'eduit que 
$T'$ ne v\'erifie pas la condition {\bf (A)}. Ceci
implique que
$\Gamma_{n}^{t^*_{n}}$ appartient \`a un arc
$\Gamma_{j_1}^{\leq t_{j_1}}$,
$j_1\not=n$. La propri\'et\'e 4.4 nous dit que l'entier $j_1$ est unique et que
$\Gamma_{n}^{t^*_{n}}$ est un point simple de 
$\Gamma_{j_1}$. Le fait que $T'$ ne v\'erifie pas 
la
condition {\bf (A)} implique en fait que 
$\Gamma_{n}^{t^*_{n}}$ n'est pas l'extr\'emit\'e
$\Gamma_{j_1}^{t_{j_1}}$ de
$\Gamma_{j_1}^{\leq t_{j_1}}$ et appartient donc \`a
$\Gamma_{j_{1}}^{<t_{j_{1}}}$. Le fait que $n\not\in\Sigma(T)$ implique
que $j_{1}=n-1$.

  On d\'efinit alors
$$\eqalign{t^*_{n-1}&=\min p_{n-1}\vert_{{\cal T}_{0}},\cr
{\cal T}_{1}&=\{T\in {\cal T}_{0}\enskip\vert\enskip p_{n-1}(T)=
t^*_{n-1}\}.\cr}$$
L'ensemble ${\cal T}_{1}$ est non vide, on choisit donc
$T=(t_j)_{1\leq j\leq n}\in{\cal T}_1$. Montrons le r\'esultat suivant:

\medskip
{\sl Le point
$\Gamma_{n-1}^{t^*_{n-1}}$ appartient \`a un arc
$\Gamma_{j_2}^{\leq t_{j_2}}$,
$j_2\not=n-1$, cet entier $j_2$ est unique et v\'erifie $j_{2}<n-1$.}

\medskip
Puisque $n-1\in\Sigma(T)$ on sait, d'apr\`es le \fullref{lem5.2},
que
$t_{n-1}=t^*_{n-1}>-K$. D'apr\`es la \fullref{prop4.4}, on sait qu'il existe
$t'_{n-1}\in]-K,t^*_{n-1}[$ tel que
$\Gamma_{n-1}\vert_{]t'_{n-1},t^*_{n-1}]}$
ne rencontre aucun arc $\Gamma_j$, $j\not=n-1$. Par d\'efinition de
$t^{*}_{n-1}$, on sait que la configuration $T'$ obtenue \`a
partir de $T$ en rempla\c cant $t^*_{n-1}$ par $t'_{n-1}$ n'appartient
pas \`a ${\cal T}_{0}$.
Elle v\'erifie  bien \'evidemment {\bf (A)}, puisque c'est le cas de
$T$. Puisque les coordonn\'ees de $T'$ sont toutes entre $-K$
et $K$ et puisque sa $n$-i\`eme coordonn\'ee est
$t^*_n$, c'est donc qu'il existe $j\not=n$ tel que $j\not\in\Sigma(T')$.
On en d\'eduit que
$\Gamma_{n-1}^{t^*_{n-1}}$ appartient \`a un arc
$\Gamma_{j_2}^{\leq t_{j_2}}$,
$j_2\not=n-1$, que $j_{2}$ est unique et que 
l'entier $j\not=n$ qui n'est plus dans 
$\Sigma(T')$
est soit
$n-1$ soit
$j_2$. Or on a vu pr\'ec\'edemment que le point $\Gamma_{n}^{t^*_{n}}$
n'est pas l'extr\'emit\'e
$\Gamma_{n-1}^{t^*_{n-1}}$. Il appartient donc \`a
$\Gamma_{n-1}^{<t^*_{n-1}}$ et par cons\'equent appartient \`a
$\Gamma_{n-1}^{\leq t'_{n-1}}$. On en d\'eduit 
que $n-1\in\Sigma(T')$. C'est donc l'entier $j_2$
qui n'est pas dans $\Sigma(T')$, ainsi $j_2\not=n$. On en d\'eduit que
$j_{2}<n-1$.

  On d\'efinit alors
$$\eqalign{t^*_{j_2}&=\max p_{j_2}\vert_{{\cal T}_{1}},\cr
{\cal T}_{2}&=\{T\in {\cal T}_{1}\enskip\vert\enskip p_{j_2}(T)=
t^*_{j_2}\}.\cr}$$ L'ensemble ${\cal T}_{2}$ est non vide, on choisit donc
$T=(t_j)_{1\leq j\leq n}\in{\cal T}_2$. Montrons la propri\'et\'e
suivante:

{\it  le point
$\Gamma_{j_{2}}^{t^*_{j_{2}}}$ appartient \`a un arc
$\Gamma_{j_3}^{\leq t_{j_3}}$,
$j_3\not=j_{2}$, cet entier $j_3$ est unique, il 
est \'egal \`a $j_{2}-1$ et on a
$\Gamma_{j_{2}}^{t^*_{j_{2}}}\in\Gamma_{j_{2}-1}^{<t_{j_{2}-1}}$.
}

\begin{figure}[ht!]\small\vspace{2mm}
\labellist
\pinlabel $\alpha_{n-1}$ [b] at 89 191
\pinlabel $\alpha_{j_2-1}$ [tl] at 145 12
\pinlabel $\alpha_{j_2}$ [l] at 189 94
\pinlabel $\alpha_n$ [tr] at 8 54
\pinlabel $\omega_{j_2-2}$ [tl] at 185 58
\pinlabel $\omega_{n-1}$ [tr] at 44 13
\pinlabel $\omega_{n-2}$ [br] at 25 158
\pinlabel $\omega_{n-3}$ [bl] at 140 177
\endlabellist
\cl{\psfig{file=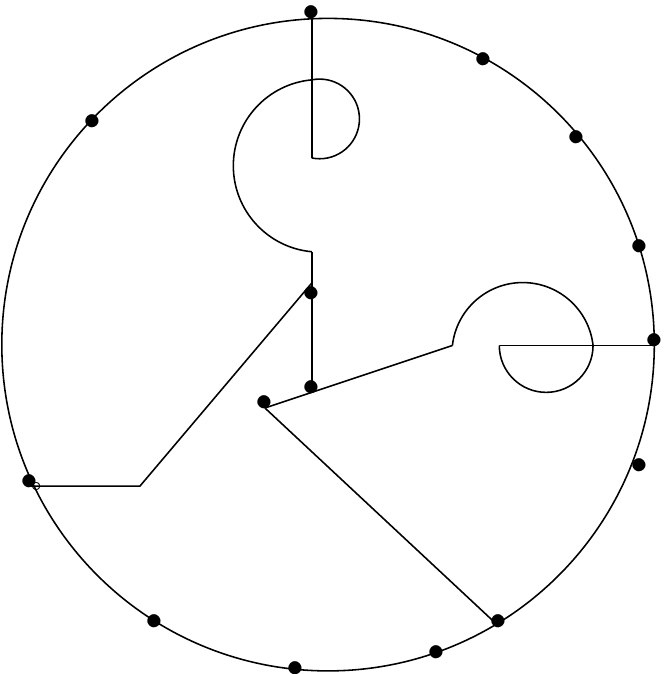,scale=70}}
\caption{}\label{fig6}
\end{figure}

On sait que $j_2\not=n-1$ et donc que $j_2+1\in\Sigma(T)$. On en d\'eduit que
$t_{j_2}=t^*_{j_2}<K$, d'apr\`es le \fullref{lem5.2}. 
Utilisant le raisonnement fait plus haut (qui 
utilise
la maximalit\'e de
$t^*_2$) on sait qu'il existe $j_3\not=j_2$,
unique, tel que
$\Gamma_{j_2}^{t^*_{j_2}}\in\Gamma_{j_{3}}^{\leq
t_{j_{3}}}$ et que
$\Gamma_{j_2}^{t^*_{j_2}}\not=\Gamma_{j_3}^{t_{j_3}}$. On \'ecrit alors
$\Gamma_{j_2}^{t^*_{j_2}}=\Gamma_{j_3}^{t''{}_{j_3}}$. Notons
$T'=(t'_j)_{1\leq j\leq n}$ la configuration obtenue \`a
partir de $T$ en rempla\c cant $t^*_{n-1}$ par 
$t'_{n-1}\in]-K,t^*_{n-1}[$, o\`u
$\Gamma_{n-1}\vert_{]t'_{n-1},t^*_{n-1}]}$
ne rencontre aucun arc $\Gamma_j$, $j\not=n-1$. Nous avons vu plus
haut que
$j_2\not\in\Sigma(T')$. Remarquons que $t''_{j_3}\leq t'_{j_3}$ dans le cas
\'eventuel o\`u
$j_3$ serait \'egal \`a $n-1$. Dans le cas contraire nous avons
l'\'egalit\'e $t''_{j_3}= t'_{j_3}$. Nous en d\'eduisons que
$\Gamma_{j_2}^{\leq t'_{j_2}}\cap \Gamma_{j_3}^{\leq t'_{j_3}}\not=\emptyset$
et par cons\'equent que
$j_3$ est n\'ecessairement \'egal \`a $j_2-1$. 
Remarquons \'egalement que les seuls arcs 
$\Gamma_j^{\leq
t_j}$,
$j\not=j_2$, rencontr\'es par
$\Gamma_{j_2}^{\leq t^*_{j_2}}$ sont
$\Gamma_{n-1}^{\leq t^*_{n-1}}$ et $\Gamma_{j_3}^{\leq t^*_{j_3}}$
(voir \fullref{fig6}).

En continuant ce processus, nous
allons montrer:

\begin{lem}\label{lem5.4}Si 
${\cal T}_{-1}$ ne contient aucune configuration
d'ordre $n$, on peut construire:

\bei
\item
une suite strictement d\'ecroissante d'entiers positifs $(j_r)_{r\geq
0}$, v\'erifiant $j_0=n$ et $j_{2r+1}=j_{2r}-1$;

\item
une suite de r\'eels $(t^*_{j_r})_{r\geq0}$;

\item
une suite d\'ecroissante de parties
compactes
$({\cal T}_{r})_{r\geq 0}$ de $\cal T$;
\vspace{-6pt}
\eni
uniquement d\'efinies par les propri\'et\'es suivantes:
\bee
\item $t^*_{j_{2r}}=\max p_{j_{2r}}\vert_{{\cal 
T}_{2r-1}}$ et
${\cal T}_{2r}=\{T\in {\cal 
T}_{2r-1}\enskip\vert\enskip p_{j_{2r}}(T)= 
t^*_{j_{2r}}\}$;
\item
$t^*_{j_{2r+1}}=\min p_{j_{2r+1}}\vert_{{\cal 
T}_{2r}}$ et
${\cal T}_{2r+1}=\{T\in {\cal 
T}_{2r}\enskip\vert\enskip p_{j_{2r+1}}(T)= 
t^*_{j_{2r+1}}\}$;
\item
il existe
un unique entier $j_{2r+2}\not=j_{2r+1}$ dans $\{1,\dots,n\}$ tel que
$\Gamma_{j_{2r+1}}^{t^*_{j_{2r+1}}}\in\Gamma_{j_{2r+2}}$.
\vspace{-6pt}
\ene
De plus, on a:
\bee
\addtocounter{enumi}{3}
\item
si 
$T=(t_j)_{1\leq j\leq n}\in{\cal T}_{2r}$, l'arc
$\Gamma_{j_{2r}}^{\leq t_{j_{2r}}}=\Gamma_{j_{2r}}^{\leq
t^*_{j_{2r}}}$ est disjoint de $\Gamma_{j}^{\leq 
t_j}$ si $j\not\in\{j_{2r-1}, j_{2r}, j_{2r+1},
\}$ (si $j\not\in\{n-1, n\}$ dans le cas o\`u
$r=0$) ainsi que de $\Gamma_{j_{2r-1}}^{<
t_{j_{2r-1}}}=\Gamma_{j_{2r-1}}^{<
t^*_{j_{2r-1}}}$;
\item
si 
$T=(t_j)_{1\leq j\leq n}\in{\cal T}_{2r}$, le 
point
$\Gamma_{j_{2r}}^{t_{j_{2r}}}=\Gamma_{j_{2r}}^{t^*_{j_{2r}}}$ appartient \`a
$\Gamma_{j_{2r+1}}^{< t_{j_{2r+1}}}$;
\item
si 
$T=(t_j)_{1\leq j\leq n}\in{\cal T}_{2r+1}$,
le point 
$\Gamma_{j_{2r+1}}^{t_{j_{2r+1}}}=\Gamma_{j_{2r+1}}^{t^*_{j_{2r+1}}}$ 
appartient \`a
$\Gamma_{j_{2r+2}}^{\leq t_{j_{2r+2}}}$.
\ene\end{lem}

\demo
On peut
supposer que les ensembles ${\cal T}_{s}$ ont \'et\'e d\'efinis
jusqu'\`a l'ordre $2r-1$, que les
entiers $(j_s)$ ont \'et\'e d\'efinis jusqu'\`a 
l'ordre  $2r$ et que les assertions du 
th\'eor\`eme sont v\'erifi\'ees
jusqu'\`a cet ordre. C'est le cas, on l'a vu, 
pour $r=1$. On d\'efinit alors $t^*_{j_{2r}}$ et 
${\cal
T}_{2r}$ par la condition {\rm(i)}. On pose ensuite
$j_{2r+1}=j_{2r}-1$ et on d\'efinit $t^*_{j_{2r+1}}$ et ${\cal
T}_{2r+1}$ par la condition {\rm(ii)}. On va commencer par montrer que
{\rm(iv)}, puis {\rm(v)}, sont v\'erifi\'ees. On 
montrera ensuite que c'est \'egalement le cas de 
{\rm(iii)} et {\rm(vi)}. Enfin, on prouvera que 
$j_{2r+2}$ est strictement plus petit que
$j_{2r+1}$.

\proof[Preuve de (iv)]Fixons 
$T=(t_j)_{1\leq j\leq n}\in{\cal
T}_{2r}$. Puisque $j_{2r-1}\not=n$ on sait que 
$j_{2r-1}\in\Sigma(T)$ et donc que
$t_{j_{2r-1}}=t^*_{j_{2r-1}}>-K$, d'apr\`es le \fullref{lem5.2}. Choisissons
$t'_{j_{2r-1}}\in]-K,t^*_{j_{2r-1}}[$ assez proche de
$t^*_{j_{2r-1}}$ pour que l'arc 
$\Gamma_{j_{2r-1}}\vert_{[t'_{j_{2r-1}},t^*_{j_{2r-1}}[}$
ne rencontre aucun arc $\Gamma_j$, $j\not=j_{2r-1}$. La configuration
$T'$ obtenue \`a partir de $T$ en rempla\c cant 
$t^*_{j_{2r-1}}$ par $t'_{j_{2r-1}}$ n'\'etant 
plus dans
${\cal T}_{2r-2}$, il existe donc $j\not=n$ qui 
n'est plus dans $\Sigma(T)$. Il s'agit soit de
$j_{2r-1}$ soit de $j_{2r}$ puisque 
$\Gamma_{j_{2r-1}}^{t^*_{j_{2r-1}}}\in\Gamma_{j_{2r}}^{\leq
t^*_{j_{2r}}}$ (d'apr\`es l'hypoth\`ese de r\'ecurrence {\rm(vi)}).
  L'hypoth\`ese de r\'ecurrence {\rm(v)} nous dit que
$\Gamma_{j_{2r-2}}^{t^*_{j_{2r-2}}}\in\Gamma_{j_{2r-1}}^{< 
t^*_{j_{2r-1}}}$ et donc que
$j_{2r-1}\in\Sigma(T')$. On en d\'eduit que 
$j_{2r}\not\in\Sigma(T')$. Ceci implique, d'une 
part que
$\Gamma_{j_{2r}}^{\leq 
t^*_{j_{2r}}}\cap\Gamma_{j}^{\leq t_j}=\emptyset$ 
si $j\not\in\{j_{2r-1},
j_{2r}, j_{2r+1},
\}$, d'autre part que $\Gamma_{j_{2r}}^{\leq
t^*_{j_{2r}}}\cap\Gamma_{j_{2r-1}}^{<
t^*_{j_{2r-1}}}=\emptyset$.

\proof[Preuve de (v)]Fixons 
$T=(t_j)_{1\leq j\leq n}\in{\cal
T}_{2r}$. Puisque $j_{2r}+1\not=n$ on sait que 
$j_{2r}+1\in\Sigma(T)$ et donc que
$t_{2j}=t^*_{2j}<K$, d'apr\`es le \fullref{lem5.2}. Choisissons
$t'_{j_{2r}}\in]t^*_{j_{2r}},K[$ assez proche de $t^*_{j_{2r}}$ pour que
$\Gamma_{j_{2r}}\vert_{]t^*_{j_{2r}},t'{}_{j_{2r}}]}$
ne rencontre aucun arc $\Gamma_j$, 
$j\not=j_{2r}$. La configuration $T'$ obtenue
\`a partir de $T$ en rempla\c cant $t^*_{j_{2r}}$ 
par $t'{}_{j_{2r}}$ ne v\'erifie plus
{\bf (A)} par maximalit\'e de $t^*_{j_{2r}}$. Il 
existe donc $j\not=j_{2r}$ tel que
$\Gamma_{j_{2r}}^{t^*_{j_{2r}}}$ appartient \`a 
un arc $\Gamma_j^{\leq t_j}$, $j\not=j_{2r}$, et
ce point, qui est un point simple de $\Gamma_j$, 
n'est pas l'extr\'emit\'e  $\Gamma_j^{t_j}$. 
Utilisant
le point {\rm(iv)} qui vient juste d'\^etre 
prouv\'e, on en d\'eduit que $j=j_{2r+1}$ et que 
{\rm(v)}
est vrai.

\proof[Preuve de (iii) et de (vi)] Fixons
$T=(t_j)_{1\leq j\leq n}\in{\cal T}_{2r+1}$. On 
vient de d\'emontrer {\rm(v)} et on sait donc que
$t_{j_{2r+1}}=t^*_{j_{2r+1}}>-K$. Choisissons
$t'_{j_{2r+1}}\in]-K,t^*_{j_{2r+1}}[$ assez proche de
$t^*_{j_{2r+1}}$ pour que l'arc 
$\Gamma_{j_{2r+1}}\vert_{[t'_{j_{2r+1}},t^*_{j_{2r+1}}[}$
ne rencontre aucun arc $\Gamma_j$, $j\not=j_{2r+1}$. La configuration
$T'$ obtenue \`a partir de $T$ en rempla\c cant 
$t^*_{j_{2r+1}}$ par $t'_{j_{2r+1}}$ n'\'etant 
plus dans
${\cal T}_{2r}$, il existe donc $j\not=n$ qui 
n'est pas dans $\Sigma(T')$. Il existe donc
$j_{2r+2}\not=j_{2r+1}$ tel que
$\Gamma_{j_{2r+1}}^{t^*_{j_{2r+1}}}\in\Gamma_{j_{2r+2}}^{\leq 
t_{j_{2r+2}}}$. L'entier $j_{2r+2}$
est le seul entier $j\not =j_{2r+1}$ tel que
$\Gamma_{j_{2r+1}}^{t^*_{j_{2r+1}}}\in\Gamma_j$. 
Nous venons de d\'emontrer {\rm(iii)} et {\rm 
(vi)}.

\proof[Preuve de l'in\'egalit\'e 
$j_{2r+2}<j_{2r+1}$] Nous gardons les 
notations
du paragraphe pr\'ec\'edent. Uti\-li\-sant de 
nouveau le point {\rm(v)} qui vient juste 
d'\^etre prouv\'e, on
en d\'eduit que
$j_{2r+1}\in\Sigma(T')$ et donc que l'entier 
$j\not=n$ qui n'appartient pas \`a $\Sigma(T')$ 
est
$j_{2r+2}$. En particulier $j_{2r+2}\not=n$. 
Puisque $j_{2r+2}\not\in\Sigma(T')$, l'arc
$\Gamma_{j_{2r+2}}^{\leq t_{j_{2r+2}}}$ ne 
rencontre aucun arc autre que lui m\^eme, que
$\Gamma_{j_{2r+1}}^{\leq t^*_{j_{2r+1}}}$ et que 
$\Gamma_{j_{2r+2}-1}^{\leq t_{j_{2r+2}-1}}$. On en
d\'eduit que $j_{2r+2}\not=j_{s}$, si $1\leq s \leq 2r+1$. En effet tout arc
$\Gamma_{j_s}^{\leq t_{j_s}}$ rencontre 
$\Gamma_{j_{s-1}}^{\leq t_{j_{s-1}}}$ et on sait 
que
$j_{s-1}\not\in\{j_s,j_s-1,j_{2r+1}\}$, puisque 
$j_{s-1}>j_s\geq j_{2r+1}$. Il  reste \`a montrer 
que
$j_{2r+2}\not\in
\{j_{2s}+1,\dots,j_{2s-1}-1\}$ si $s\leq r$. 
M\^eme dans le cas o\`u $s=r$, on sait que
$\Gamma_{j_{2s}}^{t_{j_{2s}}}=\Gamma_{j_{2s}}^{t^*_{j_{2s}}}$ appartient \`a
$\Gamma_{j_{2s+1}}^{<t_{j_{2s+1}}}=\Gamma_{j_{2s+1}}^{<t^*_{j_{2s+1}}}$.
On sait d'autre part que 
$\Gamma_{j_{2s-1}}^{t_{j_{2s-1}}}=\Gamma_{j_{2s-1}}^{t^*_{j_{2s-1}}}$
appartient \`a
$\Gamma_{j_{2s}}^{\leq t^*_{j_{2s}}}$. D'apr\`es 
la \fullref{prop4.4}, on en d\'eduit que
  $\Gamma_{j_{2s-1}}^{t^*_{j_{2s-1}}}\in\Gamma_{j_{2s}}^{<t^*_{j_{2s}}}$. 
Ceci implique que
$\Gamma_{j_{2s-1}}^{<t^*_{j_{2s-1}}}\cup\Gamma_{j_{2s}}^{<t^*_{j_{2s}}}$ 
est connexe et s\'epare
dans $\D$ tout arc $\Gamma_{j'}^{\leq -K}$, $j'\in
\{j_{2s}+1,\dots,j_{2s-1}-1\}$, de tout arc 
$\Gamma_{j''}^{\leq -K}$, $j''\not\in
\{j_{2s},\dots,j_{2s-1}\}$. Puisque, par 
hypoth\`ese, les arcs $\Gamma_j^{<t_j}$ sont 
disjoints deux
\`a deux, 
$\Gamma_{j_{2s-1}}^{<t^*_{j_{2s-1}}}\cup\Gamma_{j_{2s}}^{<t^*_{j_{2s}}}$
s\'epare
$\Gamma_{j'}^{<t_{j'}}$ de 
$\Gamma_{j''}^{<t_{j''}}$. On en d\'eduit que 
$\Gamma_{j'}^{\leq
t_{j'}}\cap\Gamma_{j''}^{\leq 
t_{j''}}=\emptyset$. En effet si ces deux arcs 
s'intersectaient, ils
devraient le faire sur
$\Gamma_{j_{2s-1}}^{<t^*_{j_{2s-1}}}\cup\Gamma_{j_{2s}}^{<t^*_{j_{2s}}}$, 
ce qui contredirait la
\fullref{prop4.4}. Puisque $j_{2r+1}\not\in\{j_{2s},\dots,j_{2s-1}\}$ et puisque
$\Gamma_{j_{2r+1}}^{\leq 
t_{j_{2r+1}}}\cap\Gamma_{j_{2r+2}}^{\leq 
t_{j_{2r+2}}}\not=\emptyset$, on
en d\'eduit que $j_{2r+2}\not\in\{j_{2s}+1,\dots,j_{2s-1}-1\}$.
\endproof

\begin{figure}[ht!]\small\vspace{2mm}
\labellist
\pinlabel $\alpha_3$ [br] at 47 155
\pinlabel $\alpha_3$ [br] at 284 154
\pinlabel $\alpha_1$ [tl] at 115 7
\pinlabel $\alpha_1$ [tl] at 351 7
\pinlabel $\alpha_2$ [lb] at 155 115
\pinlabel $\alpha_2$ [lb] at 392 115
\pinlabel $\alpha_4$ [r] at 8 47
\pinlabel $\alpha_4$ [r] at 244 47
\pinlabel $\omega_4$ [l] at 155 47
\pinlabel $\omega_3$ [t] at 47 7
\pinlabel $\omega_2$ [br] at 8 114
\pinlabel $\omega_1$ [bl] at 115 154
\endlabellist
\cl{\psfig{file=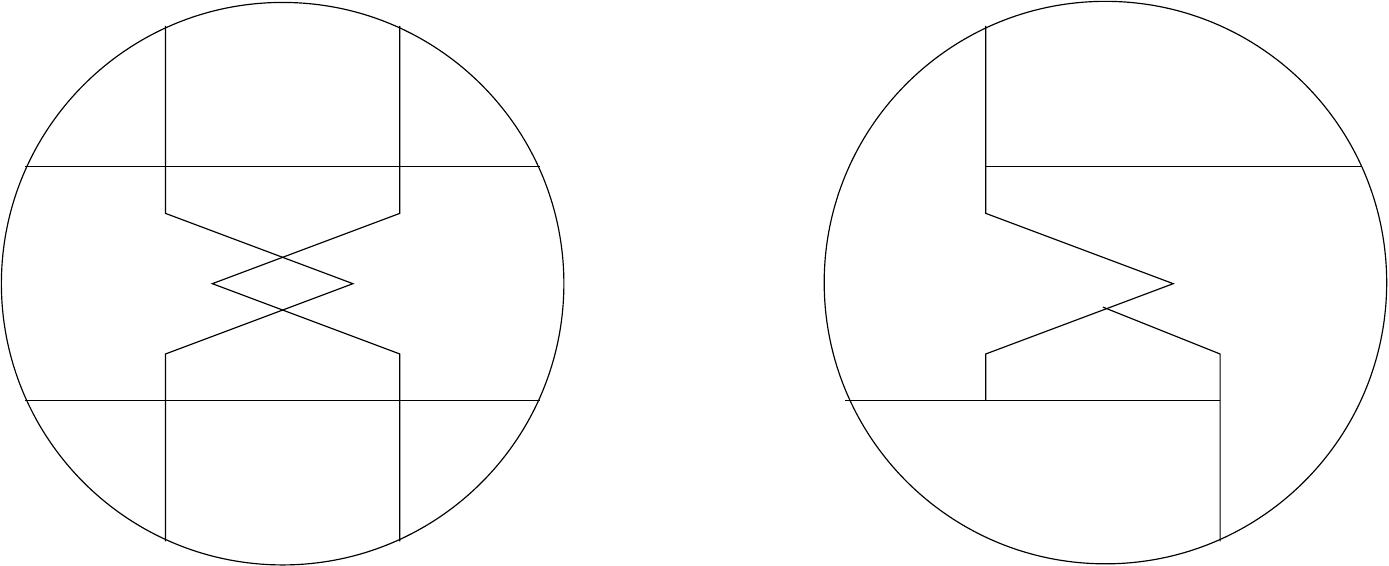,scale=70}}
\caption{}\label{fig7}
\end{figure}

Si $T=(t_i)_{i\in\Z/n\Z}\in{\cal T}_*$ est une 
configuration d'ordre $n$, on sait que $\vert
t_i\vert< K$, pour tout $i\in\Z/n\Z$, d'apr\`es 
le \fullref{lem5.2}. L'ensemble des configurations
$T\in{\cal T}_*$ d'ordre
$n$, qui est ferm\'e, est donc compact. 
Consid\'erons l'ordre lexicographique d\'efini 
sur $\cal T$ par
un ordre quelquonque sur
$\Z/n\Z$. On peut alors trouver une configuration 
maximale $T^0=(t_i^0)_{i\in\Z/n\Z}$ pour
cet ordre dans l'ensemble des configurations de ${\cal T}_*$
qui sont d'ordre $n$. On v\'erifie facilement,
utilisant les arguments vus dans la preuve pr\'ec\'edente, que pour tout
$i$, il existe $i'\not=i$ unique tel que 
$\Gamma_i^{t^0_i}\in \Gamma_{i'}$ et que
$\Gamma_i^{t^0_i}\in\Gamma_{i'}^{<t^0_{i'}}$. On 
se fixera dor\'enavant une telle configuration
et on posera $\Gamma_i^-=\Gamma_i^{\leq t_i^0}$ et $\sigma_i^0=\Gamma^{t_i^0}$.
On d\'efinira \'egalement
$\Gamma_i^+=\Gamma_i^{\geq K+1/2}$. Sur la \fullref{fig7} nous donnons un
exemple avec $n=4$, o\`u l'ordre sur $\Z/4\Z$ v\'erifie $1<2<3<4$. Il n'y a
qu'une configuration d'ordre $n=4$ maximale possible dans cet exemple.

\begin{prop}\label{prop5.5}
Nous avons les propri\'et\'es suivantes:

\bei

\item $\Gamma_i^+\cap \Gamma_{i'}^-=\Gamma_i^+\cap
\Gamma_{i'}^+=\emptyset$, si
$i'\not=i$;

\item $\Gamma_i^+\cap \Gamma_{i}^-=\emptyset$;

\item $(\Gamma_i^-\setminus\{\sigma_i^0\})\cap
(\Gamma_{i'}^-\setminus\{\sigma_{i'}^0\})=\emptyset$, si
$i'\not=i$;
\item  $\Gamma_i^{\leq -K}\cap
\Gamma_{i'}^-=\emptyset$, si $i'\not=i$;

\item l'ensemble 
$\Gamma_i^-\setminus 
\bigcup_{i'\not=i}\Gamma_{i'}^-$ est inclus dans
une composante connexe de ${\bf D}\setminus \bigcup_{i'\not=i}\Gamma_{i'}^-$;

\item pour tout $i\in\Z/n\Z$, 
il existe $i'\not\in\{i,i-1\}$ tel que
$\Gamma_i^-\cap \Gamma_{i'}^-\not=\emptyset$;

\item si
$\Gamma_i^-\cap\Gamma_{i'}^-\not=\emptyset$, avec 
$i'\not\in\{i,i-1\}$, alors la composante connexe
de
${\bf D}\setminus(\Gamma_i^-\cup\Gamma_{i'})$ qui 
contient $\Gamma_{i-1}^{\leq -K}$ est diff\'erente
de celle qui contient $\Gamma_{i-1}^+$.\eni
\end{prop}

\begin{figure}[ht!]\small\vspace{2mm}
\labellist
\pinlabel $\alpha_3$ [br] at 47 155
\pinlabel $\alpha_1$ [tl] at 115 7
\pinlabel $\alpha_2$ [lb] at 155 115
\pinlabel $\alpha_4$ [r] at 8 47
\endlabellist
\cl{\psfig{file=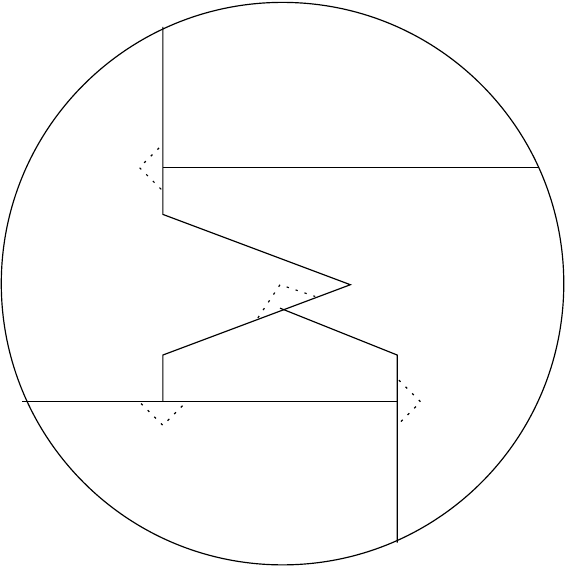,scale=70}}
\caption{}\label{fig8}
\end{figure}

\demo Le seul 
point qui n'est pas imm\'ediat ou qui n'a pas
\'et\'e mentionn\'e pr\'ec\'edemment est le 
cinqui\`eme point. Il n'est cependant pas tr\`es 
difficile \`a
prouver. Si le point
$\Gamma_i^t$,
$t<t_i^0$, n'appartient \`a aucun 
$\Gamma_{i'}^-$, $i'\not=i$, il n'y a en effet 
aucune difficult\'e \`a
perturber l'arc $\Gamma_i$ sur $[-K,t]$ pour 
qu'il aboutisse en $\Gamma_i^t$ tout en
\'evitant $\bigcup_{i'\not=i}\Gamma_{i'}^{-}$. En 
effet on peut contourner chaque point de
$\Gamma_i\vert_{[-K,t]}$ qui est sur un $\Gamma_{i'}^-$, $i'\not=i$,
au voisinage de ce point (voir \fullref{fig8}).\endproof

\section{Construction du graphe mod\`ele}\label{sec6}

Nous allons construire dans cette section un graphe $G$ qui sera le
mod\`ele unidimensionnel de notre d\'ecomposition 
en briques finale. Rappelons que l'extr\'emit\'e 
$\sigma_i^0$ de
$\Gamma_i^{-}$  appartient \`a un unique arc 
$\Gamma_{i'}^{-}$, $i'\not =i$, que
c'est un point simple de $\Gamma_i$ et de 
$\Gamma_{i'}$, qu'il n'appartient pas \`a $O_i$ 
et que
$\sigma_i^0\not=\sigma_{i'}^0$. On
appelera {\it singularit\'e} de
$\Gamma_i^{-}$ tout point  qui appartient \`a 
$O_i$ ou qui est l'extr\'emit\'e d'un arc 
$\Gamma_{i'}^-$.
Les singularit\'es \'etant des points simples de $\Gamma_i$, on peut les
param\'etrer sous la forme
$\sigma_i^l=\Gamma_{i}(t_i^l)$, o\`u 
$(t_i^l)_{l\leq 0}$ est une suite strictement 
croissante
uniquement d\'etermin\'ee. On se fixe alors $L$ 
suffisament grand pour que $t_i^{-L}\leq -K$. Les
singularit\'es $\sigma_i^l$, $l\leq -L$, sont 
alors n\'ecessairement des points de l'orbite 
$O_i$.
Remarquons qu'on ne peut pas trouver $i_0$ et $i_1$ tels que
$$\sigma_{i_0}^{-1}=\sigma_{i_1}^0,\enskip\sigma_{i_0}^{0}=\sigma_{i_1}^{-1}.$$ 
En effet, supposons que $i_0>i_1$ (pour l'ordre
sur $\Z/n\Z$ qui nous a permi de d\'efinir $T^0$). Si
$t>t_i^0$ est proche de
$t_i^0$, la configuration $T'=(t_i)_{i\in \Z/n\Z}$ o\`u $t'_{i_0}=t$,
$t'_{i_1}=t_{i_1}^{-1}$ et $t'_i=t_i^0$ si $i\not\in\{i_0,i_1\}$ v\'erifie 
les conditions
{\bf (A)} et {\bf (B)} et est plus grande que 
$T^0$ pour l'ordre lexicographique, ce qui
contredit la maximalit\'e de $T^0$.

D\'efinissons maintenant une suite
$(t^l)_{l\geq 1}$ en posant
$$\begin{cases}\,
   t^1=K+1/2&\\
\,t^l=K+l-1& \text {si $l\geq 2$}
\end{cases}
$$
et d\'efinissons les singularit\'es
$\sigma_i^l=\Gamma_i(t^l)$ de $\Gamma_i^+$. La 
\fullref{prop3.1} nous permet d'affirmer l'existence
d'une suite de disques de s\'ecurit\'e 
$(V_i^l)_{l\not =0}$ adapt\'ee \`a notre nouvelle 
indexation. Plus
pr\'ecis\'ement:

\begin{lem}\label{lem6.1}
On peut construire, pour tout $i\in\Z/n\Z$, une suite
$(V_i^l)_{l\not=0}$ de disques ouverts telle que:

\bei
\item chaque $V_i^l$ contient le segment
$\Gamma_{i}\vert_{[t^l_i,t^{l+1}_i]}$;
\item $V_i^l\cap 
V_{i'}^{l'}=\emptyset$, si $i\not=i'$ et $l>0$;
\item $V_i^l\cap V_{i}^{l'}=\emptyset$, si  $l'<0<l$;
\item $f(V_i^l)\cap V_i^{l'}=\emptyset$, si $l'<l$;
\item la suite $(\overline 
V_i^l)_{l>0}$ converge vers $\{\omega_i\}$;
\item la suite $(\overline 
V_i^l)_{l<0}$ converge vers $\{\alpha_i\}$.
\eni\end{lem}

On d\'efinit maintenant deux suites {\it d'arcs
d'attache}
$(\alpha_i^l)_{l<0}$ et $(\alpha_i^l)_{l>0}$ de la fa\c con suivante:

\item $\alpha_i^{-1}=\Gamma_i\vert_{[t_i^{-1},t^0_i]}$;

\item $\alpha_i^{1}=\Gamma_i\vert_{[t_i^{1},t^2_i]}$;

\item si $l<-1$, alors 
$\alpha_i^{l}=\Gamma_i\vert_{[t_i^{l},t'{}_i^{l}]}$, 
en notant
$t'{}_i^{l}$ le premier instant o\`u
$\Gamma_i\vert_{[t_i^{l},t_i^{l+1}]}$ rencontre $\bigcup_{l< m<0} \alpha_i^m$;

\item si $l>1$, alors 
$\alpha_i^{l}=\Gamma_i\vert_{[t'{}_i^{l},t_i^{l+1}]}$, 
en notant
$t'{}_i^{l}$ le dernier instant o\`u
$\Gamma_i\vert_{[t_i^{l},t_i^{l+1}]}$ rencontre $\bigcup_{0< m <l} \alpha_i^m$.

\noindent Remarquons que $\alpha_i^{l}$ est un 
segment (i.e.\ ne se r\'eduit pas \`a un point)
puisque  les singularit\'es sont des points simples de $\Gamma_i$.
Remarquons \'egalement que
$\alpha_i^{l}\subset V_i^l$.  Toujours gr\^ace \`a la \fullref{prop4.4}, on
sait que pour tout
$l<-1$, le point $\Gamma_i(t'{}_i^l)$ appartient \`a un unique arc
$\alpha_i^m$, $l<m<0$; qu'il
est diff\'erent de l'extr\'emit\'e 
$\Gamma_i(t'{}_i^m)$ et qu'il ne peut \^etre 
\'egal \`a l'extr\'emit\'e
$\sigma_i^m$ que dans le cas o\`u $m=l-1$ et 
$t'^l_i=t_i^{l-1}$. Nous avons un r\'esultat 
analogue pour les $\alpha_i^l$,
$l>0$. Parmi les singularit\'es, nous porterons 
particuli\`erement attention \`a la singularit\'e
$\sigma_{i}^*=\Gamma_i[(t_i^0)]$, o\`u $[t]$ d\'esigne la partie
enti\`ere d'un r\'eel $t$. Il s'agit de la
singularit\'e de $\Gamma_i^{-}$ la plus ``proche'' de l'extr\'emit\'e
$\sigma_{i}^0$ parmi les singularit\'es qui sont sur $O_i$.
Si on \'ecrit $\sigma_i^*=\sigma_i^{l_i}$, on peut remarquer que
$\alpha_i^l=\Gamma_{i}\vert_{[t_i^l,t_i^{l+1}]}$ si $l_i\leq l<0$.

\begin{figure}[ht!]\small\vspace{2mm}
\labellist
\pinlabel $\alpha_3$ [b] at 89 216
\pinlabel $\alpha_3$ [b] at 386 214
\pinlabel $\alpha_1$ [tl] at 163 15
\pinlabel $\alpha_1$ [tl] at 460 16
\pinlabel $\alpha_2$ [lb] at 204 162
\pinlabel $\alpha_2$ [lb] at 499 162
\pinlabel $\alpha_4$ [r] at 1 95
\pinlabel $\alpha_4$ [r] at 299 95
\pinlabel $\omega_4$ [l] at 217 95
\pinlabel $\omega_4$ [l] at 512 95
\pinlabel $\omega_3$ [t] at 89 1
\pinlabel $\omega_3$ [t] at 386 3
\pinlabel $\omega_2$ [br] at 15 162
\pinlabel $\omega_2$ [br] at 313 162
\pinlabel $\omega_1$ [bl] at 163 204
\pinlabel $\omega_1$ [bl] at 460 201
\hair1pt
\pinlabel $\sigma^0_1$ [bl] at 90 127
\pinlabel $\sigma^0_2$ [bl] at 90 165
\pinlabel $\sigma^0_3$ [br] at 86 97
\pinlabel $\sigma^0_4$ [tr] at 121 94
\endlabellist
\cl{\psfig{file=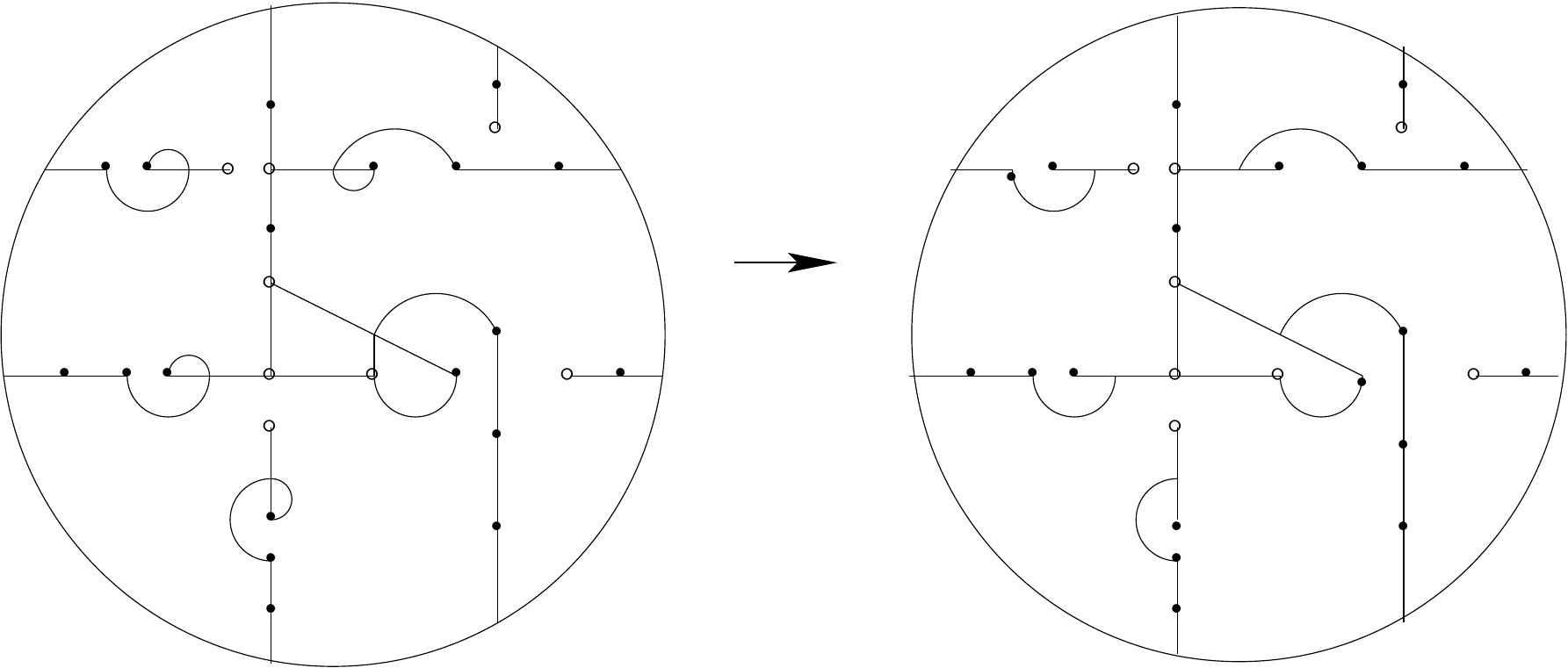,scale=68}}
\caption{}\label{fig9}
\end{figure}

On construit ainsi des arbres localement finis 
$A_i^-=\bigcup_{l<0} \alpha_i^l$ et 
$A_i^+=\bigcup_{l>0}
\alpha_i^l$: seul un nombre fini de $\alpha_i^l$
rencontrent une partie compacte de $\D$ donn\'ee. 
L'arbre $A_i^-$ est inclus dans $\Gamma_i^-$ et
contient toutes les singularit\'es de $\Gamma_i^-$. C'est sur le graphe
$G=\bigcup_{i\in\Z/n\Z}(A_i^-\cup A_i^+)$ que nous allons modeler notre
d\'e\-com\-po\-si\-tion. Nous illustrons la construction du graphe $G$ sur la
\fullref{fig9}. Les singularit\'es qui correspondent \`a des points de
l'orbite sont repr\'esent\'ees par des points noirs, celles
qui sont des extr\'emit\'es sont repr\'esent\'ees par des points blancs.

On d\'eduit facilement de la propri\'et\'e 5.5 et 
des inclusions $A_i^-\subset\Gamma_i^-$ et
$A_i^+\subset\Gamma_i^+$ le r\'esultat qui suit:

\eject
\begin{prop}\label{prop6.2}
Nous avons les propri\'et\'es suivantes:

\bei
\item $A_i^+\cap A_{i'}^-=\emptyset$, pour tous $i'$ et $i$
dans $\Z/n\Z$;
\item $A_i^+\cap A_{i'}^+=\emptyset$, si $i\not=i'$;
\item $(A_i^-\setminus\{\sigma_i^0\})\cap
(A_{i'}^-\setminus\{\sigma_{i'}^0\})=\emptyset$, si
$i\not=i'$;
\item  $\alpha_i^l\cap 
A_{i'}^-=\emptyset$, si $i\not=i'$ et $l<-L$;
\item l'ensemble 
$A_i^-\setminus \bigcup_{i'\not=i}A_{i'}^-$ est 
inclus dans une
composante connexe de\break ${\bf D}\setminus \bigcup_{i'\not=i}A_{i'}^-$;
\item pour tout $i\in\Z/n\Z$, 
il existe $i'\not\in\{i,i-1\}$ tel que
$A_i^-\cap A_{i'}^-\not=\emptyset$;
\item si
$A_i^-\cap A_{i'}^-\not=\emptyset$, avec 
$i'\not\in\{i,i-1\}$, alors la composante connexe 
de\break
${\bf D}\setminus(A_i^-\cup A_{i'}^-)$ qui 
contient les $\alpha_{i-1}^l$, $l<-L$, est 
diff\'erente de
celle qui contient $A_{i-1}^+$.\eni
\end{prop}

\section{Construction d'une d\'ecomposition adapt\'ee}\label{sec7}

Nous construirons dans cette section une d\'ecomposition en briques
libre maximale model\'ee sur le graphe 
$G=\bigcup_{i\in\Z/n\Z}(A_i^-\cup A_i^+)$.
Nous utiliserons implicitement le th\'eor\`eme de Schoenflies, ou
plut\^ot la version g\'en\'eralis\'ee qu'en donne 
Homma \cite{Ho}. La construction se fera en 
plusieurs \'etapes. Commen\c cons par la plus
simple, en rappelant que la famille de disques de s\'ecurit\'e
$(V_i^l)_{i\in\Z/n\Z, l\not=0}$ a \'et\'e d\'efinie dans le \fullref{lem6.1}.

\medskip
{\bf Construction de la famille $(b'{}_i^l)_{i\in\Z/n\Z, l>0}$}

\begin{lem}\label{lem7.1}On 
peut construire une famille 
$(b'{}_i^l)_{i\in\Z/n\Z,l>0}$
de disques ferm\'es de $\D$ telle que:

\bee
\item
chaque $b'{}_i^{l}$ est libre et inclus dans $V_i^l$;

\item${\rm Int}(b'{}_i^{l})\cap{\rm
Int}(b'{}_i^{l'})=\emptyset$, si $l\not=l'$;

\item${\rm 
Int}(b'{}_i^{1})$ contient $\sigma_{i}^{1}$ et 
$\sigma_{i}^{2}$;

\item${\rm Int}(b'{}_i^{l})$ contient
$\sigma_{i}^{l+1}$, si $l\geq2$;

\item${\rm 
Int}(\bigcup_{m\leq l}b'{}_i^{l})$ contient 
$\bigcup_{m\leq l}\alpha_i^{l}$;

\item
$\alpha_i^l$, $l>1$, rencontre un unique disque 
$b'{}_i^m$,
$m<l$, que l'on note $b'{}_i^{m[l]}$;

\item
si $b'{}_i^m$ intersecte un arc
$\alpha_i^{l}=\Gamma_i\vert_{[t'{}_i^{l},t_i^{l+1}]}$, 
$l> m$, alors $b'{}_i^m\cap
\alpha_i^{l}$ est un sous-segment non trivial 
$\Gamma_i\vert_{[t'{}_i^{l},t''{}_i^{l}]}$ de
$\alpha_i^{l}$ et $\Gamma_i\vert_{[t'{}_i^{l},t''{}_i^{l}[}\subset {\rm
Int}(b'{}^m_i)$;

\item un disque 
$b'{}_i^{l}$, $l>1$, ne rencontre aucun disque 
$b'{}_i^m$, $m<l$, autre
que
$b'{}_i^{m[l]}$;

\item$b'{}_i^{l}\cap b'{}_i^{m[l]}=\partial
b'{}_i^{l}\cap
\partial b'{}_i^{m[l]}$ est un segment.\ene\end{lem}

\demo
Puisque 
$\Gamma_i\vert_{[K,K+1]}$ est un arc de
translation, le segment
$\alpha_i^1=\Gamma_i\vert_{[K+1/2,K+1]}$ est 
libre. Nous pouvons donc construire un disque 
ferm\'e
$b'{}_i^1\subset V_i^1$ qui est libre et dont l'int\'erieur contient
$\alpha_i^1$. Nous pouvons de plus supposer que $b'{}_i^1$ ne contient
aucune singularit\'e $\sigma_i^{l}$,
$l>2$. Soyons plus pr\'ecis. Il n'y a qu'un 
nombre fini d'arcs $\alpha_i^l$, $l>1$, qui 
rencontrent
$\alpha_i^1$. Nous pouvons supposer, d'une part que ce sont les seuls
que rencontre $b'{}_i^1$,
d'autre part que
$b'{}_i^1$ v\'erifie la condition {\rm(vii)} du \fullref{lem7.1}.

Supposons maintenant que les disques
$b'{}_i^m$ ont \'et\'e construits jusqu'\`a l'ordre
$l$ et qu'ils v\'erifient les conditions du \fullref{lem7.1}. Ceci a un sens
si on remplace la condition {\rm(vi)} par la condition suivante:

{\rm(vi$'$)}\qua {\sl un arc d'attache
$\alpha_i^{l'}$, $l'>1$, rencontre au plus un disque $b'{}_i^m$,
  $m\leq \min(l, l'-1)$.}

\begin{figure}[ht!]\small\vspace{2mm}
\labellist
\pinlabel $\alpha_3$ [b] at 89 216
\pinlabel $\alpha_1$ [tl] at 163 15
\pinlabel $\alpha_2$ [lb] at 204 162
\pinlabel $\alpha_4$ [r] at 1 95
\pinlabel $\omega_4$ [l] at 217 95
\pinlabel $\omega_3$ [t] at 89 1
\pinlabel $\omega_2$ [br] at 15 162
\pinlabel $\omega_1$ [bl] at 163 204
\pinlabel $b'^1_2$ [b] at 60 165
\pinlabel $b'^1_1$ [r] at 160 183
\pinlabel $b'^1_3$ [l] at 91 61
\pinlabel $b'^1_4$ [b] at 194 97
\pinlabel $b'^2_3$ [r] at 71 45
\pinlabel $b'^3_3$ [r] at 84 19
\endlabellist
\cl{\psfig{file=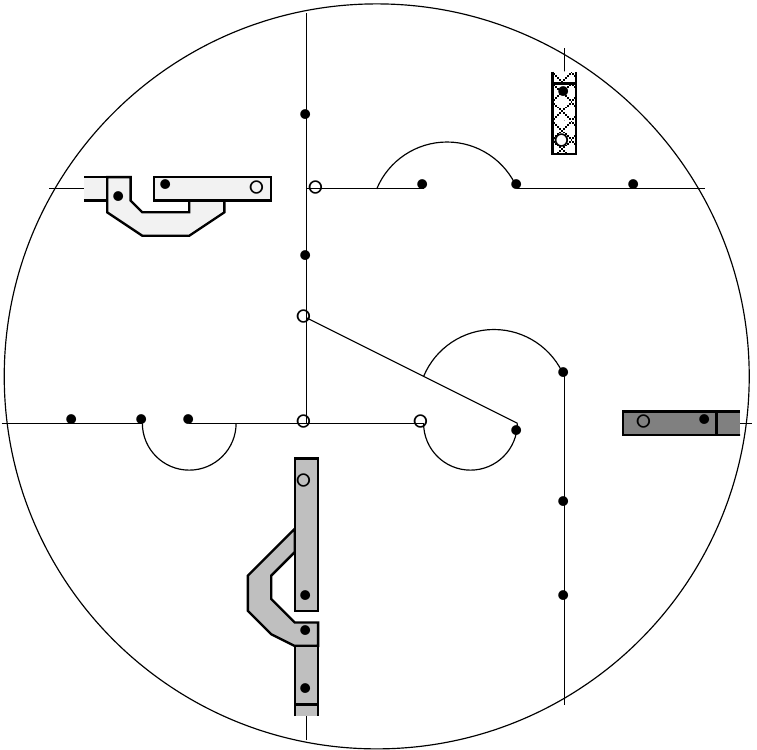,scale=80}}
\caption{}\label{fig10}
\end{figure}

Consid\'erons l'arc d'attache
$\alpha_i^{l+1}=\Gamma_i\vert_{[t'{}_i^{l+1},t_i^{l+2}]}$. Nous savons que
$\Gamma_i(t'{}_i^{l+1})$ appartient \`a un arc 
$\alpha_i^m$, $m\leq l$, et donc \`a un disque 
$b_i^{m'}$,
$m'\leq l$ d'apr\`es {\rm(v)}. Ce disque est unique d'apr\`es {\rm
(vi')}, on le note $b_i^{m[l+1]}$. D'apr\`es {\rm(viii)}, nous savons que
$\alpha_i^{l+1}\cap 
b_i^{m[l+1]}=\Gamma_i\vert_{[t'{}_i^{l+1},t''{}_i^{l+1}]}$, 
o\`u
$t'{}_i^{l+1}<t''{}_i^{l+1}<t_i^{l+2}$. Le segment
$\Gamma_i\vert_{[t''{}_i^{l+1},t{}_i^{l+2}]}$ est 
libre, puisque $\Gamma_i\vert_{[t_i^{l+1},
t_i^{l+2}]}$ est un arc de translation. Il est 
\'egalement inclus dans $V_i^{l+1}$. On peut donc
construire un disque ferm\'e $b'{}_i^{l+1}\subset 
V_i^{l+1}$ qui est libre, dont l'int\'erieur 
contient
$\Gamma_i\vert_{]t''{}_i^{l+1}, t_i^{l+2}]}$, et 
qui ne rencontre qu'un seul disque
$b'{}_i^m$, $m\leq l$, \`a savoir 
$b'{}_i^{m[l+1]}$. On peut \'egalement supposer 
que
$b'{}_i^{l+1}\cap b'{}_i^{m[l+1]}=\partial b'{}_i^{l+1}\cap \partial
b'{}_i^{m[l+1]}$ est un
segment dont
$\Gamma_i(t''{}_i^{l+1})$ n'est pas une extr\'emit\'e. Tout arc d'attache
$\alpha_i^{l'}$, $l'>l+1$, qui rencontre 
$\Gamma_i\vert_{[t''{}_i^{l+1}, t_i^{l+2}]}$ est 
disjoint de tous
les disques $b'{}_i^m$, $m\leq l$, d'apr\`es {\rm
(vii)}. On peut donc supposer que
$b'{}_i^{l+1}$ v\'erifie \'egalement la condition 
{\rm(vii)}. Remarquons maintenant que la
famille
$(b'{}_i^m)_{1\leq m\leq l+1}$ v\'erifie toutes les conditions du
lemme. Nous avons dessin\'e les disques $b'{}_{i}^l$,
$i\in\Z/n\Z$, $l>0$, sur la \fullref{fig10}.
\endproof

\noindent{\bf Construction de la famille $(b_i^l)_{i\in\Z/n\Z, l<0}$}

La construction sera similaire, \`a savoir un 
\'epaississement des arbres $A_i^-$, $i\in\Z/n\Z$.
La situation est cependant plus
d\'elicate car ces arbres ne sont pas disjoints. 
Rappelons que l'extr\'emit\'e $\sigma_i^0$ n'est
pas sur $O_i$ et que le segment $\Gamma_i\vert_{[[t_i^0],t_i^0]}$ qui joint la
singularit\'e $\sigma^*_i$ \`a $\sigma_i^0$ est 
donc libre. Il en est de m\^eme, bien s\^ur, du
sous-segment
$\alpha_i^{-1}$. Rappelons \'egalement que $\sigma_i^{0}$, qui appartient
\`a un unique arc
$\Gamma_{i_0}$, $i_0\not=i$, est sur $A_{i_0}^-$. Dans le cas o\`u
$\sigma_i^{-1}=\sigma_i^*$, on sait que
$\alpha_i^{-1}\setminus\sigma_i^0$ est disjoint 
de tout arc $\Gamma_{i'}^-$, $i'\not=i$, et
donc de $A_{i'}^-$.  Dans le cas contraire, 
$\sigma_i^{-1}$ appartient \`a un unique arc
$\Gamma_{i_1}$,
$i_1\not=i$,
et c'est
l'extr\'emit\'e
$\sigma_{i_1}^0$ de
$A_{i_1}^{-}$.  Dans ce cas ${\rm int}(\alpha_i^1)$ ne rencontre aucun arc
$\Gamma_{i'}^{-}$, $i'\not=i$. On peut donc 
construire une famille 
$(b'{}_i^{-1})_{i\in\Z/n\Z}$ de
disques ferm\'es libres de $\D$, v\'erifiant 
$\alpha_i^{-1}\subset b_i^{-1}\subset V_i^{-1}$,
dont les int\'erieurs sont disjoints deux \`a deux, et tels que
$$b_i^{-1}\cap b_{i'}^{-1}\not=\emptyset\enskip{\rm et}\enskip
i\not=i'\Longrightarrow
\sigma_i^{0}=\sigma_{i'}^{-1} \enskip{\rm ou} 
\enskip\sigma_{i}^{-1}=\sigma_{i'}^{0}.$$
On va \^etre plus pr\'ecis. Pour toute singularit\'e
$\sigma_{i}^l=\Gamma_i(t_i^l)$, $l<0$, choisissons
$\varepsilon_i^l>0$ tel que:

\bei
\item $t_i^l+\varepsilon_i^l<t'{}_i^l$;
\item
$\Gamma_i\vert_{[t_i^l-\varepsilon_i^l,t_i^l+\varepsilon_i^l]}\subset V_i^l$;
\item
$\Gamma_i\vert_{[t_i^l-\varepsilon_i^l,t_i^l
+\varepsilon_i^l]}\cap\left(\bigcup_{i'\not=i}\Gamma_{i'}\right)=\emptyset$,
si $\sigma_{i}^l\in O_{i}$;
\item
$\Gamma_i\vert_{[t_i^l-\varepsilon_i^l,t_i^l
+\varepsilon_i^l]}\cap\left(\bigcup_{i'\not=i}\Gamma_{i'}\right)=\{\sigma_i^l\}$,
si $\sigma_{i}^l\not\in O_{i}$.
\eni

\begin{lem}\label{lem7.2}On 
peut constuire une famille 
$(b_i^{-1})_{i\in\Z/n\Z}$ de
disques ferm\'es de $\D$ et une famille $(\eta_i^{-1})_{i\in\Z/n\Z}$ de r\'eels
$>0$, telles
que:

\bee
\item$b_i^{-1}$ contient $\alpha_i^{-1}$;
\item$b_i^{-1}$ est 
libre et inclus dans $V_i^{-1}$;
\item${\rm Int}(b_i^{-1})\cap {\rm
Int}(b_{i'}^{-1})=\emptyset$, si
$i\not=i'$;
\item si $b_i^{-1}$ intersecte un arc d'attache
$\alpha_i^l=\Gamma_{i}\vert_ {[t_i^l,t'{}_i^l]}$, 
$l\leq -2$, alors $b_i^1\cap \alpha_i^l$ est un
sous-segment $\Gamma_i\vert_{[t''{}_i^l,t'{}_i^l]}$ de
$\alpha_i^l$, avec $t_i^l<t''{}_i^l< t'{}_i^l$;
\item
$0<\eta_i^{-1}<\min(\varepsilon_i^{-1},\varepsilon_{i_0}^{l_0})$, 
o\`u on
\'ecrit $\sigma_i^0=\sigma_{i_0}^{l_0}$;
\item$b_i^{-1}\cap 
A_{i_0}^-=\Gamma_{i_0}\vert_{[t_{i_0}^{l_0}, 
t_{i}^{l_0}
+\eta_{i}^{-1}]}$;
\item si $\sigma_{i}^{-1}=\sigma_i^*$, alors
$b_i^{-1}\cap A_{i'}^-=b_i^{-1}\cap b_{i'}^{-1}=\emptyset$, si
$i'\not\in\{i,i_0\}$;
\item si $\sigma_{i}^{-1}=\sigma_i^*$,
alors $\sigma_{i}^{-1}\in{\rm Int}(b_i^{-1})$;
\item si 
$\sigma_{i}^{-1}= \sigma_{i_1}^0$ est une 
extr\'emit\'e,
  alors
$b_i^{-1}\cap A_{i'}^-=b_i^{-1}\cap b_{i'}^{-1}=\emptyset$, si
$i'\not\in\{i,i_0,i_1\}$;
\item si 
$\sigma_{i}^{-1}= \sigma_{i_1}^0$ est une 
extr\'emit\'e, alors
$\partial b_i^{-1}$ contient $\Gamma_{i}\vert_{[t_{i}^{-1}-\eta_i^{-1},
t_{i}^{-1} +\eta_{i}^{-1}]}$;
\item si
$\sigma_{i}^{-1}= \sigma_{i_1}^0$ est une 
extr\'emit\'e, alors
$$b_i^{-1}\cap b_{i_1}^{-1}=\partial
b_i^{-1}\cap
\partial 
b_{i_1}^{-1}=\Gamma_{i_1}\vert_{[t_{i_1}^{l_1}, 
t_{i_1}^{l_1} +\min(\eta_{i}^{-1}, 
\eta_{i_1}^{-1})
]}.$$\ene\end{lem}

\begin{figure}[ht!]\small\vspace{2mm}
\labellist
\pinlabel $\alpha_3$ [b] at 89 216
\pinlabel $\alpha_1$ [tl] at 163 15
\pinlabel $\alpha_2$ [lb] at 204 162
\pinlabel $\alpha_4$ [r] at 1 95
\pinlabel $\omega_4$ [l] at 217 95
\pinlabel $\omega_3$ [t] at 89 1
\pinlabel $\omega_2$ [br] at 15 162
\pinlabel $\omega_1$ [bl] at 163 204
\pinlabel $b^{-1}_1$ [bl] at 110 118
\pinlabel $b^{-1}_2$ [b] <4pt,0pt> at 101 165
\pinlabel $b^{-1}_3$ [r] at 85 115
\pinlabel $b^{-1}_4$ [t] <4pt,0pt> at 106 87
\endlabellist
\cl{\psfig{file=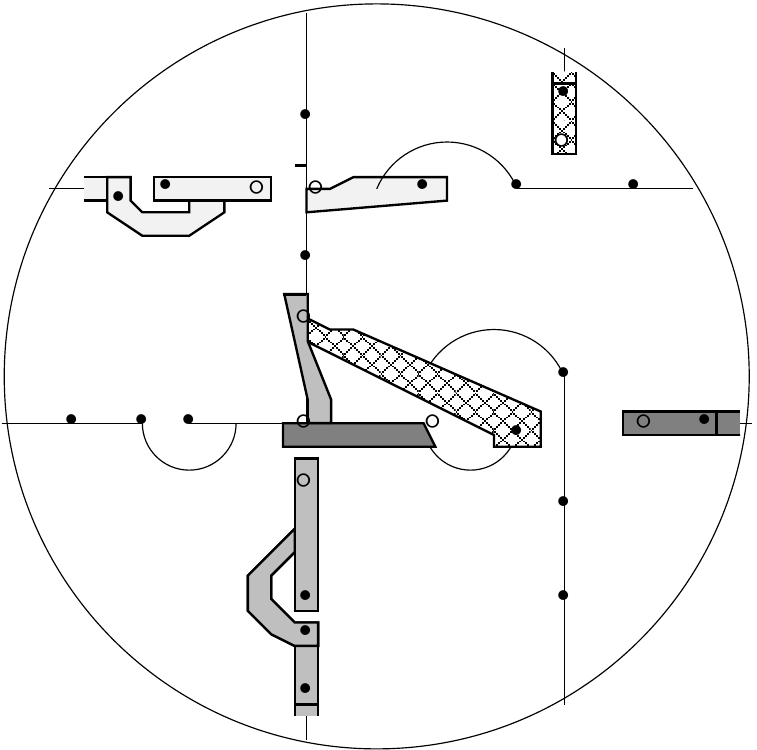,scale=80}}
\caption{}\label{fig11}
\end{figure}

Nous avons dessin\'e les disques $b_{i}^{-1}$,
$i\in\Z/n\Z$, sur la \fullref{fig11}.

 On peut alors \'etendre la famille
$(b_i^{-1})_{i\in \Z/n\Z}$:

\begin{lem}\label{lem7.3}On 
peut \'etendre la famille $(b_i^{-1})_{i\in 
\Z/n\Z}$ en une
famille
$(b_i^l)_{i\in\Z/n\Z, l<0}$ de disques ferm\'es 
de $\D$ et la famille $(\eta_i^{-1})_{i\in 
\Z/n\Z}$ en une famille
$(\eta_i^l)_{i\in\Z/n\Z, l<0}$ de r\'eels $>0$ telles que:
\bee
\item$0<\eta_i^{l}<\varepsilon_i^{l}$;
\item chaque 
$b_i^{l}$ est libre et inclus dans $V_i^l$;
\item $\sigma_{i}^{l}\in{\rm
Int}(b_i^{l})$, si $\sigma_{i}^l\in
O_i$;
\item$\sigma_{i}^{l}\in\partial b_i^{l}$, si
$\sigma_{i}^l\not\in
O_i$,
plus pr\'ecis\'ement, $\Gamma_{i}\vert_{[t_{i}^{l}-\eta_i^l,
t_{i}^l +\eta_{i}^l]}\subset\partial b_i^l$;
\item $\bigcup_{m\geq 
l}b_i^{m}$ contient $\bigcup_{m\geq 
l}\alpha_i^{m}$;
\item ${\rm Int}(b_i^{l})\cap {\rm
Int}(b_{i'}^{l'}) =\emptyset$, si
$(i,l)\not=(i',l')$;
\item ${\rm Int}(b_i^{l})\cap A_{i'}^-
=\emptyset$, si
$i\not=i'$;
\item un arc 
d'attache $\alpha_i^l$, $l<-1$, rencontre un 
unique disque $b_i^m$,
$m>l$, que l'on note $b_i^{m[l]}$;
\item si $b_i^m$ intersecte un arc d'attache
$\alpha_i^{l}=\Gamma_i\vert_{[t_i^{l},t'{}_i^{l}]}$, $l<m$, alors $b_i^m\cap
\alpha_i^{l}$ est un sous-segment 
$\Gamma_i\vert_{[t''{}_i^{l},t'{}_i^{l+1}]}$ de
$\alpha_i^{l}$, o\`u $t_i^l<t''{}_i^l<t_i^l$;
\item un disque 
$b_i^{l}$, $l<-1$, ne rencontre aucun disque 
$b_i^m$, $m>l$, autre que
$b_i^{m[l]}$;
\item $b_i^{l}\cap b_i^{m[l]}=\partial b_i^{l}\cap
\partial b_i^{m[l]}$ est un segment;
\item un disque 
$b_i^{l}$, $l<-1$, ne rencontre aucun disque 
$b_{i'}^{l'}$, $i'\not=i$, si
  $\sigma_{i}^l\in O_i$;
\item un disque 
$b_i^{l}$, $l<-1$, ne rencontre aucun disque 
$b_{i'}^{l'}$, $i'\not=i$,
autre que
$b_{i_1}^{-1}$, si $\sigma_{i}^l$ est l'extr\'emit\'e
$\sigma_{i_1}^0$, et on a
$b_i^{l}\cap b_{i_1}^{-1}=\partial b_i^{l}\cap
\partial b_{i_1}^{-1}=\Gamma_i\vert_{[t_i^l,t_i^l+\eta_{i_1}^{-1}]}$.
\ene\end{lem}

\demo On raisonne par r\'ecurrence, comme dans la preuve du
\fullref{lem7.1}. On suppose donc qu'une famille $(b_i^m)_{i\in\Z/n\Z,
m\geq l}$ a \'et\'e construite qui v\'erifie toutes les hypoth\`eses
du \fullref{lem7.3} en rempla\c cant {\rm(viii)} par:

\medskip{\rm(viii$'$)}\qua {\sl un arc d'attache $\alpha_i^{l'}$,
$l'<0$, rencontre au plus un disque $b_i^m$,
$m\geq \max(l, l'+1)$.}

\medskip Fixons $i\in\Z/n\Z$ et
consid\'erons l'arc d'attache 
$\alpha_i^{l-1}=\Gamma_i\vert_{[t_i^{l-1},t'{}_i^{l-1}]}$. 

Nous savons que
$\Gamma_i(t'{}_i^{l-1})$ appartient \`a un arc $\alpha_i^m$,
$m\geq l$, et donc \`a un disque $b_i^{m'}$,
$m'\geq l$, d'apr\`es {\rm(v)}. D'apr\`es {\rm(viii')}, ce disque est
le seul disque
$b_i^m$, $m\geq l$ que rencontre $\alpha_i^{l-1}$.  On le note
$b_i^{m[l-1]}$.
D'apr\`es {\rm(ix)}, nous savons que $\alpha_i^{l-1}\cap
b_i^{m[l-1]}=\Gamma_i\vert_{[t''{}_i^{l-1},t'{}_i^{l-1}]}$, o\`u
$t_i^{l-1}<t''{}_i^{l-1}<t'{}_i^{l-1}$.
  Le segment
$\Gamma_i\vert_{[t_i^{l-1},t''{}_i^{l-1}]}$ est 
libre et inclus dans $V_i^{l-1}$. Dans le cas o\`u
$\sigma_i^{l-1}\in O_i$, le segment 
$\Gamma_i\vert_{[t_i^{l-1},t''{}_i^{l-1}]}$ ne 
rencontre aucun
arbre
$A_{i'}^-$,
$i'\not=i$, et aucun disque
$b_{i'}^{m}$, $i'\not=i$,
$m\geq l$. Dans le cas o\`u $\sigma_i^{l-1}$ est 
l'extr\'emit\'e $\sigma_{i_1}^0$, ce segment
  ne rencontre aucun arbre
$A_{i'}^-$,
$i'\not\in \{i,i_1\}$, et n'intersecte 
$A_{i_1}^-$ qu'en son extr\'emit\'e 
$\sigma_{i_1}^0$. De
m\^eme, il
ne rencontre aucun disque $b_{i'}^m$, $i'\not=i$, 
$m\geq l$, autre que $b_{i_1}^{-1}$ et
on a
$b_{i_1}^{-1}\cap\Gamma_i\vert_{[t_i^{l-1},t''{}_i^{l-1}]}=
\Gamma_i\vert_{[t_i^{l-1},t_i^{l-1}+\eta_{i_1}^{-1}]}$. 
On peut donc construire une
famille de disques libres
$(b_i^{l-1})_{i\in\Z/n\Z}$, disjoints deux \`a deux, telle que

\bei
\item
$b_i^{l-1}$ est inclus dans
$V_i^{l-1}$,

\item $b_{i}^{l-1}$ contient 
$\Gamma_i\vert_{[t_i^{l-1},t''{}_i^{l-1}]}$;

\item $b_i^{l-1}$ ne rencontre
aucun disque $b_i^{m}$, $m\geq l$, autre que $b_{i}^{m[l-1]}$;

\item $b_i^{l-1}\cap 
b_{i}^{m[l-1]}=\partial b_i^{l-1}\cap \partial
b_{i}^{m[l-1]}$
est un segment dont $\Gamma(t''{}_i^{l-1})$ n'est pas une extr\'emit\'e;

\item $b_i^{l-1}$ ne rencontre aucun disque
$b_{i'}^{m'}$,  $i'\not=i$, $m'\geq l$, dans le 
cas o\`u $\sigma_i^{l-1}\in O_i$;

\item ${\rm Int}(b_i^{l-1})$ 
contient 
$\Gamma_i\vert_{[t_i^{l-1},t''{}_i^{l-1}[}$, dans
le cas o\`u $\sigma_i^{l-1}\in O_i$;

\item $b_i^{l-1}$  ne rencontre aucun
disque
$b_{i'}^{m'}$, $i'\not=i$, $m'\geq l$, autre que 
$b_{i_1}^{-1}$ dans le cas o\`u
$\sigma_i^{l-1}=\sigma_{i_1}^0$;

\item  $\partial b_i^{l-1}$ contient
$\Gamma_{i}\vert_{[t_{i}^{l+1}-\eta_i^{l-1}, t_{i}^{l-1} +\eta_{i}^{l-1}]}$ si
$\sigma_i^{l-1}\not\in O_i$.
\eni

Puisque par hypoth\`ese, tout arc attache $\eta_i^{l'}$, $l'<l-1$ qui rencontre
$\Gamma_i\vert_{[t_i^{l-1},t''{}_i^{l-1}]}$ est 
disjoint de toute brique $b_i^m$, $m\geq l$, on
peut supposer de plus que la condition {\rm(x)} est \'egalement remplie.

\begin{figure}[ht!]\small\vspace{2mm}
\labellist
\pinlabel $\alpha_3$ [b] at 89 216
\pinlabel $\alpha_1$ [tl] at 163 15
\pinlabel $\alpha_2$ [lb] at 204 162
\pinlabel $\alpha_4$ [r] at 1 95
\pinlabel $\omega_4$ [l] at 217 95
\pinlabel $\omega_3$ [t] at 89 1
\pinlabel $\omega_2$ [br] at 15 162
\pinlabel $\omega_1$ [bl] at 163 204
\pinlabel $b^{-1}_1$ [bl] at 106 118
\pinlabel $b^{-2}_1$ [t] at 133 79
\pinlabel $b^{-3}_1$ [b] at 150 122
\pinlabel $b^{-4}_1$ [l] at 166 79
\pinlabel $b^{-5}_1$ [l] at 166 52
\endlabellist
\cl{\psfig{file=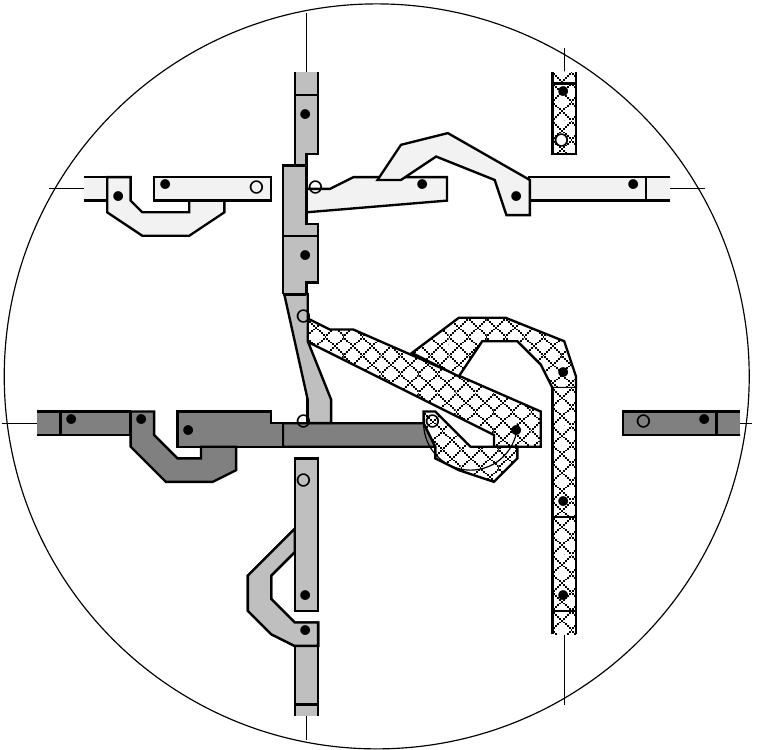,scale=80}}
\caption{}\label{fig12}
\end{figure}

On peut v\'erifier que les conditions du lemme 
sont maintenant v\'erifi\'ees jusqu'\`a l'ordre 
$l-1$.
Nous avons dessin\'e les disques $b_{i}^l$,
$i\in\Z/n\Z$, $l<0$, sur la \fullref{fig12}. \endproof

\medskip{\bf Remarque}\qua\'Ecrivons 
$\sigma_i^*=\sigma_i^{l_i}$. Puisque les $V_i^l$ 
peuvent \^etre
choisis arbitrairement proches des $\alpha_i^l$ 
et puisque $\Gamma_i\vert_{[[t_i^0], t_i^0]}$ est
libre, on peut supposer que
$\bigcup_{l_i\leq l<0}b_i^{l_i}$ (qui est un disque ferm\'e) est libre.

\medskip{\bf Construction de la famille $(b'{}_i^l)_{i\in\Z/n\Z, l<0}$}

On va modifier de fa\c con naturelle la famille 
$(b_i^l)_{i\in\Z/n\Z, l<0}$. On d\'efinit,
pour tout
$r<0$, une famille $(b_{i,r}^l)_{i\in\Z/n\Z, 
l<0}$ de disques libres par les relations
de r\'ecurrence suivantes, o\`u
$m[l]\in\{l+1,\dots, -1\}$ est d\'efini par l'assertion  {\rm(viii)}
du \fullref{lem7.2}:

{\sl\bei

\item $b_{i,-1}^l=b_i^l$, pour tout $l<0$;
\item $b_{i,r}^l=b_i^l$, pour tout $l<r$;

\item  si $b_i^r\cup
b_{i,r+1}^{m[r]}$ n'est pas libre, on pose
$$b_{i,r}^l=b_{i,r+1}^l, \enskip {\it pour \enskip tout }\enskip l<0;$$
 si $b_i^r\cup
b_{i,r+1}^{m[r]}$ est libre,  on pose
$$
\begin{cases}
\, b_{i,r}^r=b_i^r\cup
b_{i,r+1}^{m[r]}, &\\
\, b_{i,r}^l=b_i^r\cup b_{i,r+1}^{m[r]} ,& \text {si $l>r$ et
$b_{i,r+1}^l=b_{i,r+1}^{m[r]}$,} \\
\, b_{i,r}^l=b_{i,r+1}^l, & \text {si $l>r$ et
$b_{i,r+1}^l\not=b_{i,r+1}^{m[r]}$.}
\end{cases}
$$\eni
}

 Chaque disque $b_{i,r}^l$ est une 
r\'eunion de disques $b_i^m$, $m\geq l$, et 
contient
$b_i^l$. Puisque
$b_i^l$ est inclus dans le disque de s\'ecurit\'e
$V_i^l$ et puisque
$f(V_i^m)\cap V_i^l=\emptyset$, si $m>l$, on peut 
remarquer, dans la construction pr\'ec\'edente que
$f(b_i^r)\cap b_{i,r+1}^{m[r]}\not=\emptyset$, dans le cas o\`u  $b_i^r\cup
b_{i,r+1}^{m[r]}$ n'est
pas libre. On en d\'eduit:

\begin{lem}\label{lem7.4} 
Pour tout $r<-1$ et tout $l\in\{r,\dots,-2\}$, il 
existe
$m\in\{l+1,\dots,-1\}$ 
tel que  $b'{}^l_{i,r}=b'{}^m_{i,r}$ ou $f(b'{}^l_{i,r})\cap
b'{}^m_{i,r}\not=\emptyset$.
\end{lem}

Remarquons maintenant:
\begin{lem}\label{lem7.5}Deux cas sont possibles:
\qua {\rm (i)}\qua $f$ est r\'ecurrent;

\qua {\rm (ii)}\qua tout disque 
$b_{i,r}^l$ contient au plus $n-1$ disques 
$b_i^{m}$.
\end{lem}

\demo Chaque disque $b_i^m$, $m<0$, contient
$\sigma_i^m$. Parmi ces singularit\'es, il y en a 
au plus $n-1$ qui sont l'extr\'emit\'e d'un graphe
$A_{i'}^-$, $i'\not=i$. Ainsi un disque 
$b_{i,r}^l$ qui contient au moins $n$ disques 
$b_i^{m}$
contient deux points de l'orbite $O_i$ et 
rencontre donc l'un de ses it\'er\'es. Puisque 
$b_{i,r}^l$
est libre, $f$ est n\'ecessairement r\'ecurrent, 
d'apr\`es la \fullref{prop1.4}.\endproof

\begin{figure}[ht!]\small\vspace{2mm}
\labellist
\pinlabel $\alpha_3$ [b] at 89 216
\pinlabel $\alpha_1$ [tl] at 163 15
\pinlabel $\alpha_2$ [lb] at 204 162
\pinlabel $\alpha_4$ [r] at 1 95
\pinlabel $\omega_4$ [l] at 217 95
\pinlabel $\omega_3$ [t] at 89 1
\pinlabel $\omega_2$ [br] at 15 162
\pinlabel $\omega_1$ [bl] at 163 204
\pinlabel $b'{}^{-1}_1{=}b'{}^{-2}_1$ [t] <-5pt,1pt> at 133 79
\pinlabel $b'{}^{-3}_1$ [b] at 150 122
\pinlabel $b'{}^{-4}_1$ [l] at 166 79
\pinlabel $b'{}^{-5}_1$ [l] at 166 52
\endlabellist
\cl{\psfig{file=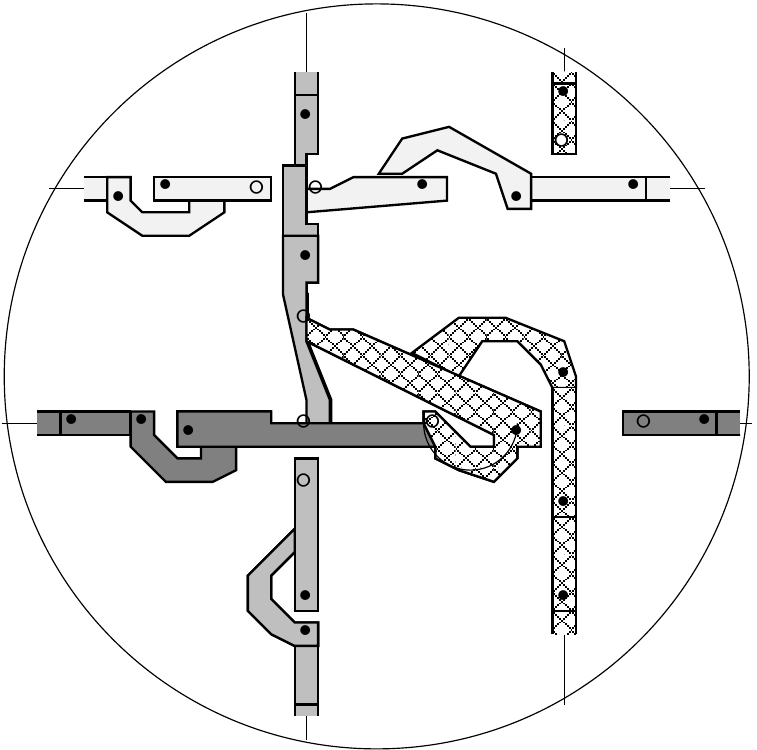,scale=80}}
\caption{}\label{fig13}
\end{figure}

Nous supposerons dor\'enavant que tout disque 
$b_{i,r}^l$ contient au plus $n-1$ disques 
$b_i^{m}$.
La suite $(b_{i,r}^l)_{r<0}$ est alors 
stationnaire. On d\'efinit donc une famille de 
disques
$(b'{}_i^l)_{i\in\Z/n\Z,l<0}$ en notant $b'{}_i^l$  la valeur de
$b_{i,r}^l$, pour $r$ assez petit.
Nous avons dessin\'e les disques $b'{}_{i}^l$,
$i\in\Z/n\Z$, $l<0$, sur la \fullref{fig13}.

La remarque qui suit le \fullref{lem7.3} implique le r\'esultat suivant:

\begin{lem}\label{lem7.6}Si 
on \'ecrit $\sigma_i^*=\sigma_i^{l_i}$, les
briques
$b'{}_i^{m}$, $l_i\leq m<0$, sont toutes 
\'egales. En particulier, $b'{}_i^{-1}$ contient
$\sigma_i^*$.\end{lem}

Gr\^ace au \fullref{lem7.4}, on sait que pour tout $l<-1$, il existe
$m\in\{l+1,\dots,-1\}$ tel que   $b'{}^l_{i}=b'{}^m_{i}$ ou $f(b'{}^l_{i})\cap
b'{}^m_{i}\not=\emptyset$. On en d\'eduit:

\begin{lem}\label{lem7.7}
Supposons que $l<-1$ et que $b'{}_i^{l}\not= 
b'{}_i^{-1}$.
Il existe alors une cha\^\i ne de disques choisis 
dans les $b'{}_i^{l}$, $l\leq m<0$, qui joint
$b'{}_i^{l}$ \`a $b'{}_i^{-1}$.\end{lem}

  Utilisant maintenant le fait que chaque 
$b'{}_i^l$, $l>0$, contient  $\sigma_i^{l+1}$
(qui est un point de $O_i$), le fait que  $b'{}_i^{-1}$ contient
$\sigma_i^*$ et le r\'esultat que l'on vient de noter, on obtient:

\begin{lem}\label{lem7.8}Soient $l$ et $l'$ deux entiers
non nuls tels que $l<l'$. Il existe une
cha\^\i ne de disques choisis dans les 
$b'{}_i^{m}$, $l\leq m\leq l'$, qui joint
$b'{}_i^{l}$ \`a $b'{}_i^{l'}$ dans les cas suivants:

\bei
\item $l<l'=-1$ et $b'{}_i^{l}\not= b'{}_i^{-1}$;

\item $l<0<l'$

\item $0<l<l'$.\eni\end{lem}

\noindent{\bf Construction de la famille $(b''{}_i^l)_{i\in\Z/n\Z, l\not=0}$}

\begin{figure}[ht!]\small\vspace{2mm}
\labellist
\pinlabel $\alpha_3$ [b] at 89 216
\pinlabel $\alpha_1$ [tl] at 163 15
\pinlabel $\alpha_2$ [lb] at 204 162
\pinlabel $\alpha_4$ [r] at 1 95
\pinlabel $\omega_4$ [l] at 217 95
\pinlabel $\omega_3$ [t] at 89 1
\pinlabel $\omega_2$ [br] at 15 162
\pinlabel $\omega_1$ [bl] at 163 204
\pinlabel $b''{}^{-1}_1\!{=}b''{}^{-2}_1$ [t] <-3pt,2pt> at 133 79
\pinlabel $b''{}^{-3}_1$ [l] at 165 128
\pinlabel $b''{}^{-4}_1$ [l] at 166 79
\pinlabel $b''{}^{-5}_1$ [l] at 166 52
\pinlabel* $b''{}^{-1}_2$ [t] at 145 154
\pinlabel $b''{}^{-1}_3$ [r] <2pt,-1pt> at 81 137
\pinlabel $b''{}^{-1}_4$ [b] at 56 107
\endlabellist
\cl{\psfig{file=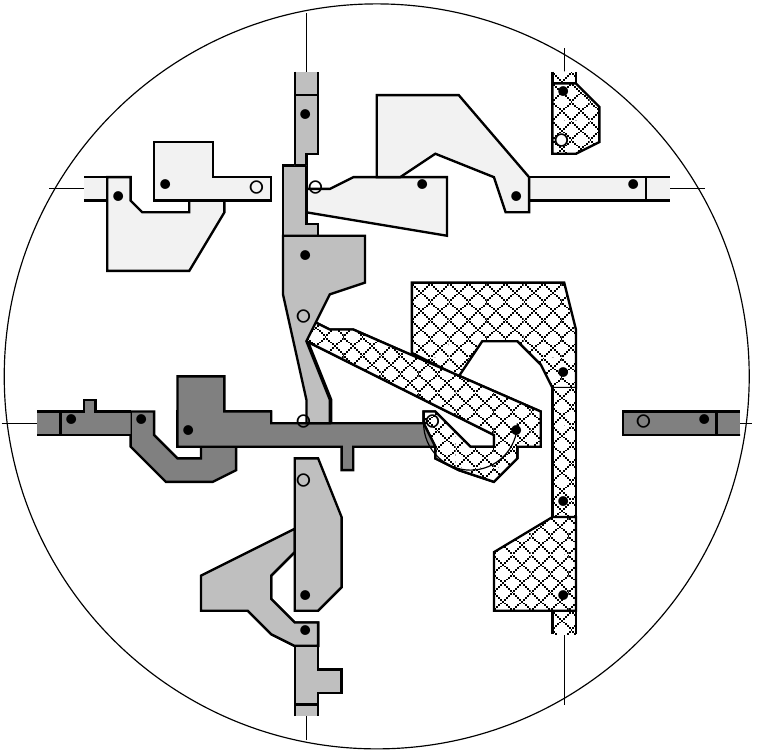,scale=80}}
\caption{}\label{fig14}
\end{figure}

On peut remarquer que l'ensemble $\bigcup_{i\in\Z/n\Z, l\not =0}\partial
b'{}_i^l$ est le squelette
d'une d\'e\-com\-po\-si\-tion en briques ${\cal 
D}'$ de ${\bf D}\setminus{\rm Fix}(f)$ dont les 
$b'{}_i^l$
sont des briques libres.  On peut d\'ecomposer les
briques de
${\cal D}'$ qui ne sont pas libres pour obtenir une d\'ecomposition libre de
$\D\setminus{\rm Fix}(f)$. On peut ensuite consid\`erer une
sous-d\'ecomposition libre maximale ${\cal 
D}''=(S'',A'',B'')$ de cette d\'ecomposition.
Cha\-que disque
$b'{}_i^l$ est alors contenu dans une brique
$b''{}_i^l$. Il faut noter qu'on peut
avoir
$b''{}_i^l=b''{}_{i'}^{l'}$ m\^eme si 
$b'{}_i^l\not=b'{}_{i'}^{l'}$ et m\^eme si
$i\not=i'$. Nous avons dessin\'e les disques $b''{}_{i}^l$,
$i\in\Z/n\Z$, $l\not=0$, sur la \fullref{fig14}.

\section{D\'emonstration du caract\`ere r\'ecurrent de $f$}\label{sec8}

Nous allons prouver dans cette section que $f$ 
est r\'ecurrent en construisant une cha\^\i ne 
ferm\'ee
de briques de ${\cal D}''$. En d'autres termes, nous allons trouver une
brique
$b''\in B''$ qui appartient \`a son futur strict, c'est-\`a-dire une
brique $b''$ telle que $b''\in b''_{>}$. Rassemblons dans une
proposition les propri\'et\'es
que nous allons utiliser. On peut supposer que la premi\`ere est
v\'erifi\'ee, d'apr\`es la \fullref{prop2.2}. Les autres se d\'eduisent
des r\'esultats de la section pr\'ec\'edente et du fait que
$b''{}_i^l$ contient $b'{}_i^l$.

\begin{prop}\label{prop8.1}Nous avons les
r\'esultats suivants:

\bee
\item pour tout $b''\in 
B''$, les ensembles $b''{}_{\geq}$, 
$b''{}_{\leq}$,
$b''{}_{>}$ et
$b''{}_{<}$ sont connexes;
\item$b''{}_i^{l'}\in
(b''{}_i^{l})_{\geq}$,  si $0<l<l'$, si
$l<l'=-1$ ou si $l<0<l'$;
\item$\sigma_i^{l+1}\in{\rm Int}( b''{}_i^{l})$
si $l>0$;
\item$\sigma_i^{l}\in b''{}_i^{l}$ si
  $l<0$;
\item$\sigma''{}_i^{l}\in{\rm Int}(
b''{}_i^{l})$, si
$\sigma_i^{l}\in O_i$ et $l<0$;
\item$\sigma_i^*\in{\rm Int}( b''{}_i^{-1})$.
\ene\end{prop}

Nous prouverons donc:

\begin{prop}\label{prop8.2}Il existe $i\in\Z/n\Z$
tel que $b''{}_i^{-1}\in(b''{}_i^{-1})_{>}$.\end{prop}

\demo Fixons
$i\in\Z/n\Z$. On sait, d'apr\`es la \fullref{prop6.2}, qu'il existe
$i'\not\in\{i, i+1\}$ tel que
$A_{i+1}^-\cap A_{i'}^-\not=\emptyset$ et que 
$A_{i+1}^-\cup A_{i'}^-$ est donc connexe.
Par construction des briques $b{}_i^l$, $b'{}_i^l$ et $b''{}_i^l$, on sait que
$$A_{i+1}^-\cup
A_{i'}^-\subset X=\left(\bigcup 
_{l<0}b''{}_{i+1}^l\right)\cup\left(\bigcup_{l<0}
b''{}_{i'}^l\right).$$ On sait de plus que $X$ 
est connexe comme r\'eunion de briques qui 
rencontrent
$A_{i+1}^-\cup A_{i'}^-$. Nous allons montrer que 
le futur $(b''{}_i^{-1})_{\geq}$
de $b''{}_i^{-1}$ contient au moins une brique 
qui est dans $X$. On sait, d'apr\`es la 
\fullref{prop6.2},
que les points $\sigma_i^{l}$, $l>0$, sont s\'epar\'es de $\sigma_i^*$
par $A_{i+1}^-\cup A_{i'}^-$.
Puisque
$(b''{}_i^{-1})_{\geq}$ est connexe, la partie
${\rm 
Int}\left((b''{}_i^{-1})_{\geq}\right)\subset{\bf 
D}$ est \'egalement connexe. Puisque
la partie ${\rm
Int}\left((b''{}_i^{-1})_{\geq}\right)$ contient les points $\sigma_i^{l}$,
$l>0$, ainsi que le point $\sigma_i^*$,
elle rencontre
$A_{i+1}^-\cup A_{i'}^-$. Ceci implique que l'une au moins des briques de
$(b''{}_i^{-1})_{\geq}$ est une brique de
$X$. En utilisant l'assertion {\rm
(ii)} de la \fullref{prop8.2}, on en d\'eduit que le 
futur de $b''{}_i^{-1}$ contient
$b''{}_{i+1}^{-1}$ ou $b''{}_{i'}^{-1}$. En 
conclusion, on sait que pour tout $i\in\Z/n\Z$,
il existe $i'\not=i$ tel que $b''{}_{i'}^{-1}\in (b''{}^{-1}_{i})_{\geq}$.

Il vaut maintenant envisager deux cas:

\bei
\item les briques 
$b''{}_{i'}^{-1}$, $i\in\Z/n\Z$, sont toutes
distinctes;
\item il existe $i_0$ et $i_1$ 
distincts tels que $b''{}_{i_0}^{-1}=
b''{}_{i_1}^{-1}$.
\eni

Dans le premier cas, le fait que pour tout $i\in\Z/n\Z$,
il existe $i'\not=i$ tel que $b''{}_{i'}^{-1}\in (b''{}^{-1}_{i})_{\geq}$
implique qu'il existe $i\in\Z/n\Z$ tel que $b''{}_{i}^{-1}\in
(b''{}^{-1}_{i})_{>}$, la proposition est d\'emontr\'ee.

\begin{figure}[ht!]\small\vspace{2mm}
\labellist
\pinlabel $\alpha_3$ [b] at 89 216
\pinlabel $\alpha_1$ [tl] at 163 15
\pinlabel $\alpha_2$ [lb] at 204 162
\pinlabel $\alpha_4$ [r] at 1 95
\pinlabel $\omega_4$ [l] at 217 95
\pinlabel $\omega_3$ [t] at 89 1
\pinlabel $\omega_2$ [br] at 15 162
\pinlabel $\omega_1$ [bl] at 163 204
\pinlabel $b''$ [b] at 75 97
\endlabellist
\cl{\psfig{file=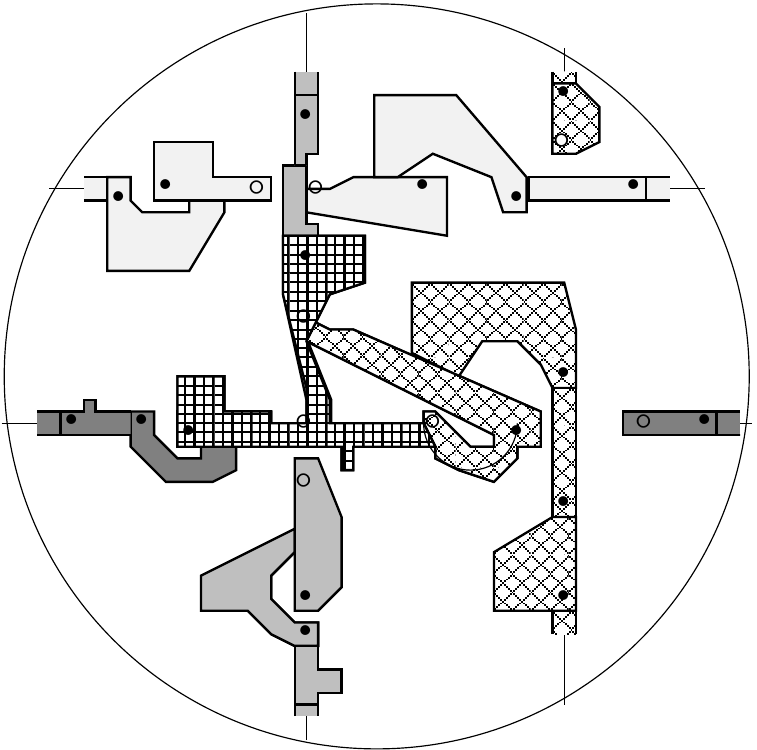,scale=80}}
\caption{}\label{fig15}
\end{figure}

Pla\c cons nous maintenant dans le second cas et 
supposons qu'il existe $i_0$ et $i_1$ distincts 
tels que $b''{}_{i_0}^{-1}=
b''{}_{i_1}^{-1}$.  Nous allons prouver que
$b''\in b''_{>}$, o\`u
$b''=b''{}_{i_0}^{-1}=b''{}_{i_1}^{-1}$.  On peut 
bien s\^ur supposer que $b''$ ne rencontre aucun 
de
ses it\'er\'es stricts et rencontre donc toute 
orbite $O_i$ en au plus un point.
Nous avons illustr\'e cette situation sur la \fullref{fig15}.

 Commen\c cons par prouver que pour tout 
$i\in\Z/n\Z$, il existe $l_i<0$ tel que 
$b''{}_i^l\in
b''_{<}$ si
$l\leq l_i$. Pour cela, notons
$I$ l'ensemble des indices satisfaisant cette 
condition. Nous savons d\'ej\`a que
$b''{}_{i}^l\not=b''$ pour
$l$ petit puisque la singularit\'e
$\sigma_{i}^l$ (qui est dans $b''{}_i^l$) 
appartient alors \`a $O_i$. Remarquons que $I$ 
contient $i_0$
et
$i_1$. En effet, on sait que les briques 
$b''{}_{i_0}^l$ et $b''{}_{i_1}^l$\ appartiennent 
\`a
$(b''{}^{-1})_{\leq}$ et qu'elles sont distinctes 
de $b''{}^{-1}$ si $l$ est petit.
Pour prouver que $I=\Z/n\Z$, il
suffit donc de montrer que
$i-1\in I$ d\`es que
$i\in I$. Posons $i'=i_1$ si $i\not=i_1$ et
$i'=i_2$ si
$i=i_1$. Raisonnons par l'absurde et supposons 
que $i-1\not\in I$. Puisque $b''{}_{i-1}^l\in 
(b''{}_{i-1}^{l'})_{\leq}$, si
$l<0<l'$, on en d\'eduit que $b''_{<}$ ne 
contient aucune brique $b''{}_{i-1}^{l'}$, 
$l'>0$. Puisque
$b''{}_{<}$ est connexe (assertion {\rm(i)} de la 
\fullref{prop8.1}) et contient toutes les briques 
$b''{}_{i}^{l}$ et $b''{}_{i'}^{l}$,
pour $l<0$ assez petit, l'ensemble $b''{}_{<}$ 
s\'epare les points $\sigma_{i-1}^{-l}$ et
$\sigma_{i-1}^{l+1}$, d\`es que $l$ est assez 
grand. Or le futur $(b''{}_{i-1}^{-l})_{\geq}$ est
connexe et contient $b''{}_{i-1}^l$. La partie ${\rm
Int}\left((b''{}_{i-1}^{-l})_{\geq}\right)\subset\D$ est connexe, elle contient
$\sigma_{i-1}^{-l}$ et
$\sigma_{i-1}^{l+1}$, et rencontre donc l'ensemble $b''_{<}\subset{\bf D}$. Les
parties
$(b''{}_{i-1}^{-l})_{\geq}$ et
$b''_{<}$ de $B''$ ont donc au moins une brique en commun, ce qui implique que
$b''{}_{i-1}^{-l}\in b''_{<}$. Ainsi $i-1\in I$, ce qui contredit
l'hypoth\`ese.

Un argument tr\`es proche permet alors de montrer que $b''\in b''{}_{>}$.
L\`a encore, raisonnons par l'absurde et supposons que
$b''{}_{<}$ et $b''{}_{\geq }$ n'ont aucune 
brique en commun. Utilisant de nouveau la
connexit\'e de
$b''{}_{<}$ et le fait que
$I=\Z/n\Z$, on peut remarquer que si $l$ est
grand, $b''{}_{i_0}^{l}$ et $b''{}_{i_1}^{l}$, 
qui sont des briques de $b''{}_{\geq}$ et donc
qui n'appartiennent pas \`a $b''{}_{<}$, sont alors dans
des composantes connexes diff\'erentes de
$B''\setminus (b''{}_{<})$. Ceci est impossible 
puisque  $b''{}_{\geq}$ est connexe, contient
$b''{}_{i_0}^{l}$ et
$b''{}_{i_1}^{l}$ et n'a aucune brique en commun avec $b''{}_{<}$.
\endproof

\bibliographystyle{gtart}
\bibliography{link}

\end{document}